\documentclass[preprint,3p,times]{elsarticle}

\usepackage{example}
\usepackage{tabularray}
\usepackage[dvipsnames]{xcolor}
\usepackage{placeins}
\newcommand{\norm}[1]{\left\lVert #1\right\rVert}
\newcommand{\abs}[1]{\left\lvert #1\right\rvert}
\newcommand{\Var}{\operatorname{Var}}
\newtheorem{proposition}{Proposition}[section]
\def\al{\alpha}

\def\x{\mathbf{x}}

\def\L{\mathcal{L}}
\def\N{\mathcal{N}}

\def\E{\mathbb{E}}

\def\al{\alpha}
\def\dta{\partial_t^\alpha}
\def\d{\mathrm{d}}
\def\btau{\bar\partial_\tau^\al}

\begin{document}

\begin{frontmatter}

\title{Transformed Diffusion-Wave fPINNs: Enhancing Computational Efficiency for Time-Fractional Diffusion-Wave Equations}
\author[1]{Zhengqi Zhang}
\ead{zhengqizhang229@gmail.com}
\author[2]{Jing Li\corref{cor1}}
\ead{lijing@zhejianglab.org}
\cortext[cor1]{Corresponding author}
\affiliation[1]{organization={Ant Group}, city={Hangzhou},country={China}}
\affiliation[2]{organization={Zhejiang Lab}, city={Hangzhou},country={China}}
\begin{abstract}
We propose Transformed Diffusion-Wave Fractional Physics-Informed Neural Networks (tDWfPINNs) as a targeted computational improvement for fractional PINNs applied to time-fractional diffusion-wave equations with fractional order $\alpha \in (1, 2)$. In a direct mesh-free quadrature implementation of the Caputo fractional derivative for $\alpha\in(1,2)$, the derivative of the neural-network solution must be evaluated at many shifted quadrature points, which creates a substantial computational burden for each collocation--quadrature pair. The proposed method avoids these shifted time-derivative evaluations by introducing a specific integrand transformation, so that the quadrature involves neural-network values rather than first-order time derivatives. We provide a comprehensive theoretical analysis and numerical validation of the proposed fractional-derivative scheme, establishing consistency of the transformed operator, endpoint regularity, representation-level stability, Monte Carlo and Gauss--Jacobi quadrature convergence and error estimates, and computing-efficiency gains that scale with the number of collocation points, quadrature points, spatial dimension, and network size. We conduct a comparative numerical study using both Monte Carlo and Gauss--Jacobi quadrature across several one- and two-dimensional time-fractional diffusion-wave PDEs. The experiments show that the transformed Gauss--Jacobi implementation provides a favorable accuracy--cost balance, while the quadrature size and endpoint evaluation require careful treatment. Adaptive sampling is used only as a diagnostic setting to examine the mesh-free compatibility of the proposed formulation. Overall, the results support the proposed transformation as an efficient implementation of time-fractional derivatives within mesh-free fPINNs for the diffusion-wave regime.
\end{abstract}

\end{frontmatter}

\section{Introduction}
\label{sec:intro}

In recent years, fractional and nonlocal models have attracted considerable interest owing to their broad applications in diverse disciplines, including physics, engineering, biology, and finance. Time-fractional diffusion-wave equations \eqref{def:pde}, in particular, have been extensively used to model the propagation of mechanical waves in viscoelastic media \cite{Mainardi1996, Mainardi2022}. Furthermore, such equations are instrumental in modeling the fractional dynamics of systems relevant to medicine \cite{Swain2025, BenSalah2025}. For a comprehensive overview of the applications of fractional models in biology and physics, readers are referred to \cite{Metzler2000, Metzler2014}.

Nevertheless, solving time-fractional diffusion-wave equations presents significant challenges, primarily stemming from their inherent nonlocality and wave-like characteristics. Analytical solutions involving Mittag-Leffler and Wright functions have been introduced in \cite{Sakamoto2011, Jiang2012, Agrawal2002}. To overcome these difficulties, a variety of numerical methods have been developed, including finite-difference methods \cite{Du2010, Sun2020}, convolution-quadrature methods \cite{Cuesta2006, Jin2017Correction}, collocation-type methods \cite{Zhu2019, Stynes2017}, and spectral methods \cite{Chen2020, Zayernouri2013}.

In recent years, the rapid advancement of computing resources and machine learning methodologies has brought Physics-Informed Neural Networks (PINNs) \cite{RaissiParisGE:2019:JCP, CaiMaoWangYinGE:2021:fluid, Cuomo2022:PINN:review} to the forefront of scientific computing. PINNs have attracted considerable attention due to their mesh-free collocation structure, which can be useful for complex geometries and high-dimensional input spaces. To address the computational challenges associated with incorporating fractional derivatives into the PINNs framework, various numerical approximation techniques have been proposed. Among these, fractional Physics-Informed Neural Networks (fPINNs), as proposed in \cite{Pang2019:fPINN}, have been widely adopted. The core strategy of fPINNs involves applying automatic differentiation for integer-order operators, while approximating fractional derivatives using numerical discretization schemes---typically based on finite differences---applied to the neural network outputs. For further developments in this direction, readers are referred to \cite{Lu2025, Vats2024}.

Despite their practical utility, these finite-difference-based schemes for fractional derivatives inherently rely on discretization grids, thereby weakening the mesh-free character that distinguishes PINNs. Recent theoretical investigations \cite{Ren2023, WangYuParis:2022:NTK} have shown that PINNs exhibit a hierarchical spectral bias: they preferentially learn the low-frequency components of the solution, while accurate resolution of high-frequency features is delayed. This spectral bias imposes significant challenges in accurately capturing multiscale dynamics and motivates the development of advanced training strategies, such as adaptive sampling methods \cite{LuLu:2023:RAD, Zhang2025}. In this work, adaptive sampling is also applied to auxiliarily establish the effectiveness of mesh-free methods.

Spectral methods \cite{Zhang2025Spectral, Chen2025} have been proposed to reformulate problems involving fractional derivatives by transforming the original time-space domain into the frequency domain. Although this approach effectively avoids direct numerical evaluation of fractional operators, it typically incurs elevated computational costs and may suffer from reduced accuracy due to the inverse transform. Alternatively, Monte Carlo fractional Physics-Informed Neural Networks (MCfPINNs) \cite{Guo2022:MCfPINN} leverage stochastic sampling to approximate both temporal and spatial fractional derivatives. Building on this framework, recent studies \cite{Hu2024:CMAME, Lin2025} have introduced Gauss--Jacobi quadrature rules to improve computational efficiency by replacing Monte Carlo integration with deterministic quadrature rules. These mesh-free quadrature formulations can evaluate the fractional residual at scattered collocation points and are therefore compatible with adaptive point redistribution.

In this work, we consider an initial-boundary value problem for a diffusion-wave equation with order $\alpha \in (1, 2)$:
\begin{equation}
    \label{def:pde}
    \begin{aligned}
        \dta u(t,\x) + \N(u(t,\x),f(t,\x)) &= 0,~~ (t,\x)\in (0,T]\times\Omega,\\
        u(t,\x)&=0,~~(t,\x)\in(0,T]\times\partial\Omega\\
        u(0,\x)=a(\x),\partial_tu(0,\x)&=b(\x), ~~\x\in\Omega,
    \end{aligned}
\end{equation}
where $T>0$ denotes the final time, and $\Omega\subset\mathbb{R}^d$ ($d=1,2,3$) is a bounded spatial domain with boundary $\partial\Omega$. $\N$ is a time-independent differential operator acting on the solution $u$ and the source term $f$. $\dta u $ denotes the Caputo fractional derivative of order $\al \in (1, 2)$ with respect to time $t$, defined as \cite[P. 70]{Kilbas2006},
\begin{equation*}\dta u(t) = \frac{1}{\Gamma(2-\al)}\int_0^t (t-\tau)^{1-\al}\frac{d^2}{d\tau^2} u(\tau)\d\tau,
\end{equation*}
where $\Gamma$ denotes the Gamma function.

In recent years, various PINN-based methods have been proposed to address diffusion-wave problems. For example, \cite{Vats2024} employed fPINNs to solve time-dependent diffusion-wave equations, while Lu et al. \cite{Lu2025} utilized the $L2$ scheme to approximate Caputo fractional derivatives for $\alpha\in(0,2)$. Chen et al. \cite{Chen2025} applied the Laplace transform to reformulate the problem described by \eqref{def:pde}, and Lin et al. \cite{Lin2025} integrated tensor neural networks with Gauss--Jacobi quadrature to achieve high-accuracy estimation. Nevertheless, there remains a significant gap in the development of efficient mesh-free quadrature implementations for Caputo derivatives of order $\alpha \in (1, 2)$ inside PINN residuals. A natural direct extension of the MCfPINNs framework---originally proposed for $\alpha\in(0,1)$---uses either the Monte Carlo (MC) method \cite{Guo2022:MCfPINN} or Gauss--Jacobi (GJ) quadrature \cite{Hu2024:CMAME}. However, for diffusion-wave equations with $\alpha \in (1, 2)$, this direct quadrature construction contains shifted time derivatives of the neural-network solution at the quadrature points. During PINN training, these shifted derivative evaluations create a substantial automatic-differentiation burden for every collocation--quadrature pair, and the cost grows with the number of collocation points and quadrature points.

To address the aforementioned challenges, we propose transformed Diffusion-Wave fractional Physics-Informed Neural Networks (tDWfPINNs). The key idea is to transform the Caputo residual so that, for quadrature-based fPINNs with $\alpha\in(1,2)$, the quadrature evaluates shifted neural-network values rather than shifted time derivatives. This transformed residual preserves the original Caputo operator while substantially improving the computing efficiency of fractional-derivative evaluation inside the fPINN residual. The resulting advantage becomes more pronounced as the number of collocation points, quadrature points, spatial dimension, and network size increase, and the numerical experiments show that this reduction in computational cost can be obtained without degrading the accuracy level of the PINN solution. At the same time, the practical performance still depends on suitable choices of quadrature size, endpoint treatment, and other training hyperparameters.

We summarize the main contributions of this paper as follows:
\begin{itemize}
    \item \textbf{Transformed Caputo-residual construction for $\alpha \in (1,2)$}:  We propose an integrand transformation for the Caputo derivative in time-fractional diffusion-wave fPINNs. Compared with the direct mesh-free quadrature construction, the transformed formulation avoids shifted time-derivative evaluations of the neural-network solution and thereby improves the computing efficiency of the fractional residual.
    \item \textbf{Theoretical and computational analysis}: We conduct a comprehensive theoretical analysis and numerical validation. We prove equivalence between the transformed formulations and the original Caputo derivative, establish removable endpoint singularities and a representation-level stability bound, and derive Monte Carlo and Gauss--Jacobi quadrature convergence/error estimates. We further validate these estimates through derivative-level numerical tests and verify the predicted computing-efficiency gains.
    \item \textbf{Comprehensive study on a suite of time-fractional PDEs}: We compare the resulting four numerical schemes across one- and two-dimensional time-fractional diffusion-wave problem. The results identify the transformed Gauss--Jacobi implementation as the most favorable accuracy--cost choice in the tested settings, while also showing the advantages of mesh-free qudrature methods.
\end{itemize}

This paper is organized as follows. In Section \ref{sec:setup}, we introduce the framework of Physics-Informed Neural Networks (PINNs) and review the original numerical approaches based on MCfPINNs. Section \ref{sec:methodology} presents the proposed transformed formulation, together with the operator-level analysis, quadrature error estimates, computational-complexity estimates, and a concise summary of derivative-level implementation implications for fractional derivatives of order $\alpha \in (1, 2)$. In Section \ref{sec:experiments}, we evaluate the proposed method through one- and two-dimensional PINN experiments, comparing the quadrature strategies and standard fPINNs baselines in terms of accuracy and computing efficiency. Detailed derivative-level diagnostics and complexity-derivation details are provided in the appendices.

\section{Preliminaries}
\label{sec:setup}
 In this section, we introduce the framework of PINNs related to the time-fractional PDEs \eqref{def:pde} and the original numerical approaches based on the MCfPINNs.

 We begin by considering the abstract formulation of the diffusion-wave equation given in equation~\eqref{def:pde}. Given the initial conditions $a$ and $b$, the source term $f$, and a time-independent differential operator $\N$, our objective is to approximate the solution $u(t,x)$ in $(0,T]\times\Omega$.
\subsection{Physics-Informed Neural Networks (PINNs)}
We briefly review the PINNs framework for solving equation~\eqref{def:pde}, based on the formulation introduced in \cite{RaissiParisGE:2019:JCP}. Let $u(t,x;\theta)$ denote a neural network representation with parameter set $\theta$ as an approximation to the true solution $u(t,x)$. PINNs are trained by minimizing a loss function that enforces agreement with the underlying PDE, initial conditions, and boundary conditions.

\textbf{PDE loss}: Let $\{t_i^{in}, \x_{i}^{in}\}_{i=1}^{N_{in}} \subset (0,T]\times\Omega$ be a set of collocation points in the domain. The PDE loss is defined as
\begin{equation*}\L_{in}(\theta)= \frac{1}{N_{in}}\sum_{i=1}^{N_{in}}\Bigg|\btau u(t_{i}^{in}, \x_{i}^{in};\theta)-\N(u(t_{i}^{in},\x_{i}^{in};\theta),f(t_{i}^{in},\x_{i}^{in}))\Bigg|^2
\end{equation*}
where $\btau u(t_{i}^{in},\x_{i}^{in};\theta)$ denotes a numerical approximation of the Caputo fractional derivative $\dta u(t_{i}^{in},\x_{i}^{in};\theta)$.

\textbf{Boundary Loss}: Let $\{t_i^{bd}, \x_{i}^{bd}\}_{i=1}^{N_{bd}} \subset (0,T]\times\partial\Omega$ be a set of boundary points. The boundary loss is given by
\begin{equation*}\L_{bd}(\theta) = \frac{1}{N_{bd}}\sum_{i=1}^{N_{bd}}\Bigg| u(t_i^{bd},\x_i^{bd};\theta)\Bigg|^2.
\end{equation*}

\textbf{Initial Loss}: Let $\{\x_i^{init}\}_{i=1}^{N_{init}} \subset \Omega$ be a set of points for the initial conditions at $t=0$. The loss of initial conditions is defined as
\begin{equation*}\L_{init}(\theta) = \frac{1}{N_{init}}\sum_{i=1}^{N_{init}}\Bigg|u(0,\x_i^{init};\theta)-a(\x_i^{init})\Bigg|^2+\frac{1}{N_{init}}\sum_{i=1}^{N_{init}}\Bigg|\partial_tu(0,\x_i^{init};\theta)-b(\x_i^{init})\Bigg|^2.
\end{equation*}
Together, these components define the composite loss function used to train the neural network approximation in the PINNs framework.
\textbf{Time-fractional diffusion-wave problem}: $\theta$ is obtained by minimizing the following loss function:
    \begin{equation}
    \label{eqn:loss-forward}
    \theta =\arg\min_{\theta\in\Theta}\L_F(\theta)=\arg\min_{\theta\in\Theta}\L_{in}(\theta)+\lambda_{bd}\L_{bd}(\theta)+\lambda_{init}\L_{init}(\theta),
    \end{equation}
    where the weights $\lambda_{in}$, $\lambda_{bd}$ and $\lambda_{init}$ are specified before training.

\subsection{Improved MCfPINNs for $\al\in(0,1)$}
\label{subsec:IMCPINN}
PINNs have demonstrated significant promise in solving a wide range of PDEs with integer-order derivatives, largely due to the availability of automatic differentiation tools in modern machine learning frameworks such as TensorFlow, PyTorch, and JAX. Their ease of implementation and reliable approximation capabilities have contributed to their growing popularity in scientific computing. However, encoding fractional derivatives directly in PINNs is nontrivial. A fundamental challenge arises from the failure of the classical chain rule in the context of fractional calculus, which impedes the use of automatic differentiation for computing fractional-order derivatives. To address this, a numerical approach based on the substitute form of fractional derivatives for $\alpha\in(0,1)$ has been introduced in \cite{Guo2022:MCfPINN, Hu2024:CMAME}, fully exploiting the tensor-paralleled and mesh-free properties of PINNs.

Firstly, we introduce a lemma that provides an alternative representation of the Caputo fractional derivative. A detailed proof can be found in \cite[Lemma 2.10]{Jin2021:Book}:
\begin{lemma}
\label{lem:represent}
    Let $\al\in(n-1,n)$, $n\in\mathbb{N}$. Given a sufficiently smooth function $f$, we have
    \begin{equation*}
        \begin{aligned}
        \dta f(t) &= \frac{f^{(n-1)}(t)-f^{(n-1)}(0)}{\Gamma(n-\al)t^{\al-n+1}}+\frac{\al-n+1}{\Gamma(n-\al)}\int_0^t \frac{f^{(n-1)}(t)-f^{(n-1)}(\tau)}{(t-\tau)^{(\al-n+2)}}\d\tau\\
        & = \frac{f^{(n-1)}(t)-f^{(n-1)}(0)}{\Gamma(n-\al)t^{\al-n+1}}+\frac{\al-n+1}{\Gamma(n-\al)}t^{n-\al}\int_0^1 \frac{f^{(n-1)}(t)-f^{(n-1)}(t-t\tau)}{t\tau} \tau^{n-1-\al}\d\tau,
        \end{aligned}
    \end{equation*}
    where $f^{(k)}(t)$ denotes the $k$-th derivative of $f$ at $t$.

    For $n=1$, the transformation degenerates to
    \begin{equation}
        \label{eqn:represent_1}
        \begin{aligned}
        \dta f(t) &= \frac{f(t)-f(0)}{\Gamma(1-\al)t^{\al}}+\frac{\al}{\Gamma(1-\al)}t^{1-\al}\int_0^1 \frac{f(t)-f(t-t\tau)}{t\tau} \tau^{-\al}\d\tau\\
        &= \frac{1}{\Gamma(1-\al)}\left\{ \frac{\al}{1-\al}t^{1-\al} \E_{\tau\sim \text{Beta}(1-\al,1)}\left[\frac{f(t)-f(t-t\tau)}{t\tau}\right] +\frac{f(t)-f(0)}{t^\al}\right\}.
        \end{aligned}
    \end{equation}
\end{lemma}

The basic idea of MCfPINNs\cite{Guo2022:MCfPINN} is to approximate the Caputo fractional derivative $\dta f(t)$ using a Monte Carlo approach based on the integral representation \eqref{eqn:represent_1} in Lemma~\ref{lem:represent}. To mitigate numeric blow-ups when some samples are very small, the authors set some hyperparameters, i.e., let $\tau_\varepsilon=\max\{\tau,\varepsilon_t t^{-1}\}$ given fixed $\varepsilon_t$. This yields the approximation:
\begin{equation*}
    \dta f(t) \approx \frac{1}{\Gamma(1-\al)}\left\{\frac{\al}{1-\al}t^{1-\al} \E_{\tau\sim \text{Beta}(1-\al,1)}\left[\frac{f(t)-f(t-t\tau)}{t\tau_\varepsilon}\right] +\frac{f(t)-f(0)}{t^\al}\right\}.
\end{equation*}

However, as the fractional order $\al$ approaches 1, the Beta distribution becomes increasingly concentrated near 0, potentially introducing significant numerical errors due to this treatment.

Therefore, in \cite{Hu2024:CMAME, Lin2025}, the authors propose a new method Improved MCfPINNs using Gauss-Jacobi quadrature. These methods approximate the integral in \eqref{eqn:represent_1} of Lemma~\ref{lem:represent} using a weighted sum over carefully selected quadrature points $(w_i,\tau_i)_{i=1}^N$ (see detailed forms in \cite{GJ2018}), resulting in the approximation:
\begin{equation*}
    \dta f(t) \approx \frac{1}{\Gamma(1-\al)}\left\{\al t^{1-\al} \left(\sum_{i=1}^M w_i \frac{f(t)-f(t-t\tau_i)}{t\tau_i}\right)+\frac{f(t)-f(0)}{t^\al}\right\}
\end{equation*}
The rapid convergence of Gauss-Jacobi quadrature \cite{Chernov2012} significantly accelerates training, as it allows for a smaller number of quadrature points $M$ compared to the Monte Carlo approach when $\alpha<1$. These promising results naturally motivate the extension of this strategy to the case $\alpha>1$, which we explore in this work.

\section{Transformed Diffusion-Wave fPINNs}
\label{sec:methodology}
This section presents the main proposed methodology. Specifically, we first extend the MCfPINNs concept based on Monte Carlo and Gauss-Jacobi from $\alpha\in(0,1)$ to $\alpha \in (1, 2)$ directly. Subsequently, we present a novel transformed approach for approximating the fractional derivative $\partial^\alpha_t f$ with $\alpha \in (1, 2)$, and validate its effectiveness through numerical integration using both Monte Carlo sampling and Gauss-Jacobi quadrature. Monte Carlo integration can be interpreted as a form of numerical quadrature, where both the sampling nodes and their corresponding weights are derived from a probability density function. In this sense, the quadrature points correspond to the sampling nodes used in both Monte Carlo integration and the Gauss-Jacobi method.

For clarity, we denote by $M$ the number of quadrature points used in either the Monte Carlo or Gauss-Jacobi approach, collectively referred to as numerical integration schemes throughout this work.

\subsection{Methodology}
\label{subsec:methodology}
Lemma~\ref{lem:represent} gives an alternative representation of Caputo fractional derivatives of arbitrary orders $\alpha$. For $\alpha\in(1,2)$, the direct use of this representation contains the shifted derivative $f'(t-t\tau)$ inside the quadrature integrand. In a PINN implementation, evaluating $\partial^\alpha_t f(t_i)$ at $N$ collocation points with $M$ quadrature points would therefore require $O(NM)$ automatic-differentiation evaluations at shifted points. The following transformation removes this shifted derivative from the integral while preserving the exact Caputo operator.

\begin{theorem}[Direct and transformed representations]
\label{thm:represent_2}
Let $\alpha\in(1,2)$ and $f\in C^2([0,T])$. Then, for any fixed $t\in(0,T]$,
\begin{align}
\dta f(t)
&=\frac{1}{\Gamma(2-\alpha)}
\left\{
\frac{f'(t)-f'(0)}{t^{\alpha-1}}
+(\alpha-1)t^{2-\alpha}
\int_0^1
\frac{f'(t)-f'(t-t\tau)}{t\tau}\tau^{1-\alpha}\d\tau
\right\} \label{eqn:direct}\\[0.5em]
&=\frac{1}{\Gamma(2-\alpha)}
\Bigg\{
\frac{f'(t)-f'(0)}{t^{\alpha-1}}
-(\alpha-1)\frac{f(t)-f(0)-t f'(t)}{t^\alpha}
\nonumber\\
&\hspace{6em}
-\alpha(\alpha-1)t^{2-\alpha}
\int_0^1
\frac{f(t)-f(t-t\tau)-t\tau f'(t)}{(t\tau)^2}\tau^{1-\alpha}\d\tau
\Bigg\}.\label{eqn:transformed}
\end{align}
\end{theorem}

\begin{proof}
We write the Caputo integral in the memory variable
\[
    r=t-s,\qquad s=t-r,
\]
so that
\begin{equation*}
    \Gamma(2-\alpha)\dta f(t)
    = I(t):=\int_0^t r^{1-\alpha} f''(t-r)\d r .
\end{equation*}

First define
\[
    F_t(r):=f'(t)-f'(t-r),\qquad 0\le r\le t .
\]
Then $F_t(0)=0$, $F_t(t)=f'(t)-f'(0)$, $F_t'(r)=f''(t-r)$, and
\begin{equation}\label{eqn:F-integral-remainder}
    F_t(r)=r\int_0^1 f''(t-\theta r)\d\theta .
\end{equation}
Thus $F_t(r)=O(r)$ as $r\to0^+$. For $0<\varepsilon<t$, integration by parts on $[\varepsilon,t]$ gives
\begin{align*}
I_\varepsilon(t)
&:=\int_\varepsilon^t r^{1-\alpha}F_t'(r)\d r  \\
&=t^{1-\alpha}F_t(t)-\varepsilon^{1-\alpha}F_t(\varepsilon)
+(\alpha-1)\int_\varepsilon^t r^{-\alpha}F_t(r)\d r .
\end{align*}
By \eqref{eqn:F-integral-remainder}, $\varepsilon^{1-\alpha}F_t(\varepsilon)=O(\varepsilon^{2-\alpha})\to0$. Moreover, $r^{-\alpha}F_t(r)=O(r^{1-\alpha})$, which is integrable near zero because $\alpha<2$. Letting $\varepsilon\to0^+$ yields
\begin{equation}\label{eqn:direct-r}
I(t)=t^{1-\alpha}\bigl(f'(t)-f'(0)\bigr)
+(\alpha-1)\int_0^t\bigl(f'(t)-f'(t-r)\bigr)r^{-\alpha}\d r .
\end{equation}
Changing variables $r=t\tau$ proves \eqref{eqn:direct}.

For the transformed representation, define the second-order Taylor remainder around $t$ by
\begin{equation*}
    G_t(r):=f(t)-f(t-r)-r f'(t),\qquad 0\le r\le t .
\end{equation*}
Then $G_t'(r)=f'(t-r)-f'(t)=-F_t(r)$, and Taylor's formula with integral remainder gives
\begin{equation*}
    G_t(r)=-r^2\int_0^1(1-\theta)f''(t-\theta r)\d\theta .
\end{equation*}
Hence $G_t(r)r^{-\alpha}\to0$ and $G_t(r)r^{-\alpha-1}=O(r^{1-\alpha})\in L^1(0,t)$. Let
\[
    J(t):=\int_0^t F_t(r)r^{-\alpha}\d r .
\]
For $0<\varepsilon<t$,
\begin{align*}
J_\varepsilon(t)
&:=\int_\varepsilon^tF_t(r)r^{-\alpha}\d r
=-\int_\varepsilon^tG_t'(r)r^{-\alpha}\d r \\
&=-G_t(t)t^{-\alpha}+G_t(\varepsilon)\varepsilon^{-\alpha}
-\alpha\int_\varepsilon^tG_t(r)r^{-\alpha-1}\d r .
\end{align*}
Letting $\varepsilon\to0^+$ gives
\begin{equation}\label{eqn:J-identity}
J(t)=
-\bigl(f(t)-f(0)-t f'(t)\bigr)t^{-\alpha}
-\alpha\int_0^t\bigl(f(t)-f(t-r)-r f'(t)\bigr)r^{-\alpha-1}\d r .
\end{equation}
Substituting \eqref{eqn:J-identity} into \eqref{eqn:direct-r} gives
\begin{align*}
I(t)
&=t^{1-\alpha}\bigl(f'(t)-f'(0)\bigr)
-(\alpha-1)\bigl(f(t)-f(0)-t f'(t)\bigr)t^{-\alpha}\nonumber\\
&\quad-\alpha(\alpha-1)
\int_0^t\bigl(f(t)-f(t-r)-r f'(t)\bigr)r^{-\alpha-1}\d r .
\end{align*}
Using again $r=t\tau$ in the last integral proves \eqref{eqn:transformed}.
\end{proof}

For later use, introduce the Type-I and Type-II kernels
\begin{align*}
    K_f(t,\tau)
    &:=\frac{f'(t)-f'(t-t\tau)}{t\tau},
    \qquad \tau\in(0,1],\\
    H_f(t,\tau)
    &:=\frac{f(t)-f(t-t\tau)-t\tau f'(t)}{(t\tau)^2},
    \qquad \tau\in(0,1].
\end{align*}
Then the two representations \eqref{eqn:direct} and \eqref{eqn:transformed} in Theorem~\ref{thm:represent_2} can be written as
\begin{align}
\mathcal D_I^\alpha f(t)
&:=\frac{1}{\Gamma(2-\alpha)}
\left\{
\frac{f'(t)-f'(0)}{t^{\alpha-1}}
+(\alpha-1)t^{2-\alpha}\int_0^1K_f(t,\tau)\tau^{1-\alpha}\d\tau
\right\},\label{eqn:D-I}\\
\mathcal D_{II}^\alpha f(t)
&:=\frac{1}{\Gamma(2-\alpha)}
\left\{
\frac{f'(t)-f'(0)}{t^{\alpha-1}}
-(\alpha-1)\frac{f(t)-f(0)-t f'(t)}{t^\alpha}
-\alpha(\alpha-1)t^{2-\alpha}\int_0^1H_f(t,\tau)\tau^{1-\alpha}\d\tau
\right\}.\label{eqn:D-II}
\end{align}
By Theorem~\ref{thm:represent_2}, $\mathcal D_I^\alpha f(t)=\mathcal D_{II}^\alpha f(t)=\dta f(t)$.

The quotient kernels $K_f(t,\tau)$ and $H_f(t,\tau)$ in \eqref{eqn:D-I}--\eqref{eqn:D-II} appear singular as $\tau\to0^+$. Before using these representations for quadrature and stability analysis, we need to verify that the apparent singularities are removable and that the weighted endpoint integrals are well defined. Lemma~\ref{lem:type-kernels} gives the required continuous extensions and uniform bounds.

\begin{lemma}[Removable endpoint singularities]
\label{lem:type-kernels}
Let $f\in C^2([0,T])$ and $t\in(0,T]$. Then $K_f(t,\cdot)$ and $H_f(t,\cdot)$ admit continuous extensions to $[0,1]$ given by
\begin{equation*}
    K_f(t,0)=f''(t),
    \qquad
    H_f(t,0)=-\frac12 f''(t).
\end{equation*}
More precisely,
\begin{align}
    K_f(t,\tau)
    &=\int_0^1 f''(t-\theta t\tau)\d\theta,\label{eqn:K-remainder}\\
    H_f(t,\tau)
    &=-\int_0^1(1-\theta)f''(t-\theta t\tau)\d\theta.\label{eqn:H-remainder}
\end{align}
Consequently,
\begin{equation}\label{eqn:kernel-linf}
    \norm{K_f(t,\cdot)}_{L^\infty(0,1)}\le \norm{f''}_{L^\infty(0,T)},
    \qquad
    \norm{H_f(t,\cdot)}_{L^\infty(0,1)}\le \frac12\norm{f''}_{L^\infty(0,T)}.
\end{equation}
If $f\in C^{m+2}([0,T])$ for some integer $m\ge0$, then $K_f(t,\cdot),H_f(t,\cdot)\in C^m([0,1])$ and
\begin{align}
    \norm{\partial_\tau^m K_f(t,\cdot)}_{L^\infty(0,1)}
    &\le \frac{t^m}{m+1}\norm{f^{(m+2)}}_{L^\infty(0,T)},\label{eqn:K-der-bound}\\
    \norm{\partial_\tau^m H_f(t,\cdot)}_{L^\infty(0,1)}
    &\le \frac{t^m}{(m+1)(m+2)}\norm{f^{(m+2)}}_{L^\infty(0,T)}.\label{eqn:H-der-bound}
\end{align}
\end{lemma}

\begin{proof}
The identities \eqref{eqn:K-remainder} and \eqref{eqn:H-remainder} follow respectively from the first- and second-order Taylor integral remainders. The endpoint values and $L^\infty$ bounds are immediate. Differentiating \eqref{eqn:K-remainder} and \eqref{eqn:H-remainder} $m$ times with respect to $\tau$ gives
\begin{align*}
    \partial_\tau^mK_f(t,\tau)
    &=(-t)^m\int_0^1\theta^m f^{(m+2)}(t-\theta t\tau)\d\theta,\\
    \partial_\tau^mH_f(t,\tau)
    &=(-1)^{m+1}t^m\int_0^1(1-\theta)\theta^m f^{(m+2)}(t-\theta t\tau)\d\theta .
\end{align*}
Using $\int_0^1\theta^m\d\theta=(m+1)^{-1}$ and $\int_0^1(1-\theta)\theta^m\d\theta=((m+1)(m+2))^{-1}$ proves \eqref{eqn:K-der-bound}--\eqref{eqn:H-der-bound}.
\end{proof}

\begin{corollary}[Well-posed weighted integrals]
\label{cor:well-posed}
Under the assumptions of Lemma~\ref{lem:type-kernels}, the weighted integrals
\[
    \int_0^1K_f(t,\tau)\tau^{1-\alpha}\d\tau,
    \qquad
    \int_0^1H_f(t,\tau)\tau^{1-\alpha}\d\tau
\]
are well-defined because $\tau^{1-\alpha}\in L^1(0,1)$ for $\alpha\in(1,2)$. Moreover,
\begin{align*}
    \left|\int_0^1K_f(t,\tau)\tau^{1-\alpha}\d\tau\right|
    &\le \frac{1}{2-\alpha}\norm{f''}_{L^\infty(0,T)},\\
    \left|\int_0^1H_f(t,\tau)\tau^{1-\alpha}\d\tau\right|
    &\le \frac{1}{2(2-\alpha)}\norm{f''}_{L^\infty(0,T)}.
\end{align*}
\end{corollary}

\begin{corollary}[Stability of the transformed representations]
\label{cor:stability}
Let $1<\alpha<2$ and $f,g\in C^2([0,T])$. For every $t\in(0,T]$, the Type-I and Type-II transformed operators satisfy
\begin{equation*}
    \abs{\mathcal D_I^\alpha f(t)-\mathcal D_I^\alpha g(t)}
    =
    \abs{\mathcal D_{II}^\alpha f(t)-\mathcal D_{II}^\alpha g(t)}
    \le
    \frac{t^{2-\alpha}}{\Gamma(3-\alpha)}
    \norm{f''-g''}_{L^\infty(0,T)} .
\end{equation*}
Consequently,
\begin{equation*}
    \sup_{0<t\le T}
    \abs{\mathcal D_I^\alpha f(t)-\mathcal D_I^\alpha g(t)}
    \le
    \frac{T^{2-\alpha}}{\Gamma(3-\alpha)}
    \norm{f''-g''}_{L^\infty(0,T)},
\end{equation*}
and the same bound holds for $\mathcal D_{II}^\alpha$.
\end{corollary}

\begin{proof}
Set $v=f-g$. By Theorem~\ref{thm:represent_2},
\[
    \mathcal D_I^\alpha v(t)=\mathcal D_{II}^\alpha v(t)=\partial_t^\alpha v(t).
\]
Using the Caputo representation,
\[
    \partial_t^\alpha v(t)
    =
    \frac{1}{\Gamma(2-\alpha)}
    \int_0^t(t-s)^{1-\alpha}v''(s)\d s .
\]
Therefore,
\[
    \abs{\partial_t^\alpha v(t)}
    \le
    \frac{\norm{v''}_{L^\infty(0,T)}}{\Gamma(2-\alpha)}
    \int_0^t(t-s)^{1-\alpha}\d s
    =
    \frac{t^{2-\alpha}}{\Gamma(3-\alpha)}
    \norm{v''}_{L^\infty(0,T)} .
\]
This proves the pointwise estimate for both transformed representations, and taking the supremum over $0<t\le T$ gives the uniform bound.
\end{proof}

\begin{remark}
The estimate above is written using the Caputo representation after invoking the equivalence theorem. Therefore, the same stability constant applies to both transformed forms. In particular, the Type-II representation does not introduce an additional continuous-level instability. The apparent factor $(t\tau)^{-2}$ is compensated by the second-order Taylor remainder in the numerator. Numerical cancellation may still occur when the raw quotient is evaluated at extremely small $\tau$, which motivates the endpoint-conditioning analysis for the Monte Carlo cutoff introduced below.
\end{remark}

Based on Theorem~\ref{thm:represent_2}, the four numerical integration schemes used in this work are as follows.
\begin{corollary}[Explicit Monte Carlo and Gauss--Jacobi schemes]
\label{coro:4}

Let $M$ be the number of quadrature points. For the Monte Carlo schemes, let
$\xi_1,\ldots,\xi_M$ be independent samples from $\mathrm{Beta}(2-\alpha,1)$ and set
$\xi_{j,\varepsilon}:=\max\{\xi_j,\varepsilon\}$ for the dimensionless denominator
cutoff; equivalently, $\varepsilon=\varepsilon_t/t$ if a physical cutoff
$\varepsilon_t$ in the memory length is used.
For the Gauss--Jacobi schemes, let $(\tau_i,w_i)_{i=1}^M$ be the $M$-point
Gauss--Jacobi nodes and weights for the weight $\tau^{1-\alpha}$ on $[0,1]$.
The four formulas below define the four numerical schemes used in the
subsequent experiments. In this notation, MC/GJ specifies the quadrature rule,
while I/II specifies whether the direct or transformed representation is used.

\begin{description}
    \item[\textbf{MC-I} (Monte Carlo, Type-I direct form).]
    \begin{align*}
        \mathcal D_{I,M}^{\alpha,\mathrm{MC}} f(t)
        &:=\frac{1}{\Gamma(2-\alpha)}\left\{
        \frac{f'(t)-f'(0)}{t^{\alpha-1}}
        +\frac{\alpha-1}{2-\alpha}t^{2-\alpha}
        \frac{1}{M}\sum_{j=1}^M
        \frac{f'(t)-f'(t-t\xi_j)}{t\xi_{j,\varepsilon}}\right\}.
    \end{align*}

    \item[\textbf{GJ-I} (Gauss--Jacobi, Type-I direct form).]
    \begin{align*}
        \mathcal D_{I,M}^{\alpha,\mathrm{GJ}} f(t)
        &:=\frac{1}{\Gamma(2-\alpha)}\left\{
        \frac{f'(t)-f'(0)}{t^{\alpha-1}}
        +(\alpha-1)t^{2-\alpha}\sum_{i=1}^Mw_i
        \frac{f'(t)-f'(t-t\tau_i)}{t\tau_i}\right\}.
    \end{align*}

    \item[\textbf{MC-II} (Monte Carlo, Type-II transformed form).]
    \begin{align*}
        \mathcal D_{II,M}^{\alpha,\mathrm{MC}} f(t)
        &:=\frac{1}{\Gamma(2-\alpha)}\Bigg\{
        \frac{f'(t)-f'(0)}{t^{\alpha-1}}
        -(\alpha-1)\frac{f(t)-f(0)-tf'(t)}{t^\alpha}\nonumber\\
        &\qquad
        -\frac{\alpha(\alpha-1)}{2-\alpha}t^{2-\alpha}
        \frac{1}{M}\sum_{j=1}^M
        \frac{f(t)-f(t-t\xi_j)-t\xi_j f'(t)}{(t\xi_{j,\varepsilon})^2}\Bigg\}.
    \end{align*}

    \item[\textbf{GJ-II} (Gauss--Jacobi, Type-II transformed form).]
    \begin{align*}
        \mathcal D_{II,M}^{\alpha,\mathrm{GJ}} f(t)
        &:=\frac{1}{\Gamma(2-\alpha)}\Bigg\{
        \frac{f'(t)-f'(0)}{t^{\alpha-1}}
        -(\alpha-1)\frac{f(t)-f(0)-tf'(t)}{t^\alpha}\nonumber\\
        &\qquad
        -\alpha(\alpha-1)t^{2-\alpha}\sum_{i=1}^Mw_i
        \frac{f(t)-f(t-t\tau_i)-t\tau_i f'(t)}{(t\tau_i)^2}\Bigg\}.
    \end{align*}
\end{description}
\end{corollary}

\subsubsection{Monte Carlo error analysis}
We first analyze the sampling error of the Monte Carlo approximations of the weighted integrals in \eqref{eqn:D-I}--\eqref{eqn:D-II}. Owing to the endpoint regularity established above, the transformed kernels can be treated as bounded integrands with respect to the beta density induced by the fractional weight. This leads to unbiased estimators with the standard root-mean-square Monte Carlo rate, while the additional endpoint regularization used in floating-point implementations is analyzed separately below.

Let $\xi$ be a random variable with beta density
\begin{equation*}
    p_\alpha(\tau)=(2-\alpha)\tau^{1-\alpha},\qquad 0<\tau<1.
\end{equation*}
For a bounded function $\phi$ define
\begin{equation*}
    Q_M^{\mathrm{MC}}[\phi]
    :=\frac{1}{(2-\alpha)M}\sum_{j=1}^M\phi(\xi_j),
\end{equation*}
where $\{\xi_j\}_{j=1}^M$ are independent copies of $\xi$. Then
\begin{equation}\label{eqn:QMC-unbiased}
    \mathbb E Q_M^{\mathrm{MC}}[\phi]
    =\int_0^1\phi(\tau)\tau^{1-\alpha}\d\tau .
\end{equation}

\begin{proposition}[Monte Carlo convergence]
\label{prop:mc-convergence}
Let $f\in C^2([0,T])$ and $t\in(0,T]$. Define $\mathcal D_{I,M}^{\alpha,\mathrm{MC}}$ and $\mathcal D_{II,M}^{\alpha,\mathrm{MC}}$ by replacing the weighted integrals in \eqref{eqn:D-I}--\eqref{eqn:D-II} of Corollary~\ref{coro:4} with $Q_M^{\mathrm{MC}}$. Then both estimators are unbiased, and
\begin{align}
\left(\mathbb E\abs{\mathcal D_I^\alpha f(t)-\mathcal D_{I,M}^{\alpha,\mathrm{MC}}f(t)}^2\right)^{1/2}
&\le
\frac{(\alpha-1)t^{2-\alpha}}{\Gamma(3-\alpha)\sqrt M}
\norm{f''}_{L^\infty(0,T)},\label{eqn:MC-I-error}\\
\left(\mathbb E\abs{\mathcal D_{II}^\alpha f(t)-\mathcal D_{II,M}^{\alpha,\mathrm{MC}}f(t)}^2\right)^{1/2}
&\le
\frac{\alpha(\alpha-1)t^{2-\alpha}}{2\Gamma(3-\alpha)\sqrt M}
\norm{f''}_{L^\infty(0,T)}.\label{eqn:MC-II-error}
\end{align}
\end{proposition}

\begin{proof}
Unbiasedness follows from \eqref{eqn:QMC-unbiased} and Theorem~\ref{thm:represent_2}. For any bounded $\phi$,
\[
    \mathbb E\abs{Q_M^{\mathrm{MC}}[\phi]-\int_0^1\phi(\tau)\tau^{1-\alpha}\d\tau}^2
    =\frac{\Var(\phi(\xi))}{(2-\alpha)^2M}
    \le \frac{\norm{\phi}_{L^\infty(0,1)}^2}{(2-\alpha)^2M}.
\]
Applying this estimate with $\phi=K_f(t,\cdot)$ and $\phi=H_f(t,\cdot)$, using \eqref{eqn:kernel-linf} in Lemma~\ref{lem:type-kernels}, and using $\Gamma(3-\alpha)=(2-\alpha)\Gamma(2-\alpha)$ gives \eqref{eqn:MC-I-error}--\eqref{eqn:MC-II-error}.
\end{proof}

\paragraph*{Why a denominator cutoff is introduced}
The cutoff used above has the form
\[
    \tau_\varepsilon=\max\{\tau,\varepsilon_t t^{-1}\},
\]
where $\varepsilon_t$ is a small cutoff in the physical memory length $h=t\tau$. In the normalized variable $\tau$, this is equivalently a dimensionless cutoff
\[
    \delta=\varepsilon_t/t,\qquad
    \tau_\delta=\max\{\tau,\delta\}.
\]
Thus, in the analysis below we write the cutoff in terms of $0<\delta<1$ when the raw Type-I and Type-II quotients are evaluated near the endpoint. The estimates in Proposition~\ref{prop:mc-convergence} are obtained for the exact kernels $K_f$ and $H_f$. In exact arithmetic these kernels are bounded at $\tau=0$ because the apparent singularities are cancelled by the Taylor remainders \eqref{eqn:K-remainder}--\eqref{eqn:H-remainder} in Lemma~\ref{lem:type-kernels}. In a floating-point PINN implementation, however, the kernels are evaluated as raw difference quotients rather than through these integral remainders. Small perturbations in the numerator are amplified by the factors $(t\tau)^{-1}$ in Type-I and $(t\tau)^{-2}$ in Type-II. Thus $\tau_\delta$ serves as an endpoint-conditioning device: it caps the denominator when Monte Carlo samples are too close to zero. The price is a deterministic consistency bias, estimated in Proposition~\ref{prop:regularization-bias}. This issue is particularly relevant for Monte Carlo quadrature because
\begin{equation*}
    \mathbb P(\xi<\delta)=\int_0^\delta(2-\alpha)\tau^{1-\alpha}\d\tau=\delta^{2-\alpha}.
\end{equation*}
When $\alpha\to2$, the exponent $2-\alpha$ becomes small, so a non-negligible fraction of samples may fall into the endpoint region $(0,\delta)$.

\begin{proposition}[Endpoint conditioning under perturbed quotient evaluations]
\label{prop:delta-conditioning}
Let $f\in C^2([0,T])$ and $t\in(0,T]$. Here $\widetilde f(y)$ and
$\widetilde {f'}(y)$ denote the scalar values returned by a numerical
implementation when evaluating $f(y)$ and $f'(y)$, respectively; they are not
assumed to be generated by a single perturbed function. Equivalently, write
\[
    \widetilde f(y)=f(y)+e_0(y),
    \qquad
    \widetilde {f'}(y)=f'(y)+e_1(y),
\]
where the evaluation errors satisfy
\begin{equation*}
    \abs{e_0(y)}=\abs{\widetilde f(y)-f(y)}\le \eta_0,
    \qquad
    \abs{e_1(y)}=\abs{\widetilde {f'}(y)-f'(y)}\le \eta_1,
    \qquad 0\le y\le t .
\end{equation*}
Let $\tau_\delta=\max\{\tau,\delta\}$. Define the exact regularized kernels
\begin{align*}
    K_{f,\delta}(t,\tau)
    &:=\frac{f'(t)-f'(t-t\tau)}{t\tau_\delta},\\
    H_{f,\delta}(t,\tau)
    &:=\frac{f(t)-f(t-t\tau)-t\tau f'(t)}{(t\tau_\delta)^2},
\end{align*}
and their perturbed numerical counterparts
\begin{align*}
    \widetilde K_{f,\delta}(t,\tau)
    &:=\frac{\widetilde {f'}(t)-\widetilde {f'}(t-t\tau)}{t\tau_\delta},\\
    \widetilde H_{f,\delta}(t,\tau)
    &:=\frac{\widetilde f(t)-\widetilde f(t-t\tau)-t\tau\widetilde {f'}(t)}{(t\tau_\delta)^2}.
\end{align*}
Let
\begin{align}
    A_\alpha(\delta)
    &:=\int_0^1\frac{\tau^{1-\alpha}}{\tau_\delta}\d\tau
      =\frac{\delta^{1-\alpha}}{2-\alpha}
       +\frac{\delta^{1-\alpha}-1}{\alpha-1},\label{eqn:A-delta}\\
    B_\alpha(\delta)
    &:=\int_0^1\frac{\tau^{1-\alpha}}{\tau_\delta^2}\d\tau
      =\frac{\delta^{-\alpha}}{2-\alpha}
       +\frac{\delta^{-\alpha}-1}{\alpha},\\
    C_\alpha(\delta)
    &:=\int_0^1\frac{\tau^{2-\alpha}}{\tau_\delta^2}\d\tau
      =\frac{\delta^{1-\alpha}}{3-\alpha}
       +\frac{\delta^{1-\alpha}-1}{\alpha-1}.\label{eqn:C-delta}
\end{align}
Then
\begin{align}
\int_0^1
\abs{\widetilde K_{f,\delta}(t,\tau)-K_{f,\delta}(t,\tau)}\tau^{1-\alpha}\d\tau
&\le \frac{2\eta_1}{t}A_\alpha(\delta),\label{eqn:K-noise-bound}\\
\int_0^1
\abs{\widetilde H_{f,\delta}(t,\tau)-H_{f,\delta}(t,\tau)}\tau^{1-\alpha}\d\tau
&\le \frac{2\eta_0}{t^2}B_\alpha(\delta)+\frac{\eta_1}{t}C_\alpha(\delta).
\label{eqn:H-noise-bound}
\end{align}
\end{proposition}

\begin{proof}
For Type-I, the numerator contains two evaluations of $f'$, so
\[
    \abs{\widetilde K_{f,\delta}(t,\tau)-K_{f,\delta}(t,\tau)}
    \le \frac{2\eta_1}{t\tau_\delta}.
\]
Multiplying by $\tau^{1-\alpha}$ and integrating gives \eqref{eqn:K-noise-bound}. For Type-II,
\[
\abs{\widetilde H_{f,\delta}(t,\tau)-H_{f,\delta}(t,\tau)}
\le\frac{2\eta_0+t\tau\eta_1}{(t\tau_\delta)^2}
=\frac{2\eta_0}{t^2\tau_\delta^2}+\frac{\eta_1\tau}{t\tau_\delta^2},
\]
which gives \eqref{eqn:H-noise-bound}. The closed forms \eqref{eqn:A-delta}--\eqref{eqn:C-delta} are obtained by splitting the integrals over $(0,\delta)$ and $(\delta,1)$.
\end{proof}

\begin{proposition}[Bias induced by denominator regularization]
\label{prop:regularization-bias}
Let $f\in C^2([0,T])$, $t\in(0,T]$, and $0<\delta<1$. Define
\begin{align*}
    K_{f,\delta}(t,\tau)
    &:=\frac{f'(t)-f'(t-t\tau)}{t\tau_\delta},\\
    H_{f,\delta}(t,\tau)
    &:=\frac{f(t)-f(t-t\tau)-t\tau f'(t)}{(t\tau_\delta)^2}.
\end{align*}
Let $\mathcal D_{I,\delta}^{\alpha}$ and $\mathcal D_{II,\delta}^{\alpha}$ denote the Type-I and Type-II operators obtained by replacing $K_f,H_f$ with $K_{f,\delta},H_{f,\delta}$ in \eqref{eqn:D-I} and \eqref{eqn:D-II} of Corollary~\ref{coro:4}. Then
\begin{align}
\abs{\mathcal D_I^\alpha f(t)-\mathcal D_{I,\delta}^{\alpha}f(t)}
&\le
\frac{(\alpha-1)t^{2-\alpha}}{\Gamma(3-\alpha)}
\norm{f''}_{L^\infty(0,T)}\delta^{2-\alpha},\label{eqn:reg-bias-I}\\
\abs{\mathcal D_{II}^\alpha f(t)-\mathcal D_{II,\delta}^{\alpha}f(t)}
&\le
\frac{\alpha(\alpha-1)t^{2-\alpha}}{2\Gamma(3-\alpha)}
\norm{f''}_{L^\infty(0,T)}\delta^{2-\alpha}.\label{eqn:reg-bias-II}
\end{align}
Equivalently, if the physical cutoff is $\varepsilon$ in the memory length $r=t\tau$, then $\delta=\varepsilon/t$ and the bias is $O(\varepsilon^{2-\alpha})$.
\end{proposition}

\begin{proof}
The regularized and unregularized kernels differ only on $0<\tau<\delta$. For Type-I, the remainder formula \eqref{eqn:K-remainder} in Lemma~\ref{lem:type-kernels} implies $K_{f,\delta}(t,\tau)=(\tau/\delta)K_f(t,\tau)$ on this interval, and hence
\[
    \int_0^1\abs{K_f-K_{f,\delta}}\tau^{1-\alpha}\d\tau
    \le\norm{f''}_{L^\infty(0,T)}\int_0^\delta\tau^{1-\alpha}\d\tau
    =\frac{\delta^{2-\alpha}}{2-\alpha}\norm{f''}_{L^\infty(0,T)}.
\]
Multiplying by the Type-I prefactor in \eqref{eqn:D-I} of Corollary~\ref{coro:4} gives \eqref{eqn:reg-bias-I}. For Type-II, the remainder formula \eqref{eqn:H-remainder} in Lemma~\ref{lem:type-kernels} gives $H_{f,\delta}(t,\tau)=(\tau/\delta)^2H_f(t,\tau)$ on $(0,\delta)$, so
\[
    \int_0^1\abs{H_f-H_{f,\delta}}\tau^{1-\alpha}\d\tau
    \le\frac{1}{2}\norm{f''}_{L^\infty(0,T)}\frac{\delta^{2-\alpha}}{2-\alpha}.
\]
Multiplying by the Type-II prefactor in \eqref{eqn:D-II} of Corollary~\ref{coro:4} proves \eqref{eqn:reg-bias-II}.
\end{proof}

\begin{remark}[Bias--conditioning trade-off]
\label{rem:delta-tradeoff}
The cutoff $\delta$ should not be interpreted as part of the exact mathematical transformation. It is a numerical stabilization parameter. Proposition~\ref{prop:delta-conditioning} shows that endpoint-induced perturbations grow as $\delta\to0^+$; in particular $A_\alpha(\delta),C_\alpha(\delta)=O(\delta^{1-\alpha})$ and $B_\alpha(\delta)=O(\delta^{-\alpha})$. Proposition~\ref{prop:regularization-bias} shows the opposite effect: the deterministic consistency bias produced by replacing $\tau$ with $\tau_\delta$ is $O(\delta^{2-\alpha})$. Thus $\delta$ balances two competing errors:
\begin{equation*}
    \text{total endpoint error}
    \approx
    \underbrace{O(\delta^{2-\alpha})}_{\text{regularization bias}}
    +
    \underbrace{O(\eta_1\delta^{1-\alpha})}_{\text{Type-I derivative-evaluation perturbation}}
\end{equation*}
for Type-I, and
\begin{equation*}
    \text{total endpoint error}
    \approx
    \underbrace{O(\delta^{2-\alpha})}_{\text{regularization bias}}
    +
    \underbrace{O(\eta_0\delta^{-\alpha})+O(\eta_1\delta^{1-\alpha})}_{\text{Type-II raw-quotient perturbation}}
\end{equation*}
for Type-II, up to powers of $t$ and constants depending on $\alpha$. This explains why the cutoff is useful in Monte Carlo implementations and why the Monte Carlo schemes become increasingly sensitive as $\alpha\to2$.
\end{remark}

\subsubsection{Gauss-Jacobi quadrature error analysis}
We next turn to deterministic Gauss--Jacobi quadrature, which is tailored to the singular weight $\tau^{1-\alpha}$ in the transformed Caputo representations. The quadrature error is therefore governed by the regularity of the endpoint-extended kernels $K_f(t,\cdot)$ and $H_f(t,\cdot)$: finite smoothness yields algebraic convergence in $M$, while analytic kernels yield substantially faster convergence. This analysis explains both the efficiency of Gauss--Jacobi quadrature and its sensitivity to kernel smoothness.

Let $Q_M^{\mathrm{GJ}}$ denote the $M$-point Gauss-Jacobi quadrature rule on $[0,1]$ for the weight $\tau^{1-\alpha}$:
\begin{equation*}
    Q_M^{\mathrm{GJ}}[\phi]:=\sum_{j=1}^M w_j\phi(\tau_j)
    \approx \int_0^1\phi(\tau)\tau^{1-\alpha}\d\tau .
\end{equation*}
The rule is exact for polynomials of degree at most $2M-1$ with respect to the weighted measure $\tau^{1-\alpha}\d\tau$.

\begin{theorem}[Gauss-Jacobi error estimates]
\label{thm:GJ-convergence}
Let $t\in(0,T]$.

\smallskip
\noindent\textup{(i) Finite smoothness.} Suppose $K_f(t,\cdot),H_f(t,\cdot)\in C^r([0,1])$ for an integer $r\ge1$. Then there exists a constant $C_{\alpha,r}>0$, independent of $M$, $t$, and $f$, such that
\begin{align}
\abs{\mathcal D_I^\alpha f(t)-\mathcal D_{I,M}^{\alpha,\mathrm{GJ}}f(t)}
&\le C_{\alpha,r}t^{2-\alpha}M^{-r}\norm{\partial_\tau^r K_f(t,\cdot)}_{L^\infty(0,1)},\label{eqn:GJ-I-algebraic}\\
\abs{\mathcal D_{II}^\alpha f(t)-\mathcal D_{II,M}^{\alpha,\mathrm{GJ}}f(t)}
&\le C_{\alpha,r}t^{2-\alpha}M^{-r}\norm{\partial_\tau^r H_f(t,\cdot)}_{L^\infty(0,1)}.\label{eqn:GJ-II-algebraic}
\end{align}

\smallskip
\noindent\textup{(ii) Analytic case.} Let $\mathcal B_\rho$ be the Bernstein ellipse with parameter $\rho>1$ on $[0,1]$ in the sense of \cite{Tadmor1986}. If $K_f(t,\cdot)$ and $H_f(t,\cdot)$ admit bounded analytic continuations to $\mathcal B_\rho$, set
\[
    B_K:=\sup_{z\in\mathcal B_\rho}\abs{K_f(t,z)},
    \qquad
    B_H:=\sup_{z\in\mathcal B_\rho}\abs{H_f(t,z)}.
\]
Then there exists a constant $C_{\alpha,\rho}>0$, independent of $M$, $t$, and $f$, such that
\begin{align}
\abs{\mathcal D_I^\alpha f(t)-\mathcal D_{I,M}^{\alpha,\mathrm{GJ}}f(t)}
&\le C_{\alpha,\rho}t^{2-\alpha}B_K\rho^{-2M},\label{eqn:GJ-I-spectral}\\
\abs{\mathcal D_{II}^\alpha f(t)-\mathcal D_{II,M}^{\alpha,\mathrm{GJ}}f(t)}
&\le C_{\alpha,\rho}t^{2-\alpha}B_H\rho^{-2M}.\label{eqn:GJ-II-spectral}
\end{align}
\end{theorem}

\begin{proof}
Define the quadrature error functional
\[
    E_M^{\mathrm{GJ}}[\phi]
    :=\int_0^1\phi(\tau)\tau^{1-\alpha}\d\tau-Q_M^{\mathrm{GJ}}[\phi].
\]
Then
\begin{align}
\abs{\mathcal D_I^\alpha f(t)-\mathcal D_{I,M}^{\alpha,\mathrm{GJ}}f(t)}
&=\frac{(\alpha-1)t^{2-\alpha}}{\Gamma(2-\alpha)}
\abs{E_M^{\mathrm{GJ}}[K_f(t,\cdot)]},\label{eqn:GJ-I-exact-error}\\
\abs{\mathcal D_{II}^\alpha f(t)-\mathcal D_{II,M}^{\alpha,\mathrm{GJ}}f(t)}
&=\frac{\alpha(\alpha-1)t^{2-\alpha}}{\Gamma(2-\alpha)}
\abs{E_M^{\mathrm{GJ}}[H_f(t,\cdot)]}.\label{eqn:GJ-II-exact-error}
\end{align}
Let $\mu_0:=\int_0^1\tau^{1-\alpha}\d\tau=1/(2-\alpha)$. The Gauss-Jacobi weights are positive and satisfy $\sum_{j=1}^Mw_j=\mu_0$. If $p$ is any polynomial of degree at most $2M-1$, exactness gives $E_M^{\mathrm{GJ}}[p]=0$, hence
\begin{align}
    \abs{E_M^{\mathrm{GJ}}[\phi]}
    &=\abs{E_M^{\mathrm{GJ}}[\phi-p]}\nonumber\\
    &\le\int_0^1\abs{\phi-p}\tau^{1-\alpha}\d\tau
       +\sum_{j=1}^Mw_j\abs{\phi(\tau_j)-p(\tau_j)}\nonumber\\
    &\le2\mu_0\norm{\phi-p}_{L^\infty(0,1)}.
    \label{eqn:GJ-best-approx}
\end{align}
Taking the infimum over all such $p$ reduces the quadrature estimate to best polynomial approximation.

For finite smoothness, the Jackson's inequality
\[
    \inf_{\deg p\le 2M-1}\norm{\phi-p}_{L^\infty(0,1)}
    \le C_rM^{-r}\norm{\phi^{(r)}}_{L^\infty(0,1)}
\]
holds for $\phi\in C^r([0,1])$, with $C_r$ depending only on $r$. Applying this with $\phi=K_f(t,\cdot)$ and $\phi=H_f(t,\cdot)$ gives
\[
    \abs{E_M^{\mathrm{GJ}}[K_f(t,\cdot)]}
    \le C_{\alpha,r}M^{-r}\norm{\partial_\tau^r K_f(t,\cdot)}_{L^\infty(0,1)},
\]
and the same bound, with $K_f$ replaced by $H_f$. Multiplying by the prefactors in \eqref{eqn:GJ-I-exact-error}--\eqref{eqn:GJ-II-exact-error} proves \eqref{eqn:GJ-I-algebraic}--\eqref{eqn:GJ-II-algebraic}.

For the analytic case, let $\mathcal B_\rho$ be the Bernstein ellipse on $[0,1]$ described in the theorem, and assume $\abs{\phi}\le B_\phi$ on $\mathcal B_\rho$. The analyticity assumption implies geometric decay of the associated Chebyshev coefficients, which is proved in \cite{Tadmor1986}. Consequently, the standard polynomial-approximation estimate gives
\[
    \inf_{\deg p\le n}\norm{\phi-p}_{L^\infty(0,1)}
    \le \frac{2B_\phi}{\rho-1}\rho^{-n}.
\]
Taking $n=2M-1$ in \eqref{eqn:GJ-best-approx} yields
\[
    \abs{E_M^{\mathrm{GJ}}[\phi]}
    \le \frac{4\mu_0B_\phi}{\rho-1}\rho^{-(2M-1)}
    \le C_{\alpha,\rho}B_\phi\rho^{-2M},
\]
where the harmless extra factor $\rho$ is absorbed into the constant. Applying this estimate to $K_f(t,\cdot)$ with $B_\phi=B_K$ and to $H_f(t,\cdot)$ with $B_\phi=B_H$, then multiplying by \eqref{eqn:GJ-I-exact-error}--\eqref{eqn:GJ-II-exact-error}, proves \eqref{eqn:GJ-I-spectral}--\eqref{eqn:GJ-II-spectral}.
\end{proof}

\subsubsection{Computational complexity}
\label{subsubsec:computational-complexity}

We now make the computational saving explicit at the level of a fully connected
network evaluation. Consider one input point
$z=(t,x_1,\ldots,x_n)\in\mathbb R^d$, where $d=n+1$. Assume that
$u(z;\theta)$ is represented by a bias-free fully connected network with $L$ hidden
layers, hidden width $H$, scalar output, and a smooth activation function. We
count dense-layer multiply--accumulate operations (MACs) and activation
elementwise operations; there are no bias additions or bias-gradient reductions
in this counting model. One plain forward evaluation contains
\begin{equation}\label{eqn:forward-macs}
    A_{\rm mac}
    =Hd+(L-1)H^2+H
    =H(n+1)+(L-1)H^2+H
\end{equation}
MACs.

\paragraph{Graph-storage count}
The per-layer fp64 storage derivation is deferred to
\ref{app:gpu-storage}. At the tensor-data level, one ordinary value
graph and one time-derivative graph require
\begin{equation}\label{eqn:value-graph-storage-main}
    S_u(1)=8\bigl(d+c_\sigma LH+1\bigr)
    \quad\text{bytes},
\end{equation}
and
\begin{equation}\label{eqn:dt-graph-storage-main}
    S_{\partial_tu}(1)
    =8\bigl[2d+(c_\sigma+2)LH+1\bigr]
    \quad\text{bytes},
\end{equation}
where $c_\sigma\in[1,2]$ accounts for whether only activations or both
activations and pre-activations are retained. The larger constant in
$S_{\partial_tu}(1)$ comes from retaining the derivative pass so that loss
backpropagation can compute $\nabla_\theta \partial_t u(z;\theta)$. Thus one
shifted quadrature node changed from a Type-I time-derivative graph to a Type-II
ordinary value graph saves
\begin{equation}\label{eqn:dt-value-storage-gap-main}
    S_{\partial_tu}(1)-S_u(1)=8(d+2LH)
    \quad\text{bytes}.
\end{equation}

Type-I has $M+2$ time-derivative graphs per collocation point, whereas Type-II
has $M$ ordinary shifted value graphs and two endpoint time-derivative graphs.
Consequently,
\begin{align}
    S_I^{\rm graph}
    &=N(M+2)S_{\partial_tu}(1) \notag\\
    &=8N(M+2)\bigl[2d+(c_\sigma+2)LH+1\bigr],
    \label{eqn:typeI-graph-storage-main}\\
    S_{II}^{\rm graph}
    &=NM S_u(1)+2N S_{\partial_tu}(1) \notag\\
    &=8NM(d+c_\sigma LH+1)
      +16N\bigl[2d+(c_\sigma+2)LH+1\bigr].
    \label{eqn:typeII-graph-storage-main}
\end{align}
Subtracting gives the graph-storage saving
\begin{equation}\label{eqn:graph-storage-saving-main}
    S_I^{\rm graph}-S_{II}^{\rm graph}
    =NM\{S_{\partial_tu}(1)-S_u(1)\}
    =8NM(d+2LH)
    \quad\text{bytes}.
\end{equation}

\paragraph{Training-backward FLOPs}
The one-graph backward counts used below are derived in
\ref{app:gpu-flops}. With the same dense-layer model and a $\tanh$
activation, replacing one shifted time-derivative graph by one ordinary value
graph saves $4A_{\rm mac}+5LH$ FLOPs in the training backward pass. Therefore
\begin{align}
    \Delta F^{\rm bwd}
    &:=F_I^{\rm bwd}-F_{II}^{\rm bwd} \notag\\
    &=NM\bigl(4A_{\rm mac}+5LH\bigr) \notag\\
    &=4NMHd+4NM(L-1)H^2+4NMH+5NMLH .
    \label{eqn:bwd-flops-saving-main}
\end{align}

This explains why the transformed Type-II formulation reduces both GPU memory
and training-backward computation: it replaces the $NM$ shifted
higher-order time-derivative graphs of Type-I by $NM$ ordinary value graphs.

The analysis above shows that Type-I and Type-II are exactly consistent with
the Caputo derivative, have the same representation-level stability bound, and
admit convergent Monte Carlo and Gauss-Jacobi quadrature approximations. The
endpoint bias--conditioning trade-off is induced by denominator regularization.
The computational advantage of Type-II is that it replaces shifted
mixed-derivative backward graphs by ordinary value-backward graphs, leading to
lower GPU memory and lower training-backward FLOPs.

\subsection{Derivative-level validation and implementation diagnostics}
\label{subsec:num-valid}

Detailed derivative-level validation and implementation diagnostics are reported in \ref{app:num-valid}. These tests complement the operator and quadrature analysis above by isolating Monte Carlo sampling error, Gauss--Jacobi quadrature error, endpoint regularization bias, finite-precision cancellation, and empirical graph-storage/training-backward FLOP savings. They provide the numerical basis for choosing quadrature rules, quadrature sizes, endpoint treatment, and network-related hyperparameters used in Section~\ref{sec:experiments}.

The diagnostics show that Monte Carlo estimators follow the expected root-mean-square behavior but become sensitive to endpoint cutoff as $\alpha\to2$, whereas Gauss--Jacobi quadrature is substantially more efficient for smooth transformed kernels and degrades when kernel regularity is lost. The empirical complexity tests verify the predicted Type-II savings with respect to the number of collocation points, quadrature points, input dimension, and network size.

\section{Experiments}
\label{sec:experiments}

In this section, we use the transformed fractional-derivative formulation developed in Section~\ref{sec:methodology} to solve time-fractional partial differential equations with $\alpha \in (1, 2)$. In the following, we use \textbf{MC-I}, \textbf{MC-II}, \textbf{GJ-I}, and \textbf{GJ-II} to refer to the four numerical integration schemes in Corollary~\ref{coro:4}. The experiments are designed to assess the accuracy--cost behavior of the transformed Caputo-residual construction inside PINN training. We keep the architecture, loss formulation, and optimizer choices fixed in order to focus the comparison on whether the transformed residual preserves a comparable error level while reducing the computational cost of fractional-derivative evaluation.

When adaptive resampling is used, we adopt established resampling strategies such as \cite{Zhang2025,LuLu:2023:RAD}. They are used to examine whether the proposed mesh-free fractional residual can be evaluated naturally at adaptively redistributed collocation points.

Unless otherwise specified, we adopt the following settings in all numerical experiments. The neural network architecture consists of 7 hidden layers, each with 20 neurons, and employs the hyperbolic tangent \textit{tanh} as the activation function. This choice ensures that Theorem~\ref{thm:represent_2} and Corollary~\ref{coro:4} can be applied directly to the smooth neural network solution $u(t,\x;\theta)$. For the Monte Carlo methods, $\varepsilon=10^{-7}$ for \textbf{MC-I} and \textbf{MC-II}. To measure the accuracy of the learned solutions, we use the $L^2$ relative error defined as
\begin{equation*}
    e_r(\theta) = \frac{\sqrt{\sum_{i=1}^{N} [u(t_i,\x_i;\theta)-u^*(t_i,\x_i)]^2}}{\sqrt{\sum_{i=1}^{N} [u^*(t_i,\x_i)]^2}},
\end{equation*}
where $(t_i,\x_i)$, $i=1,2,...,N$ represent the test data set. All training and testing are conducted on an NVIDIA V100 (32 GB).

\subsection{Summary of results and findings}
The derivative-level summary in Section~\ref{subsec:num-valid}, together with the detailed diagnostics in Appendix~\ref{app:num-valid}, provides the basis for choosing quadrature rules, quadrature sizes, endpoint treatment, and network-related hyperparameters. Building on this guidance, the PDE experiments below test, across one- and two-dimensional diffusion-wave problems, whether the operator-level and implementation-level advantages established in Section~\ref{sec:methodology} translate into practical PINN training.

Regarding computational efficiency, the transformed schemes reduce the time and GPU resource usage relative to their direct Type-I counterparts while preserving a comparable level of accuracy in the tested PINN settings. This efficiency claim is therefore made relative to the direct quadrature-based residual construction inside fPINNs, not as a universal comparison with all fractional-PDE solvers.

Furthermore, the Gauss--Jacobi method typically requires significantly fewer quadrature points than the Monte Carlo method to achieve similar accuracy (as demonstrated in Section~\ref{subsec:eg1}). Together with the transformed residual, this explains why \textbf{GJ-II} gives the most favorable accuracy--cost profile in our experiments: Gauss--Jacobi quadrature reduces the number of quadrature nodes needed for smooth transformed kernels, while the transformed formulation avoids shifted time-derivative evaluations at those nodes. However, the wave component introduced by the Mittag-Leffler function can strongly influence the numerical behavior of both Gauss--Jacobi quadrature and Monte Carlo sampling. The practical accuracy therefore depends on the quadrature rule, the number of quadrature points, and stable endpoint evaluation.

\begin{remark}
    We briefly address the stability characteristics of the proposed methods. Our numerical experiments suggest two key observations:
\begin{enumerate}
    \item \textbf{Integration Method:} The Gauss--Jacobi approach consistently exhibits superior stability compared to Monte Carlo integration. This is attributed to the fact that Gauss--Jacobi quadratures are specifically designed to handle the singular weights inherent in fractional operators, whereas Monte Carlo sampling can suffer from high variance near singularities.
    \item \textbf{Type I vs. Type II:} The transformed formulation (Type II) reduces computational cost by avoiding automatic differentiation at shifted quadrature points. The apparent stronger singularity in \eqref{eqn:transformed} of Theorem~\ref{thm:represent_2} is removable at the continuous level, and Corollary~\ref{cor:stability} shows that Type-I and Type-II have the same representation-level stability bound. The remaining practical sensitivity is therefore a numerical conditioning issue of raw difference quotients and quadrature near $\tau=0$, which is quantified by the regularization-bias estimate in Proposition~\ref{prop:regularization-bias}.
\end{enumerate}

\end{remark}

We conduct numerical experiments on time-fractional Burgers equations with $\alpha \in (1, 2)$ and incorporate the RAD method from \cite{LuLu:2023:RAD}. The purpose of using RAD here is to test the compatibility between the proposed mesh-free residual and adaptive collocation. The Burgers results show that RAD does not always provide a clear improvement for diffusion-wave dynamics; even when adaptive samples identify important regions, the nonlocal memory and wave-like components can still dominate the training difficulty. This interaction between adaptive sampling and time-fractional wave-like dynamics deserves further study.

We further include a fractional ODE diagnostic to compare the mesh-free \textbf{GJ-II} formulation with a standard fPINNs baseline under adaptive sampling. Under the same 100-point setting, \textbf{GJ-II}+RAD evaluates the residual directly at adaptively selected collocation points, whereas standard fPINNs remains tied to a fixed temporal stencil. This test isolates one practical mesh-free advantage of the proposed residual construction.

Finally we include a two-dimensional experiment on an L-shaped domain. This test uses a non-rectangular spatial geometry and numerical reference data, and therefore checks whether the transformed fractional-derivative estimators remain applicable beyond one-dimensional analytic benchmarks and rectangular domains. It also compares \textbf{MC-I}, \textbf{MC-II}, \textbf{GJ-I}, \textbf{GJ-II}, and standard fPINNs baselines under controlled point-budget settings, placing the mesh-free quadrature formulation in direct comparison with an existing fPINNs practice.

\subsection{Time-fractional diffusion-wave equations in 1D domain}
\label{subsec:eg1}
In the warm-up stage, we consider the following time-fractional PDEs in $\Omega=(0,1)$ with $\al=1.5$:
\begin{equation}
    \label{eqn:eg1}
    \begin{aligned}
        \dta u(t,x) - u_{xx}(t,x) &= 0,~~ (t,\x)\in (0,1]\times\Omega,\\
        u(t,0)=u(t,1)&=0,~~t\in(0,1],\\
        u(0,x)&=2\sin(\pi x), ~~x\in\Omega,\\
        \partial_tu(0,x)&=-\sin(\pi x), ~~x\in\Omega,
    \end{aligned}
\end{equation}
with the exact solution
\begin{equation*}
     u(t,x)= (2E_{\al,1}(-\pi^2 t^\al) - t E_{\al,2}(-\pi^2 t^\al))\sin(\pi x)
\end{equation*}
here $E_{\al,\beta}(t)$ is the Mittag-Leffler function (see details in \cite{Podlubny1999,Jin2021:Book}).

\begin{figure}[htbp]
    \centering
    \includegraphics[scale=0.35]{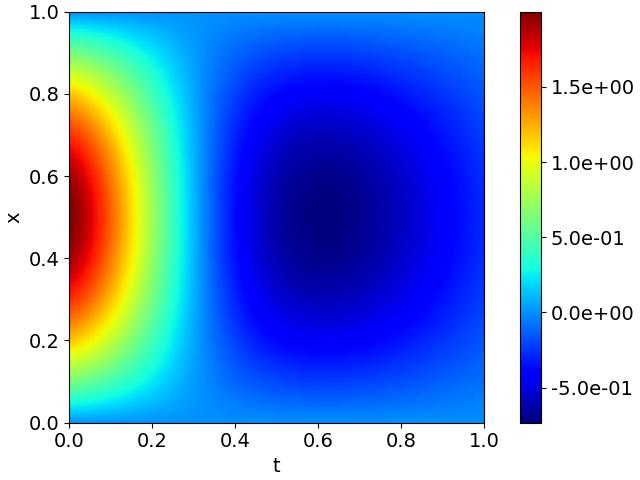}
    \caption{Solution profiles of $u(t,x)$ in equation \eqref{eqn:eg1}.}
    \label{fig:eg1-exact}
\end{figure}

The exact solution is depicted in Figure \ref{fig:eg1-exact}. In contrast to (sub)diffusion equations, the solution of the diffusion-wave equation exhibits a distinctive wave-like structure, primarily driven by the Mittag-Leffler functions. This characteristic introduces additional complexity, making such equations potentially more difficult to solve using vanilla PINNs compared to (sub)diffusion equations. To assess the computational performance of the numerical integration methods presented in Corollary~\ref{coro:4}, we conducted extensive tests with results summarized in Table \ref{tab:forward}.

For small batches, we sampled 5000 collocation points within the domain, along with 1000 points each on the spatial and temporal boundaries. In the large-batch setting, we increased the sampling density to 10000 points in the domain and 2000 points on both the spatial and temporal boundaries. We set the maximum number of iterations to be $L=10$, with each iteration comprising 5000 epochs trained using the Adam optimizer at a learning rate of $10^{-3}$. The weighting parameter in the loss function \eqref{eqn:loss-forward} was set to 1.0.

Regarding computational efficiency, our analysis reveals that the transformed schemes exhibit significantly lower and more stable computation times compared to their direct counterparts (\textbf{MC-I} and \textbf{GJ-I}), which demonstrate rapidly increasing computational costs. For each value of $M$ (number of quadrature points), the transformed schemes consistently require substantially less computational time than direct methods; even with $M=16$ for the Gauss-Jacobi method, our proposed method could save several seconds over 5000 epochs. Additionally, the memory usage is markedly reduced. For instance, with small batches at $M=1280$, \textbf{MC-II} consumes only 10.5GB of GPU memory, whereas \textbf{MC-I} requires 29.1GB. This reduced memory footprint enables the transformed MCfPINNs to accommodate denser sampling points and larger $M$ within a single GPU environment. Importantly, these gains in efficiency come without sacrificing numerical accuracy; both integration types maintain comparable error levels.

Moreover, the Gauss-Jacobi quadrature method appears to offer computational advantages over the Monte Carlo approach, as it requires fewer quadrature points to achieve comparable levels of accuracy, further reinforcing its suitability for efficient and precise fractional-order computations.

\begin{table}[htbp]

    \centering
    \begin{tabular}{|c|cc|cc|c|cc|cc|}
        \hline
        Cases    & \multicolumn{2}{c|}{Small batches}                                                                                                        & \multicolumn{2}{c|}{Large batches}                                                                                                       &          & \multicolumn{2}{c|}{Small batches}                                                                                                     & \multicolumn{2}{c|}{Large batches}                                                                                                     \\ \hline
        $M$/Type & \multicolumn{1}{c|}{\textbf{MC-I}}                                             & \textbf{MC-II}                                           & \multicolumn{1}{c|}{\textbf{MC-I}}                                            & \textbf{MC-II}                                           & $M$/Type & \multicolumn{1}{c|}{\textbf{GJ-I}}                                           & \textbf{GJ-II}                                          & \multicolumn{1}{c|}{\textbf{GJ-I}}                                           & \textbf{GJ-II}                                          \\ \hline
        $M=80$   & \multicolumn{1}{c|}{\begin{tabular}[c]{@{}c@{}}1.91e-1\\ 5m26s\end{tabular}}   & \begin{tabular}[c]{@{}c@{}}8.87e-2\\ 4m03s\end{tabular}  & \multicolumn{1}{c|}{\begin{tabular}[c]{@{}c@{}}1.12e-1\\ 8m55s\end{tabular}}  & \begin{tabular}[c]{@{}c@{}}6.16e-2\\ 5m22s\end{tabular}  & $M=16$   & \multicolumn{1}{c|}{\begin{tabular}[c]{@{}c@{}}7.74e-2\\ 2m57s\end{tabular}} & \begin{tabular}[c]{@{}c@{}}6.31e-2\\ 2m57s\end{tabular} & \multicolumn{1}{c|}{\begin{tabular}[c]{@{}c@{}}4.38e-2\\ 3m37s\end{tabular}} & \begin{tabular}[c]{@{}c@{}}3.80e-2\\ 3m19s\end{tabular} \\ \hline
        $M=160$  & \multicolumn{1}{c|}{\begin{tabular}[c]{@{}c@{}}1.10e-1\\ 8m51s\end{tabular}}   & \begin{tabular}[c]{@{}c@{}}5.76e-2\\ 5m17s\end{tabular}  & \multicolumn{1}{c|}{\begin{tabular}[c]{@{}c@{}}1.09e-1\\ 15m34s\end{tabular}} & \begin{tabular}[c]{@{}c@{}}8.61e-2\\ 8m38s\end{tabular}  & $M=32$   & \multicolumn{1}{c|}{\begin{tabular}[c]{@{}c@{}}3.92e-2\\ 3m24s\end{tabular}} & \begin{tabular}[c]{@{}c@{}}2.93e-2\\ 3m01s\end{tabular} & \multicolumn{1}{c|}{\begin{tabular}[c]{@{}c@{}}3.93e-2\\ 5m06s\end{tabular}} & \begin{tabular}[c]{@{}c@{}}3.28e-2\\ 3m39s\end{tabular} \\ \hline
        $M=320$  & \multicolumn{1}{c|}{\begin{tabular}[c]{@{}c@{}}8.20e-2\\ 15 m33s\end{tabular}} & \begin{tabular}[c]{@{}c@{}}4.11e-2\\ 8m34s\end{tabular}  & \multicolumn{1}{c|}{\begin{tabular}[c]{@{}c@{}}9.67e-2\\ 24m02s\end{tabular}} & \begin{tabular}[c]{@{}c@{}}4.44e-2\\ 10m52s\end{tabular} & $M=48$   & \multicolumn{1}{c|}{\begin{tabular}[c]{@{}c@{}}3.05e-2\\ 4m12s\end{tabular}} & \begin{tabular}[c]{@{}c@{}}4.83e-2\\ 3m21s\end{tabular} & \multicolumn{1}{c|}{\begin{tabular}[c]{@{}c@{}}4.82e-2\\ 6m19s\end{tabular}} & \begin{tabular}[c]{@{}c@{}}6.01e-2\\ 4m13s\end{tabular} \\ \hline
        $M=640$  & \multicolumn{1}{c|}{\begin{tabular}[c]{@{}c@{}}1.05e-1\\ 24m04s\end{tabular}}  & \begin{tabular}[c]{@{}c@{}}7.65e-2\\ 10m40s\end{tabular} & \multicolumn{1}{c|}{\begin{tabular}[c]{@{}c@{}}7.94e-2\\ 45m10s\end{tabular}} & \begin{tabular}[c]{@{}c@{}}5.54e-2\\ 19m45s\end{tabular} & $M=64$   & \multicolumn{1}{c|}{\begin{tabular}[c]{@{}c@{}}4.87e-2\\ 5m02s\end{tabular}} & \begin{tabular}[c]{@{}c@{}}4.46e-2\\ 3m34s\end{tabular} & \multicolumn{1}{c|}{\begin{tabular}[c]{@{}c@{}}3.21e-2\\ 7m38s\end{tabular}} & \begin{tabular}[c]{@{}c@{}}3.78e-2\\ 4m35s\end{tabular} \\ \hline
        $M=1280$ & \multicolumn{1}{c|}{\begin{tabular}[c]{@{}c@{}}7.0e-2\\ 45m10s\end{tabular}}   & \begin{tabular}[c]{@{}c@{}}5.82e-2\\ 19m44s\end{tabular} & \multicolumn{1}{c|}{OoM}                                                      & \begin{tabular}[c]{@{}c@{}}5.64e-2\\ 37m25s\end{tabular} & $M=80$   & \multicolumn{1}{c|}{\begin{tabular}[c]{@{}c@{}}6.29e-2\\ 5m27s\end{tabular}} & \begin{tabular}[c]{@{}c@{}}6.73e-2\\ 3m56s\end{tabular} & \multicolumn{1}{c|}{\begin{tabular}[c]{@{}c@{}}5.35e-2\\ 8m50s\end{tabular}} & \begin{tabular}[c]{@{}c@{}}3.66e-2\\ 5m17s\end{tabular} \\ \hline
        \end{tabular}

    \caption{$L^2$ relative error and computational time for equation \eqref{eqn:eg1} with different $M$ and number of points. For each grid, the first row is the relative error, the second row is the average computational time for 5000 epochs. OoM means out of memory.}
    \label{tab:forward}
\end{table}

\begin{figure}[htbp]
    \centering
    \begin{subfigure}{.25\textwidth}
        \centering
        \includegraphics[height=0.75\textwidth,width=1.0\textwidth]{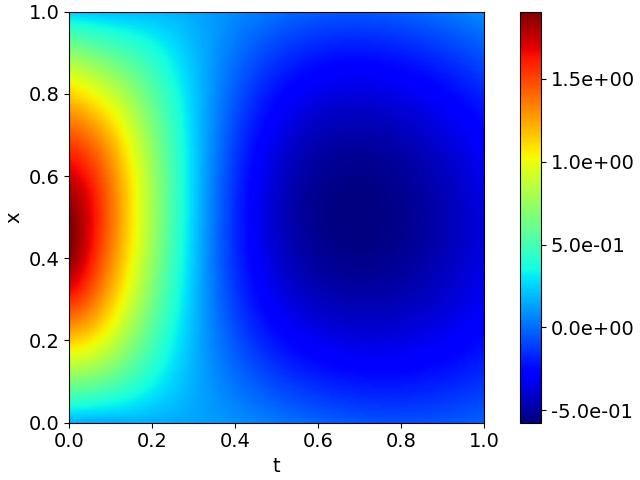}
    \end{subfigure}%
    \begin{subfigure}{.25\textwidth}
        \centering
        \includegraphics[height=0.75\textwidth,width=1.0\textwidth]{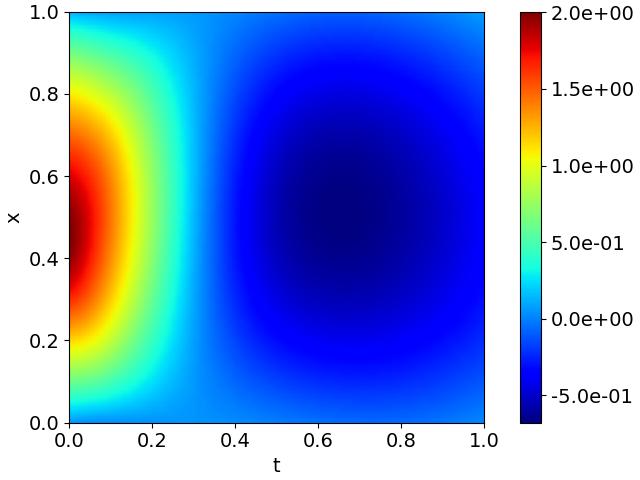}
    \end{subfigure}%
    \begin{subfigure}{.25\textwidth}
        \centering
        \includegraphics[height=0.75\textwidth,width=1.0\textwidth]{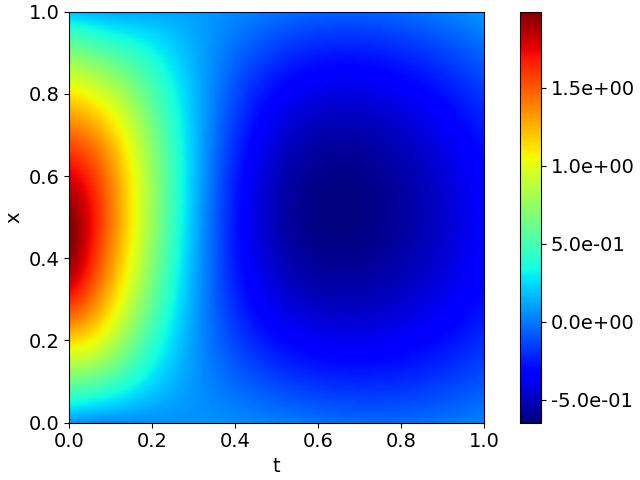}
    \end{subfigure}%
    \begin{subfigure}{.25\textwidth}
        \centering
        \includegraphics[height=0.75\textwidth,width=1.0\textwidth]{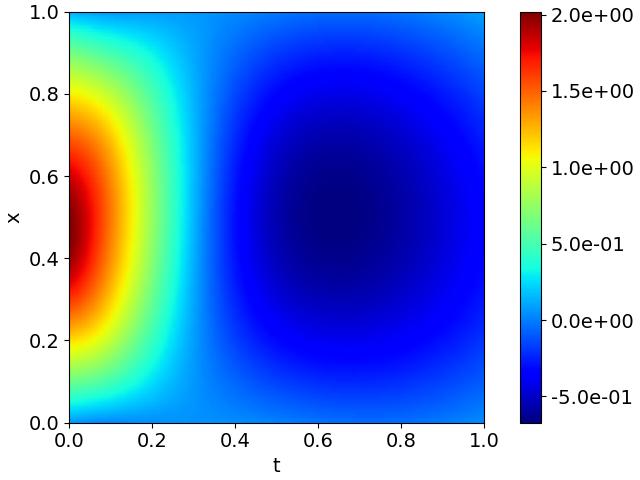}
    \end{subfigure}%
    \newline
    \raggedleft
    \begin{subfigure}{.25\textwidth}
        \centering
        \includegraphics[height=0.75\textwidth,width=1.0\textwidth]{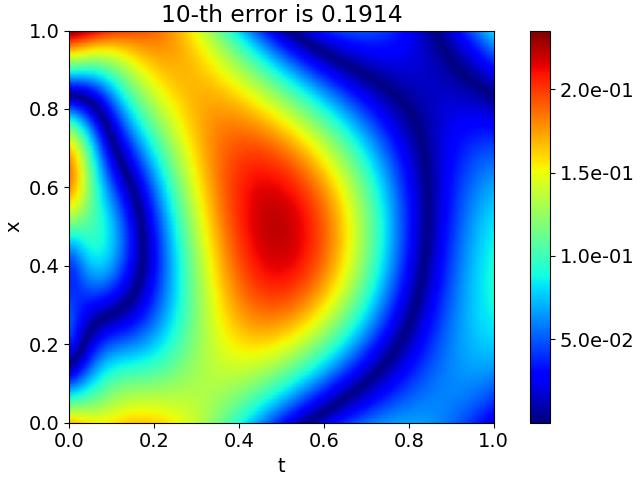}
        \caption{\textbf{MC-I}, $M=80$.}
    \end{subfigure}%
    \begin{subfigure}{.25\textwidth}
        \centering
        \includegraphics[height=0.75\textwidth,width=1.0\textwidth]{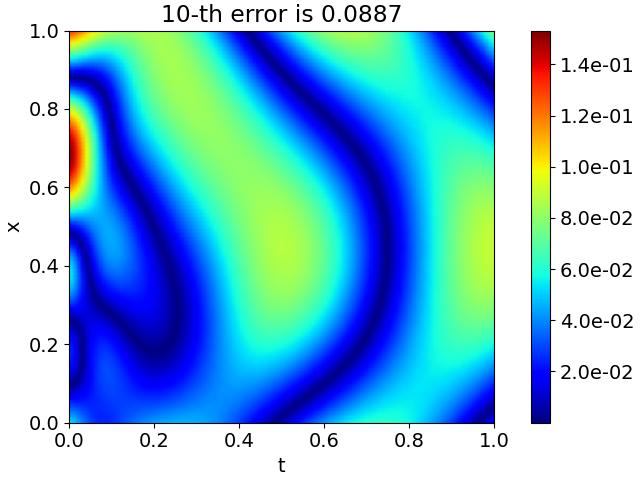}
        \caption{\textbf{MC-II}, $M=80$.}
    \end{subfigure}%
    \begin{subfigure}{.25\textwidth}
        \centering
        \includegraphics[height=0.75\textwidth,width=1.0\textwidth]{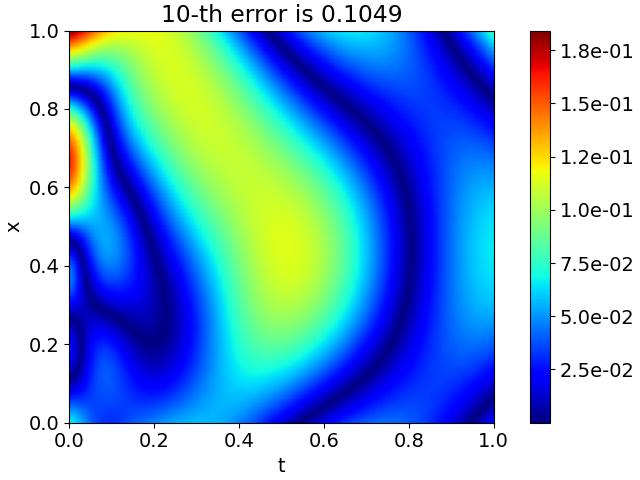}
        \caption{\textbf{MC-I}, $M=640$.}
    \end{subfigure}%
    \begin{subfigure}{.25\textwidth}
        \centering
        \includegraphics[height=0.75\textwidth,width=1.0\textwidth]{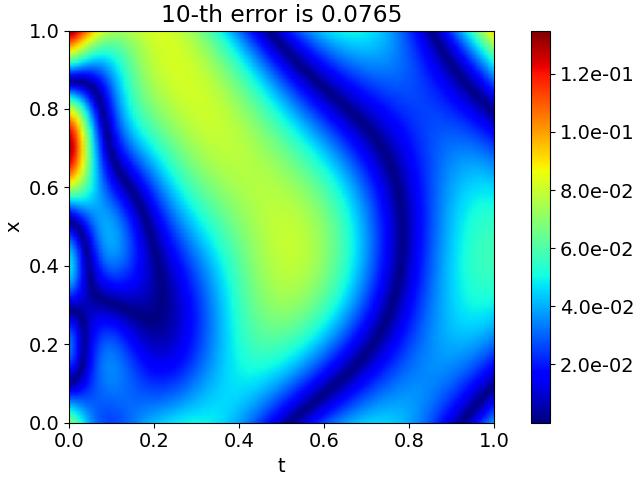}
        \caption{\textbf{MC-II}, $M=640$.}
    \end{subfigure}%
    \caption{Profiles of absolute error and neural network solutions for equation \eqref{eqn:eg1} with the Monte Carlo method for $M=80$ and $M=640$. First row: numerical solutions. Second row: absolute error.}
    \label{fig:eg1sol}
\end{figure}

In the following, we present solution profiles derived from the results in Table \ref{tab:forward}. We first compare the cases of the Monte Carlo method in Figure \ref{fig:eg1sol}, which displays the solution profiles and absolute errors for each numerical method.  The results indicate that improved accuracy in the computation of fractional derivatives leads to more effective training outcomes. At relatively low values of $M$, the Monte Carlo methods suffer from reduced accuracy, emphasizing the advantage of Gauss-Jacobi methods in achieving better performance with fewer quadrature points. This highlights the superior computational efficiency of Gauss-Jacobi quadratures in practice.

\begin{figure}[htbp]
    \centering
    \begin{subfigure}{.25\textwidth}
        \centering
        \includegraphics[height=0.75\textwidth,width=1.0\textwidth]{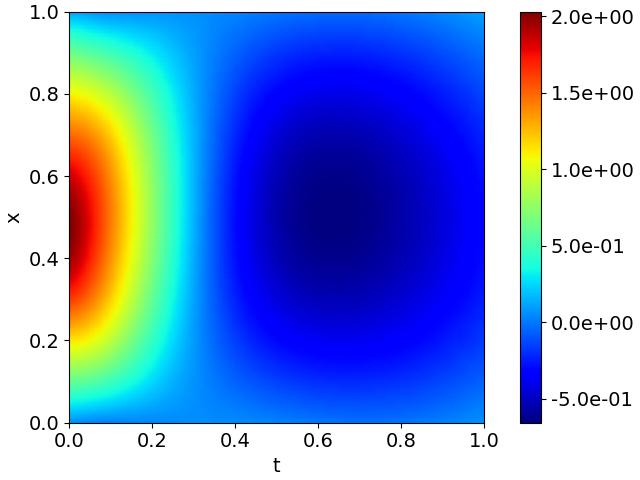}
    \end{subfigure}%
    \begin{subfigure}{.25\textwidth}
        \centering
        \includegraphics[height=0.75\textwidth,width=1.0\textwidth]{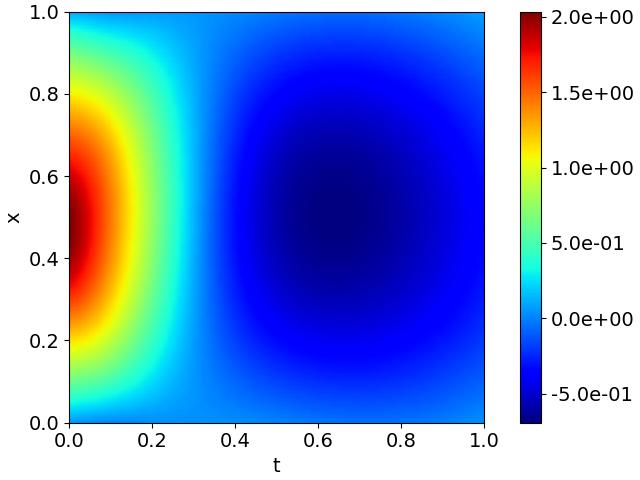}
    \end{subfigure}%
    \begin{subfigure}{.25\textwidth}
        \centering
        \includegraphics[height=0.75\textwidth,width=1.0\textwidth]{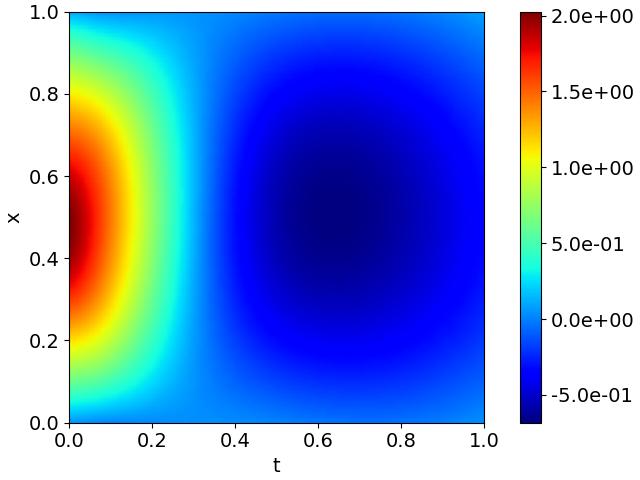}
    \end{subfigure}%
    \begin{subfigure}{.25\textwidth}
        \centering
        \includegraphics[height=0.75\textwidth,width=1.0\textwidth]{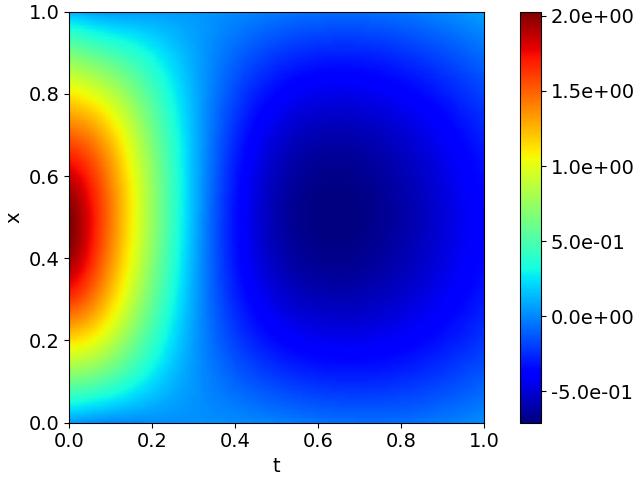}
    \end{subfigure}%
    \newline
    \raggedleft
    \begin{subfigure}{.25\textwidth}
        \centering
        \includegraphics[height=0.75\textwidth,width=1.0\textwidth]{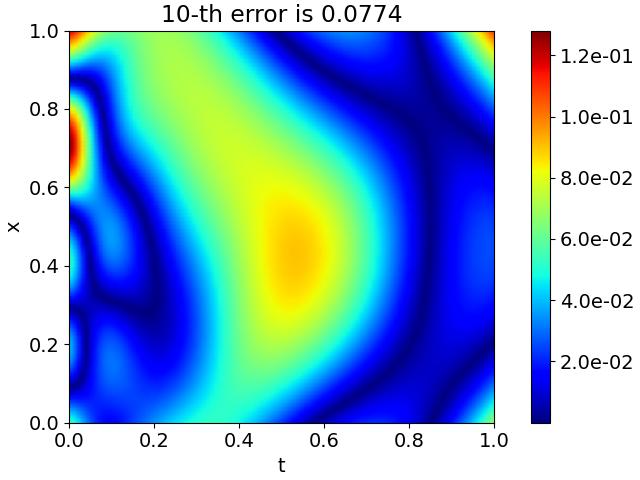}
        \caption{\textbf{GJ-I}, $M=16$.}
    \end{subfigure}%
    \begin{subfigure}{.25\textwidth}
        \centering
        \includegraphics[height=0.75\textwidth,width=1.0\textwidth]{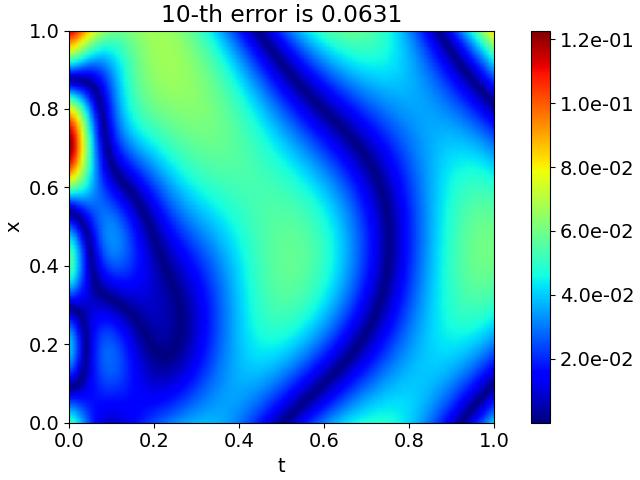}
        \caption{\textbf{GJ-II}, $M=16$.}
    \end{subfigure}%
    \begin{subfigure}{.25\textwidth}
        \centering
        \includegraphics[height=0.75\textwidth,width=1.0\textwidth]{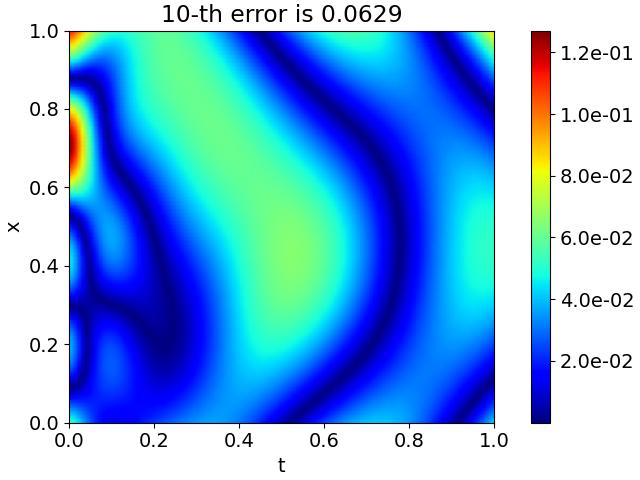}
        \caption{\textbf{GJ-I}, $M=80$.}
    \end{subfigure}%
    \begin{subfigure}{.25\textwidth}
        \centering
        \includegraphics[height=0.75\textwidth,width=1.0\textwidth]{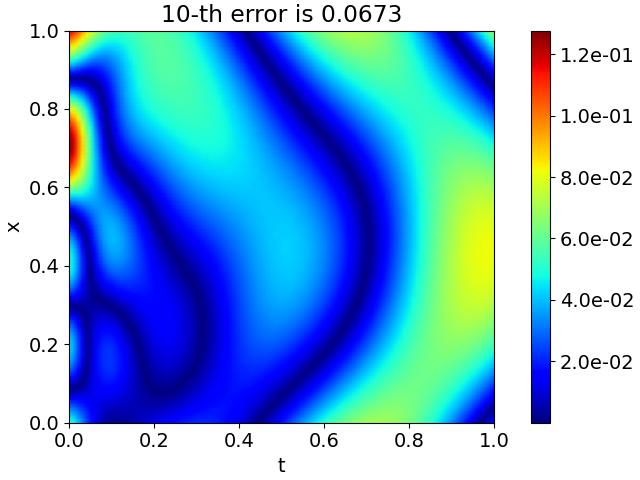}
        \caption{\textbf{GJ-II}, $M=80$.}
    \end{subfigure}%
    \caption{Profiles of absolute error and neural network solutions for equation \eqref{eqn:eg1} with the Gauss-Jacobi method for $M=16$ and $M=80$. First row: numerical solutions. Second row: absolute error.}
    \label{fig:eg1sol640}
\end{figure}

However, the Gauss-Jacobi quadrature method demonstrates remarkable efficiency in terms of the number of quadrature points required. As illustrated in Figure \ref{fig:eg1sol640}, which compares solution profiles and absolute error distributions for implementations with $M=16$ and $M=80$, the method achieves comparable accuracy levels with significantly fewer quadrature points.

\subsection{Time-fractional Diffusion-Wave equation with adjusted parameters}
\label{subsec:eg2}
Here we consider the following time-fractional PDEs in $\Omega=(0,1)$:
\begin{equation}
    \label{eqn:eg2}
    \begin{aligned}
        \dta u(t,x) - \frac{\lambda }{k^2\pi^2} u_{xx}(t,x) &= 0,~~ (t,\x)\in (0,T]\times\Omega, ~~\lambda\in\mathbb{R}^+,~~ k\in\mathbb{N}, \\
        u(t,0)=u(t,1)&=0,~~t\in(0,T],\\
        u(0,x)&=\sin(k\pi x), ~~x\in\Omega,\\
        \partial_tu(0,x)&=-0.5\sin(k\pi x), ~~x\in\Omega,
    \end{aligned}
\end{equation}
with  exact solution
\begin{equation*}
    u(t,x)=  (E_{\al, 1}(-\lambda t^\al)-0.5tE_{\al,2}(-\lambda t^\al))\sin(k\pi x).
\end{equation*}

\begin{figure}[htbp]
    \centering
    \begin{subfigure}{.25\textwidth}
        \centering
        \includegraphics[height=0.75\textwidth,width=1.0\textwidth]{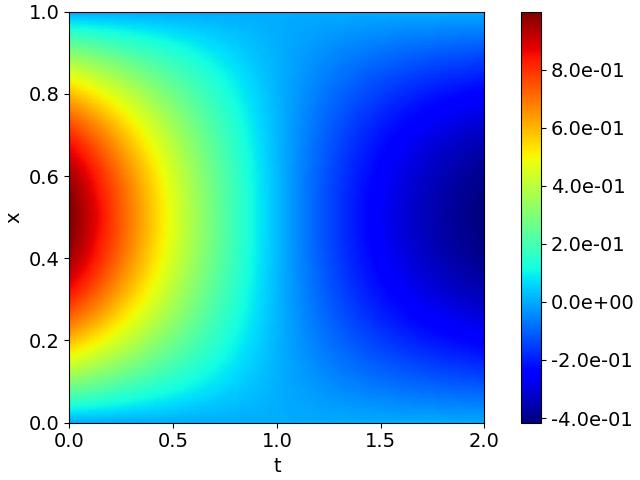}
    \end{subfigure}%
    \begin{subfigure}{.25\textwidth}
        \centering
        \includegraphics[height=0.75\textwidth,width=1.0\textwidth]{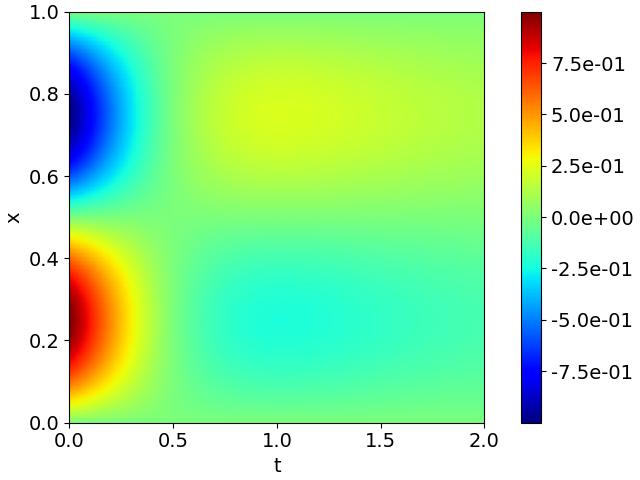}
    \end{subfigure}%
    \begin{subfigure}{.25\textwidth}
        \centering
        \includegraphics[height=0.75\textwidth,width=1.0\textwidth]{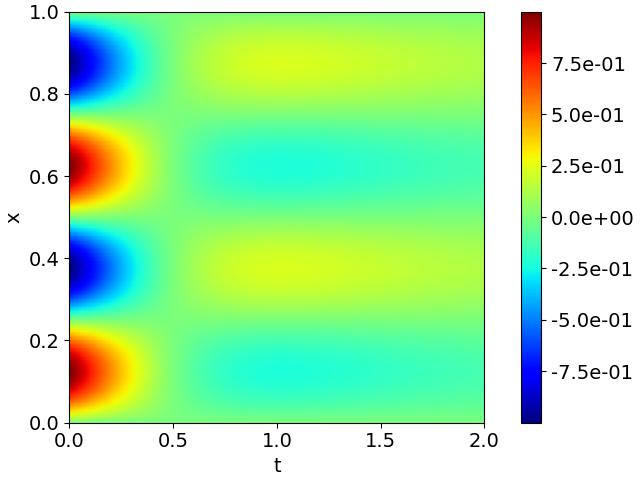}
    \end{subfigure}%
    \begin{subfigure}{.25\textwidth}
        \centering
        \includegraphics[height=0.75\textwidth,width=1.0\textwidth]{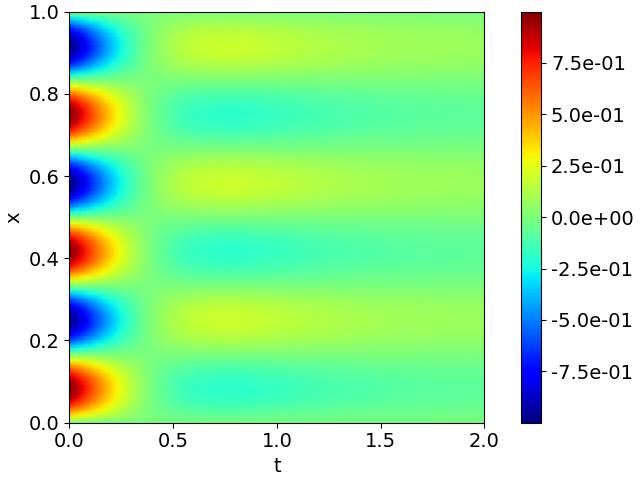}
    \end{subfigure}%
    \newline
    \raggedleft
    \begin{subfigure}{.25\textwidth}
        \centering
        \includegraphics[height=0.75\textwidth,width=1.0\textwidth]{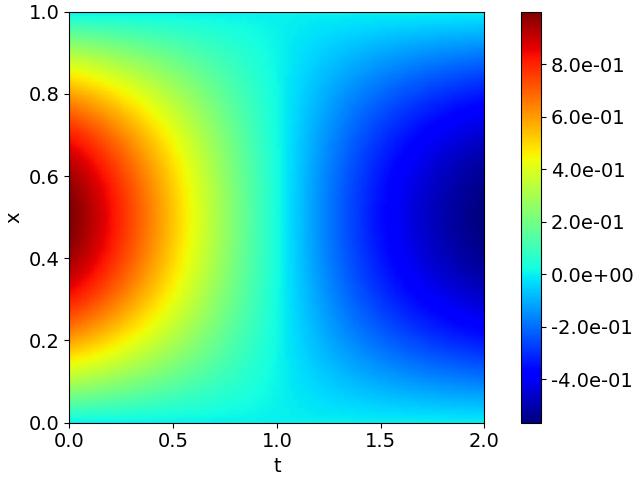}
    \end{subfigure}%
    \begin{subfigure}{.25\textwidth}
        \centering
        \includegraphics[height=0.75\textwidth,width=1.0\textwidth]{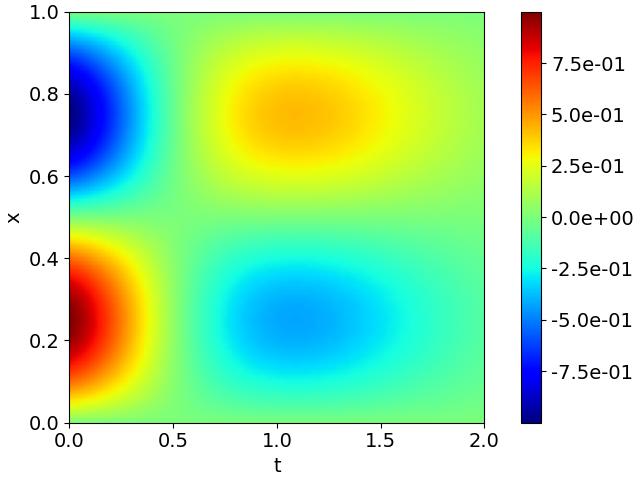}
    \end{subfigure}%
    \begin{subfigure}{.25\textwidth}
        \centering
        \includegraphics[height=0.75\textwidth,width=1.0\textwidth]{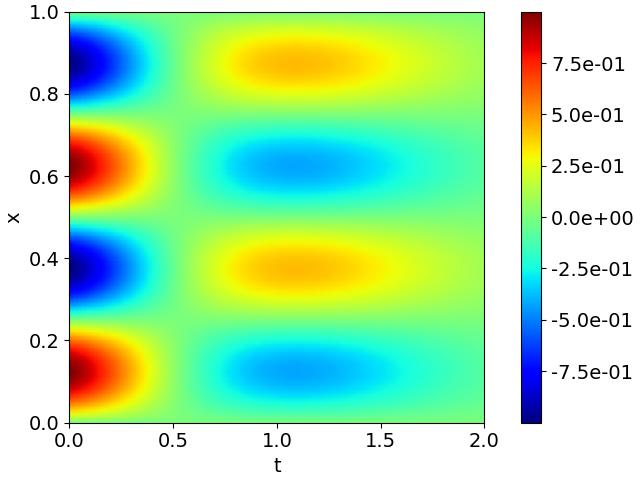}
    \end{subfigure}%
    \begin{subfigure}{.25\textwidth}
        \centering
        \includegraphics[height=0.75\textwidth,width=1.0\textwidth]{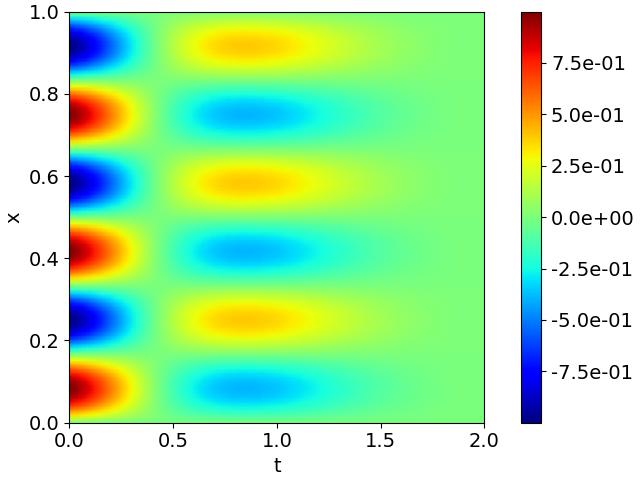}
    \end{subfigure}%
    \newline
    \raggedleft
    \begin{subfigure}{.25\textwidth}
        \centering
        \includegraphics[height=0.75\textwidth,width=1.0\textwidth]{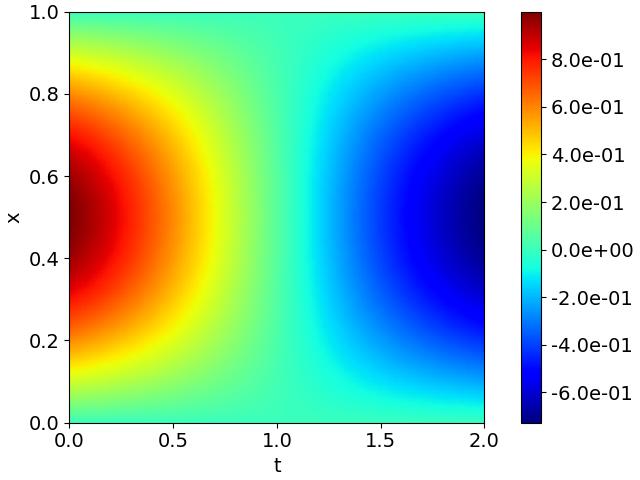}
        \caption{\textbf{$k=1$, $\lambda=1$}}
    \end{subfigure}%
    \begin{subfigure}{.25\textwidth}
        \centering
        \includegraphics[height=0.75\textwidth,width=1.0\textwidth]{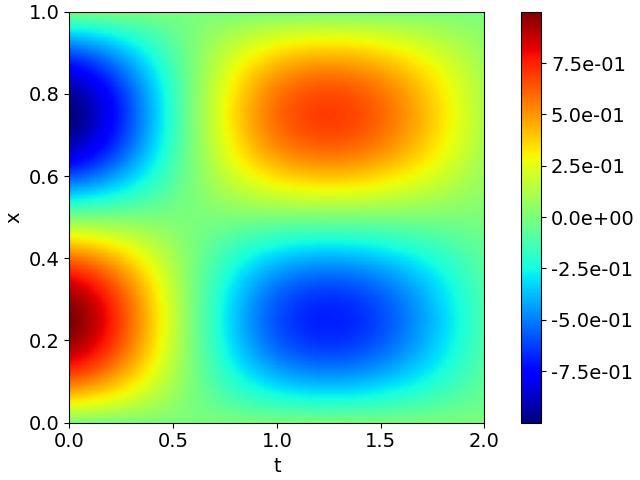}
        \caption{$k=2$, $\lambda=4$}
    \end{subfigure}%
    \begin{subfigure}{.25\textwidth}
        \centering
        \includegraphics[height=0.75\textwidth,width=1.0\textwidth]{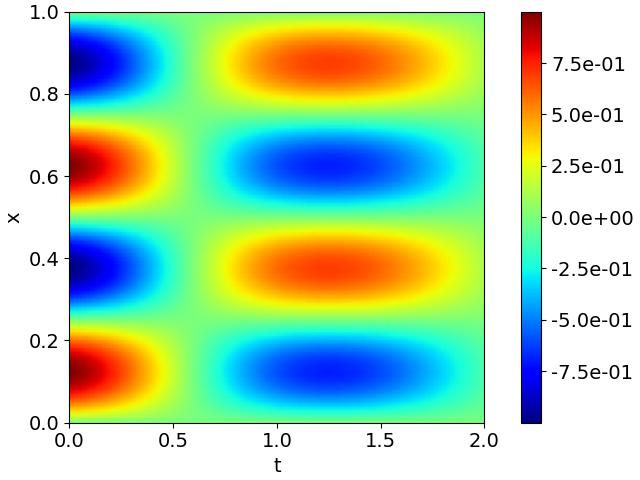}
        \caption{$k=4$, $\lambda=4$}
    \end{subfigure}%
    \begin{subfigure}{.25\textwidth}
        \centering
        \includegraphics[height=0.75\textwidth,width=1.0\textwidth]{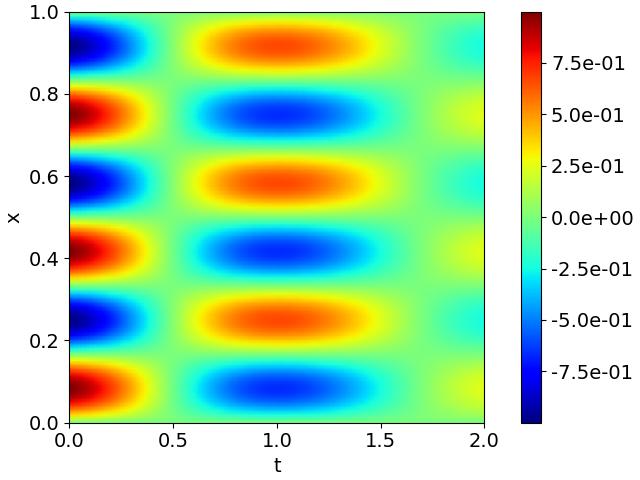}
        \caption{$k=6$, $\lambda=6$}
    \end{subfigure}%
    \caption{Exact solution of cases in Table \ref{tab:eg2-1} with different $\alpha$, $k$ and $\lambda$. In the first row, $\alpha=1.25$, in the second row, $\alpha=1.5$, in the third row, $\alpha=1.75$.}
    \label{fig:eg2-exact}
\end{figure}

Equation \eqref{eqn:eg1} can be regarded as a special case with $\lambda = \pi^2$, $k=1$, scaled initial conditions, and $T=1$. In Table \ref{tab:eg2-1}, the final time is set to $T=2$. Small batches, consistent with those in subsection \ref{subsec:eg1}, are adopted. The training procedure is configured with a maximum of  $L=3$ iterations, each consisting of 5000 epochs and optimized using the Adam optimizer with a learning rate of  $0.001$. The exact solutions for each case and each $\alpha$ are presented in Figure \ref{fig:eg2-exact}, allowing a clear comparison of the diffusion and wave components in the solutions. As $\alpha$ increases, the wave component becomes more pronounced, while larger values of the diffusion coefficient $\lambda/(k^2\pi^2)$ enhance the influence of the diffusion effects. This interplay between the wave and diffusion components helps explain the numerical behaviors observed in the subsequent results.

We present the relative errors and computational time in Table \ref{tab:eg2-1} for different numerical integration schemes and values of $M$. The proposed tDWfPINNs demonstrate accuracy on par with direct numerical integration methods, while significantly reducing computational time. These extensive numerical experiments further validate the efficiency of our approach in solving the time-fractional diffusion-wave equation.
\begin{table}
\centering

\begin{tblr}{
  cells = {c},
  cell{1}{3} = {c=2}{},
  cell{1}{5} = {c=2}{},
  cell{1}{8} = {c=2}{},
  cell{1}{10} = {c=2}{},
  cell{3}{1} = {r=2}{},
  cell{5}{1} = {r=2}{},
  cell{7}{1} = {r=2}{},
  cell{9}{3} = {c=2}{},
  cell{9}{5} = {c=2}{},
  cell{9}{8} = {c=2}{},
  cell{9}{10} = {c=2}{},
  cell{10}{1} = {r=2}{},
  cell{12}{1} = {r=2}{},
  cell{14}{1} = {r=2}{},
  vlines,
  hline{1-3,5,7,9-10,12,14,16} = {-}{},
  hline{4,6,8,11,13,15} = {2-11}{},
}
         &         & $k=1$, $\lambda=1$ &                    & $k=2$, $\lambda=4$ &                    &        & $k=1$, $\lambda=1$  &                   & $k=2$, $\lambda=4$ &                   \\
$\alpha$ &         & \textbf{MC-I}      & \textbf{MC-II}     & \textbf{MC-I}      & \textbf{MC-II}     &        & \textbf{GJ-I}       & \textbf{GJ-II}    & \textbf{GJ-I}      & \textbf{GJ-II}    \\
$1.25$   & $M=80$  & {2.16e-2\\ 5m38s}  & {3.15e-2\\ 4m04s}  & {4.58e-2\\ 5m39s}  & {4.63e-2\\ 3m59s}  & $M=16$ & {5.80e-3\\ 3m13s}   & {1.03e-2\\ 2m51s} & {4.74e-2\\ 3m11s}  & {7.06e-2\\ 3m01s} \\
         & $M=640$ & {1.93e-2\\ 24m00s} & {1.10e-2\\ 10m52s} & {6.75e-2\\ 23m57s} & {3.76e-2\\ 10m53s} & $M=80$ & {7.10e-3\\ 5m31s}   & {7.60e-3\\ 4m01s} & {6.32e-2\\ 5m35s}  & {3.96e-2\\ 3m51s} \\
$1.5$    & $M=80$  & {3.83e-2\\ 5m36s}  & {8.50e-3\\ 4m06s}  & {6.51e-2\\ 5m36s}  & {4.75e-2\\ 3m52s}  & $M=16$ & {1.20e-2\\ 3m08s}   & {1.14e-2\\ 2m59s} & {3.24e-2\\ 3m07s}  & {4.21e-2\\ 3m01s} \\
         & $M=640$ & {4.35e-2\\ 23m55s} & {2.35e-2\\ 10m49s} & {5.13e-2\\ 23m56s} & {2.48e-2\\ 10m53s} & $M=80$ & {6.10e-3\\ 5m31s}   & {1.82e-2\\ 4m04s} & {2.74e-2\\ 5m31s}  & {2.20e-2\\ 3m53s} \\
$1.75$   & $M=80$  & {7.20e-3\\ 5m32s}  & {7.10e-3\\ 4m05s}  & {6.39e-2\\ 5m34s}  & {2.59e-2\\ 3m54s}  & $M=16$ & {4.1e-3\\ 3m15s}    & {2.40e-3\\ 3m00s} & {1.64e-2\\ 3m07s}  & {2.00e-2\\ 2m59s} \\
         & $M=640$ & {4.50e-3\\ 24m02s} & {1.74e-2\\ 10m52s} & {2.64e-2\\ 23m58s} & {2.49e-2\\ 10m52s} & $M=80$ & {1.19e-2\\ 5m32s}   & {8.90e-3\\ 4m00s} & {1.82e-2\\ 5m31s}  & {9.30e-3\\ 3m55s} \\
         &         & $k=4$, $\lambda=4$ &                    & $k=6$, $\lambda=6$ &                    &        & $k=4$, $\lambda =4$ &                   & $k=6$, $\lambda=6$ &                   \\
$1.25$   & $M=80$  & {5.69e-2\\ 5m39s}  & {3.91e-2\\ 3m59s}  & {1.52e-1\\ 5m39s}  & {2.98e-1\\ 3m54s}  & $M=16$ & {6.34e-2\\ 3m07s}   & {6.23e-2\\ 2m59s} & {1.09e-1\\ 3m03s}  & {9.56e-2\\ 3m00s} \\
         & $M=640$ & {2.76e-2\\ 23m59s} & {5.13e-2\\ 10m47s} & {1.01e-1\\ 23m53s} & {1.22e-1\\ 10m52s} & $M=80$ & {3.80e-2\\ 5m36s}   & {5.22e-2\\ 4m01s} & {1.04e-1\\ 5m31s}  & {1.20e-1\\ 4m01s}  \\
$1.5$    & $M=80$  & {4.64e-2\\ 5m41s}  & {3.94e-2\\ 4m06s}  & {6.41e-2\\ 5m33s}  & {5.68e-2\\ 4m10s}  & $M=16$ & {3.60e-2\\ 3m07s}   & {3.96e-2\\ 3m00s} & {4.69e-2\\ 3m05s}  & {5.41e-2\\ 2m56s} \\
         & $M=640$ & {4.65e-2\\ 23m52s} & {6.76e-2\\ 10m52s} & {8.74e-2\\ 23m54s} & {5.28e-2\\ 10m54s} & $M=80$ & {2.84e-2\\ 5m33s}   & {2.25e-2\\ 3m58s} & {4.29e-2\\ 5m39s}  & {8.75e-2\\ 4m07s} \\
$1.75$   & $M=80$  & {9.74e-2\\ 5m35s}  & {3.78e-2\\ 3m55s}  & {6.10e-2\\ 5m42s}  & {8.12e-2\\ 4m05s}  & $M=16$ & {7.48e-2\\3m07s}    & {3.37e-2\\3m02s}  & {5.45e-2\\ 3m05s}  & {3.44e-2\\2m59s}  \\
         & $M=640$ & {2.59e-2\\ 23m58s} & {3.88e-2\\ 10m51s} & {5.28e-2\\ 23m55s} & {6.71e-2\\ 10m52s} & $M=80$ & {3.62e-2\\ 5m31s}   & {5.13e-2\\ 3m54s} & {4.76e-2\\ 5m30s}  & {1.99e-2\\ 4m04s}
\end{tblr}
\caption{Results of different numerical integrations with different parameters in equation \eqref{eqn:eg2}. In each grid, the first row is the $L^2$ relative error, the second row is the computational time for each case within 5000 epochs.}
\label{tab:eg2-1}
\end{table}

As the training process for PINNs progresses, our experiments consistently demonstrate that the Gauss-Jacobi method outperforms the Monte Carlo approach. However, we observe that the performance of both methods tends to deteriorate as the equation parameter increases, implying an increase in solution complexity.

Examining Table \ref{tab:eg2-1}, we observe that for cases involving larger values of $\lambda$ and $k$, models with higher fractional order $\alpha$ tend to produce more accurate results than those with smaller $\alpha$. This trend is likely attributable to the strong decay in the Mittag-Leffler scaling factor for small $\alpha$ when $\lambda$ is large, as demonstrated in Figure \ref{fig:eg2-exact}.

\begin{figure}[htbp]
    \centering
    \begin{subfigure}{.25\textwidth}
        \centering
        \includegraphics[height=0.75\textwidth,width=1.0\textwidth]{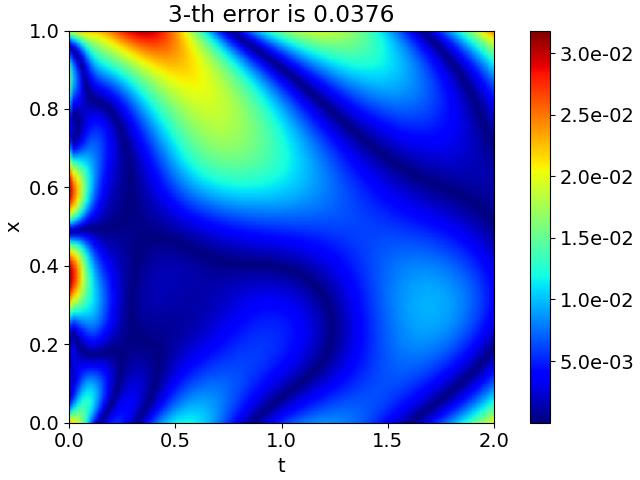}
    \end{subfigure}%
    \begin{subfigure}{.25\textwidth}
        \centering
        \includegraphics[height=0.75\textwidth,width=1.0\textwidth]{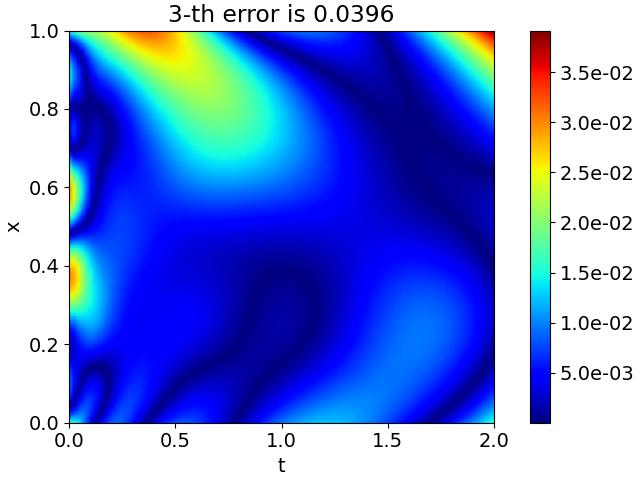}
    \end{subfigure}%
    \begin{subfigure}{.25\textwidth}
        \centering
        \includegraphics[height=0.75\textwidth,width=1.0\textwidth]{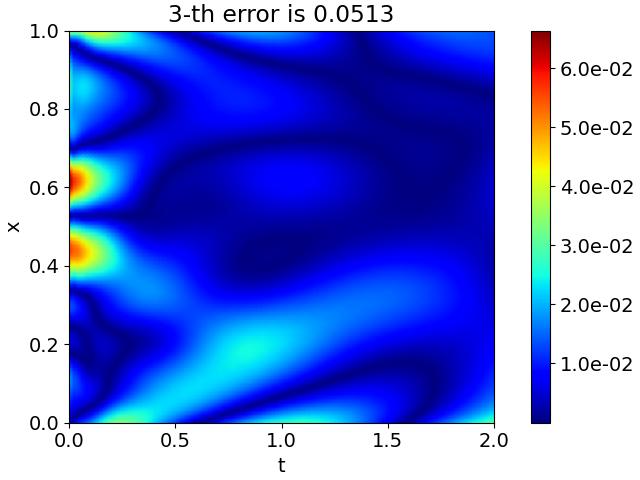}
    \end{subfigure}%
    \begin{subfigure}{.25\textwidth}
        \centering
        \includegraphics[height=0.75\textwidth,width=1.0\textwidth]{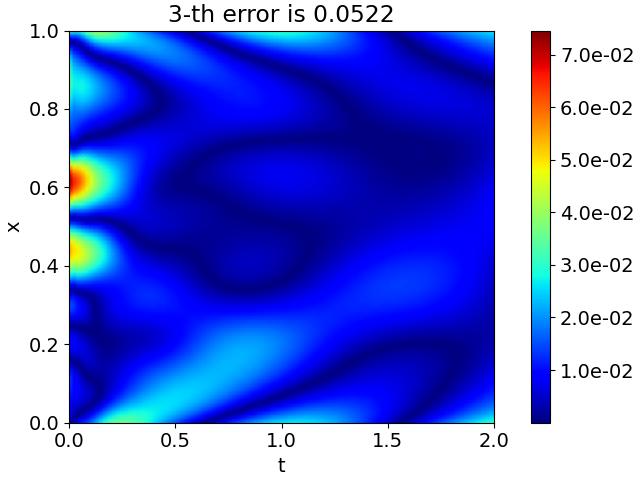}
    \end{subfigure}%
    \newline
    \raggedleft
    \begin{subfigure}{.25\textwidth}
        \centering
        \includegraphics[height=0.75\textwidth,width=1.0\textwidth]{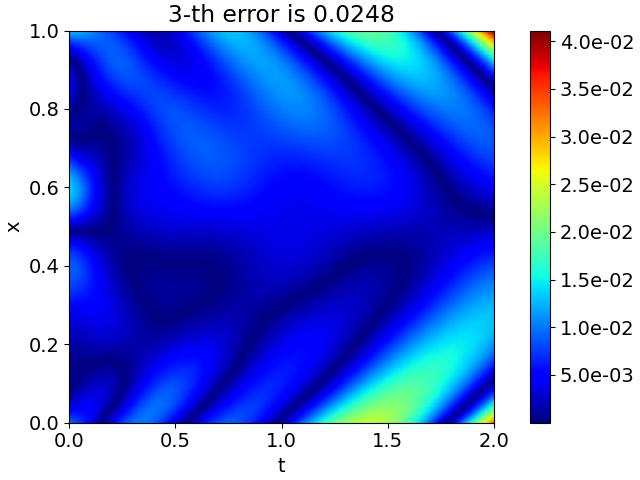}
    \end{subfigure}%
    \begin{subfigure}{.25\textwidth}
        \centering
        \includegraphics[height=0.75\textwidth,width=1.0\textwidth]{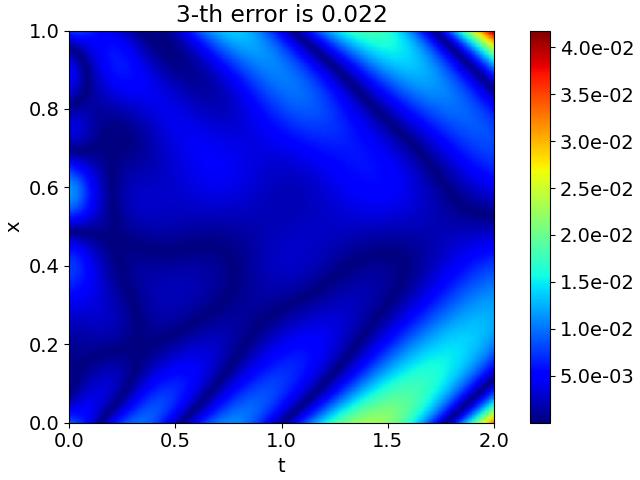}
    \end{subfigure}%
    \begin{subfigure}{.25\textwidth}
        \centering
        \includegraphics[height=0.75\textwidth,width=1.0\textwidth]{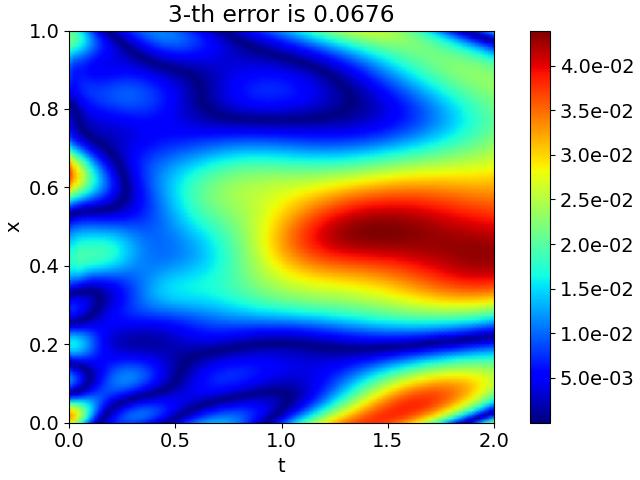}
    \end{subfigure}%
    \begin{subfigure}{.25\textwidth}
        \centering
        \includegraphics[height=0.75\textwidth,width=1.0\textwidth]{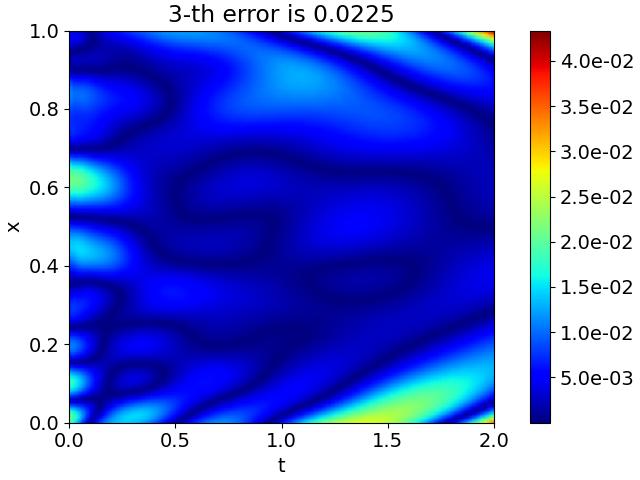}
    \end{subfigure}%
    \newline
    \raggedleft
    \begin{subfigure}{.25\textwidth}
        \centering
        \includegraphics[height=0.75\textwidth,width=1.0\textwidth]{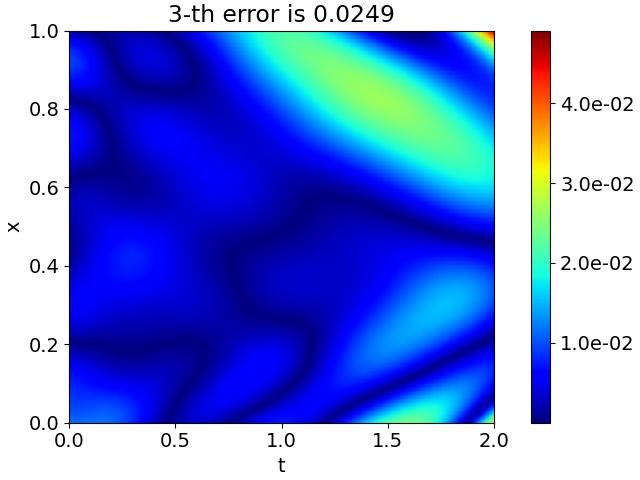}
        \caption{\textbf{MC-II}, $k=2$, $\lambda=4$, $M=640$.}
    \end{subfigure}%
    \begin{subfigure}{.25\textwidth}
        \centering
        \includegraphics[height=0.75\textwidth,width=1.0\textwidth]{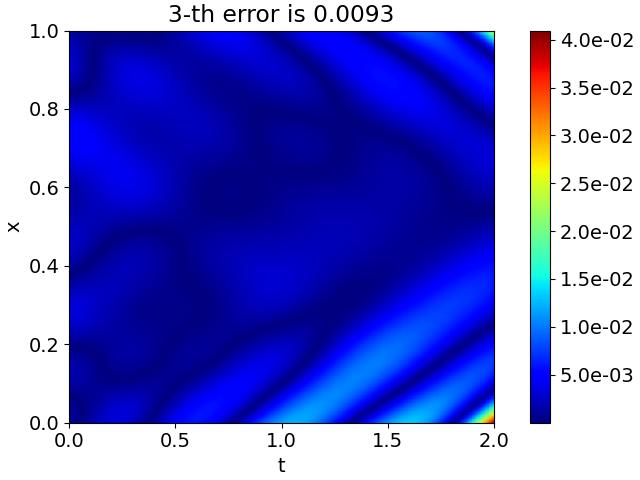}
        \caption{\textbf{GJ-II}, $k=2$, $\lambda=4$, $M=80$.}
    \end{subfigure}%
    \begin{subfigure}{.25\textwidth}
        \centering
        \includegraphics[height=0.75\textwidth,width=1.0\textwidth]{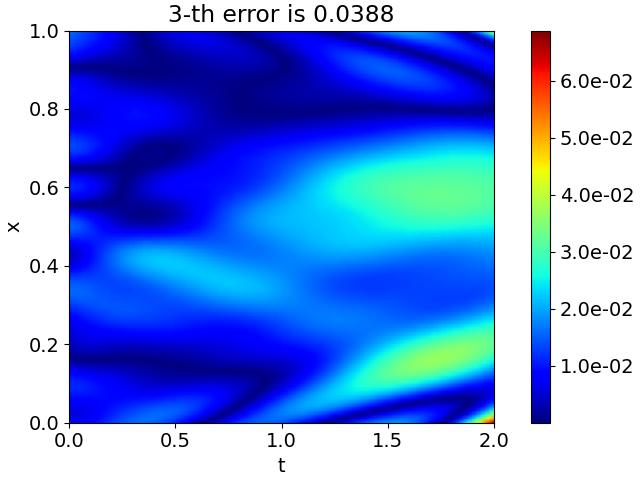}
        \caption{\textbf{MC-II}, $k=4$, $\lambda=4$, $M=640$.}
    \end{subfigure}%
    \begin{subfigure}{.25\textwidth}
        \centering
        \includegraphics[height=0.75\textwidth,width=1.0\textwidth]{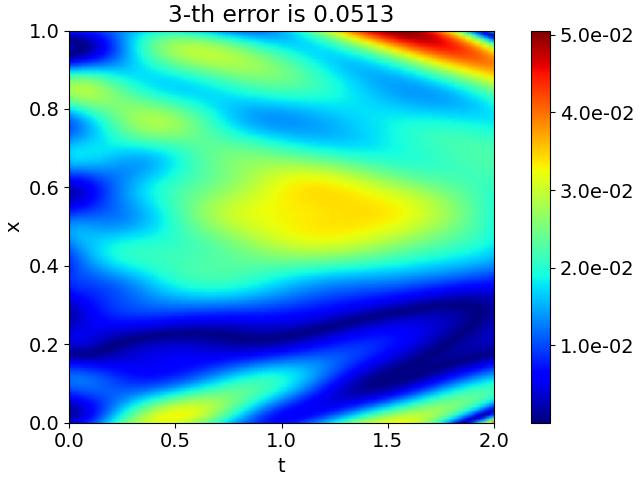}
        \caption{\textbf{GJ-II}, $k=4$, $\lambda=4$, $M=80$.}
    \end{subfigure}%
    \caption{Absolute errors of cases in Table \ref{tab:eg2-1} with different $\alpha$, $k$ and $\lambda$. In the first row, $\alpha=1.25$, in the second row, $\alpha=1.5$, in the third row, $\alpha=1.75$.}
    \label{fig:eg2-err-1}
\end{figure}

To further elucidate this behavior, Figure \ref{fig:eg2-err-1} presents the absolute error distributions of \textbf{MC-II} and \textbf{GJ-II} for cases with $k=2, \lambda=4$ and $k=4, \lambda=4$, evaluated across various $\alpha$ values when $M=640$ and $M=80$. The error patterns revealed in this figure lead to two significant observations. First, the error patterns obtained using the Gauss-Jacobi method tend to exhibit higher frequency content than those from the Monte Carlo method. This is indicative of a more effective representation of the underlying solution features and aligns with theoretical insights derived from neural tangent kernel analysis, such as those presented in \cite{WangYuParis:2022:NTK}. Second, significant error concentrations are consistently observed near the initial temporal layer and spatial boundaries. This suggests that the design of future training strategies, such as boundary-focused sampling or adaptive loss weighting, should prioritize these regions to improve overall approximation quality.

In the following experiments, the number $M$ is fixed at 80, and we employ the transformed schemes exclusively. This choice is motivated by earlier findings demonstrating that the transformed methods attain comparable accuracy to their direct counterparts while significantly reducing computational overhead.

\begin{figure}[htbp]
        \centering
        \begin{subfigure}{.25\textwidth}
            \centering
            \includegraphics[height=0.75\textwidth,width=1.0\textwidth]{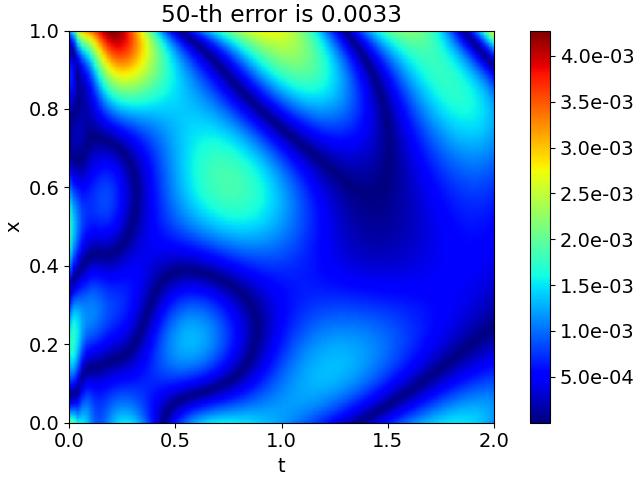}
        \end{subfigure}%
        \begin{subfigure}{.25\textwidth}
            \centering
            \includegraphics[height=0.75\textwidth,width=1.0\textwidth]{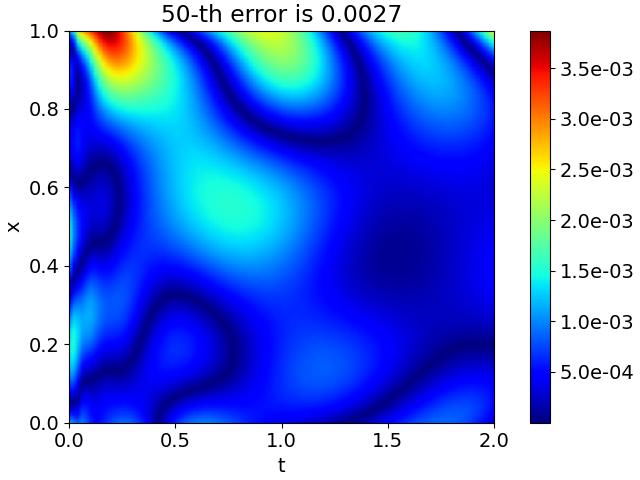}
        \end{subfigure}%
        \begin{subfigure}{.25\textwidth}
            \centering
            \includegraphics[height=0.75\textwidth,width=1.0\textwidth]{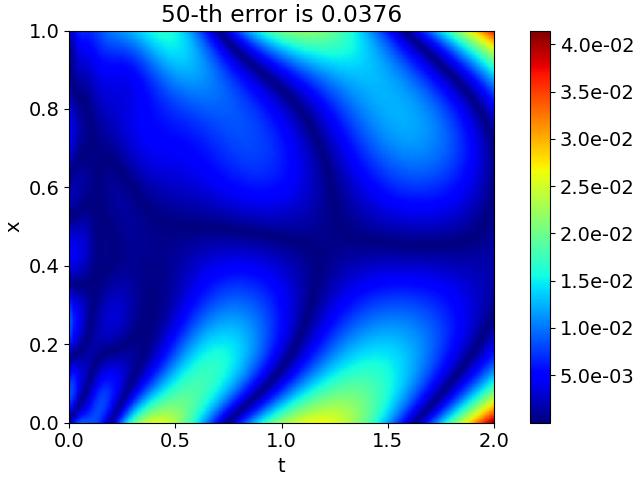}
        \end{subfigure}%
        \begin{subfigure}{.25\textwidth}
            \centering
            \includegraphics[height=0.75\textwidth,width=1.0\textwidth]{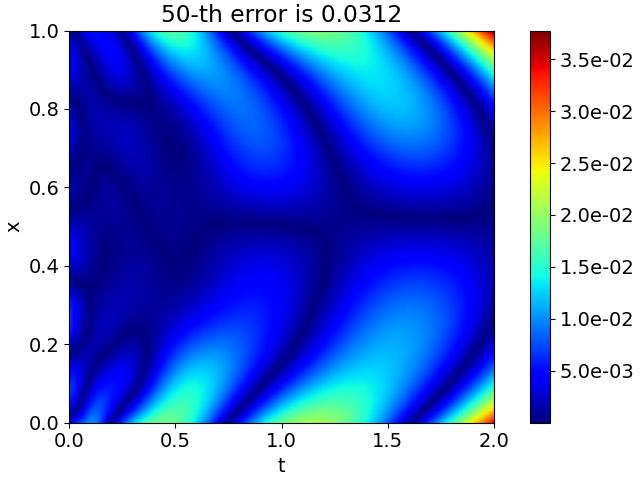}
        \end{subfigure}%
        \newline
        \raggedleft
        \begin{subfigure}{.25\textwidth}
            \centering
            \includegraphics[height=0.75\textwidth,width=1.0\textwidth]{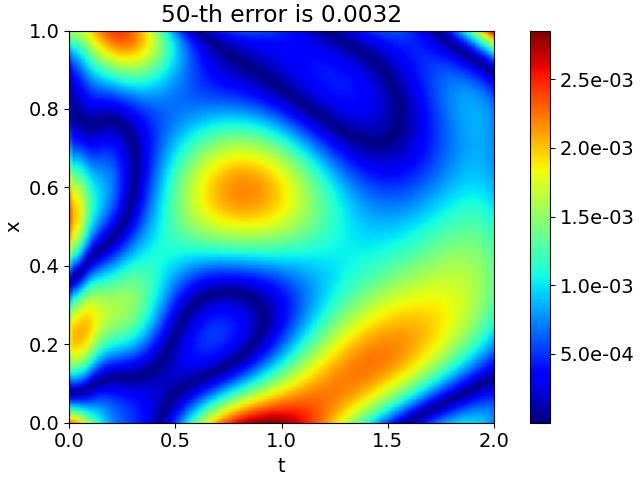}
        \end{subfigure}%
        \begin{subfigure}{.25\textwidth}
            \centering
            \includegraphics[height=0.75\textwidth,width=1.0\textwidth]{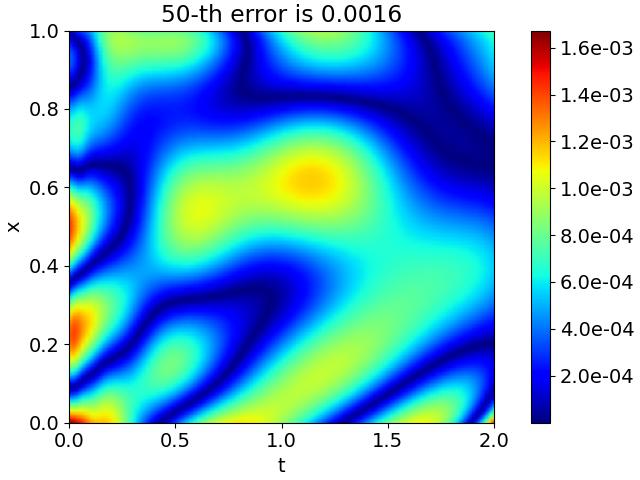}
        \end{subfigure}%
        \begin{subfigure}{.25\textwidth}
            \centering
            \includegraphics[height=0.75\textwidth,width=1.0\textwidth]{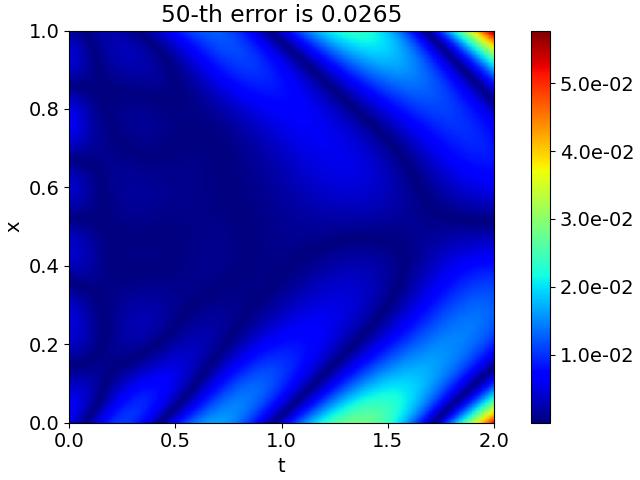}
        \end{subfigure}%
        \begin{subfigure}{.25\textwidth}
            \centering
            \includegraphics[height=0.75\textwidth,width=1.0\textwidth]{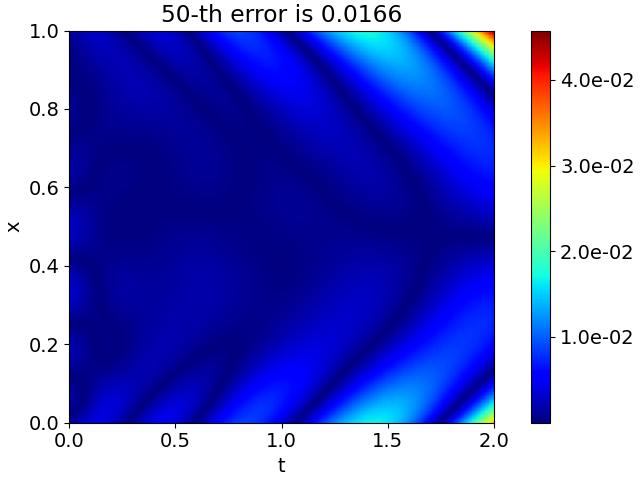}
        \end{subfigure}%
        \newline
        \raggedleft
        \begin{subfigure}{.25\textwidth}
            \centering
            \includegraphics[height=0.75\textwidth,width=1.0\textwidth]{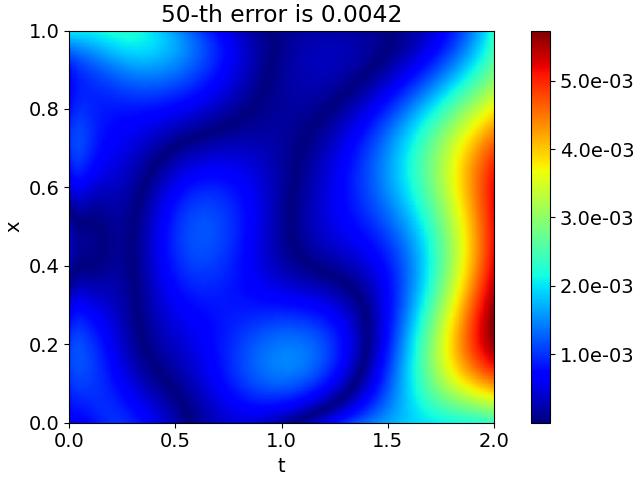}
            \caption{\textbf{MC-II}, $k=1$, $\lambda=1$.}
        \end{subfigure}%
        \begin{subfigure}{.25\textwidth}
            \centering
            \includegraphics[height=0.75\textwidth,width=1.0\textwidth]{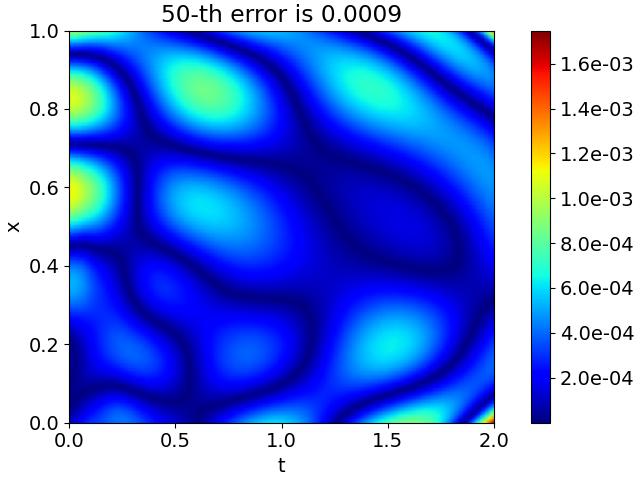}
            \caption{\textbf{GJ-II}, $k=1$, $\lambda=1$.}
        \end{subfigure}%
        \begin{subfigure}{.25\textwidth}
            \centering
            \includegraphics[height=0.75\textwidth,width=1.0\textwidth]{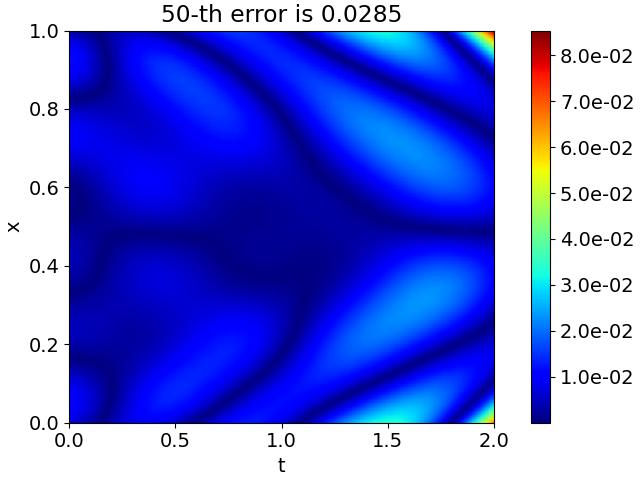}
            \caption{\textbf{MC-II}, $k=2$, $\lambda=4$.}
        \end{subfigure}%
        \begin{subfigure}{.25\textwidth}
            \centering
            \includegraphics[height=0.75\textwidth,width=1.0\textwidth]{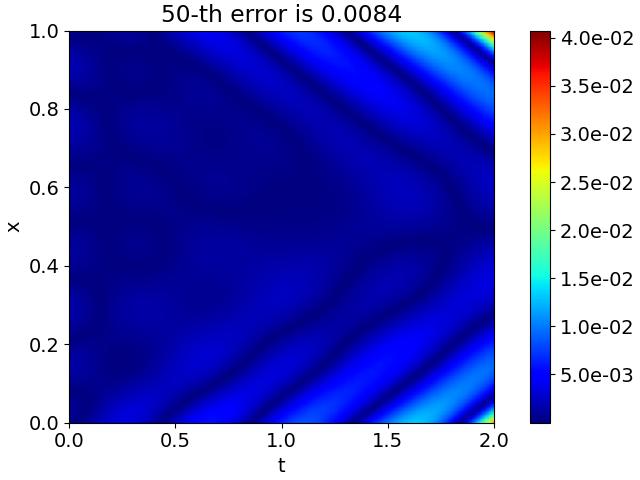}
            \caption{\textbf{GJ-II}, $k=2$, $\lambda=4$.}
        \end{subfigure}%
        \newline
        \raggedleft
        \begin{subfigure}{.25\textwidth}
            \centering
            \includegraphics[height=0.75\textwidth,width=1.0\textwidth]{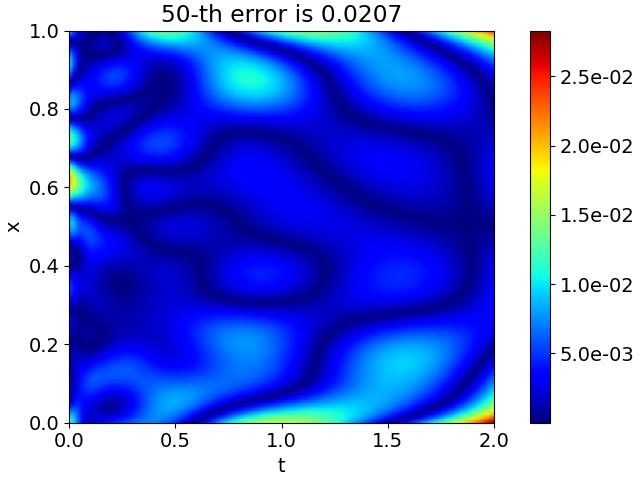}
        \end{subfigure}%
        \begin{subfigure}{.25\textwidth}
            \centering
            \includegraphics[height=0.75\textwidth,width=1.0\textwidth]{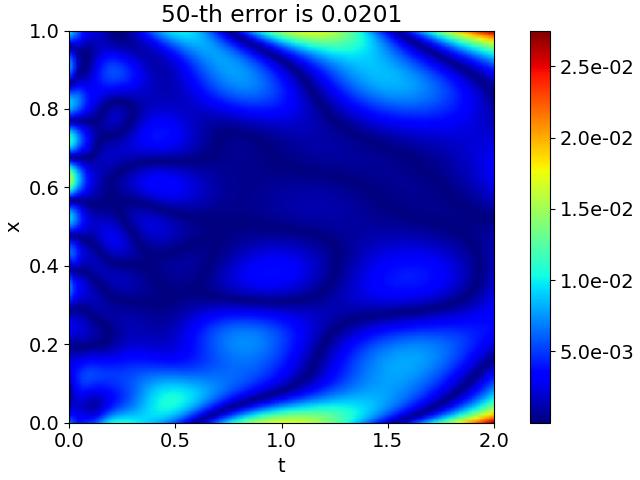}
        \end{subfigure}%
        \begin{subfigure}{.25\textwidth}
            \centering
            \includegraphics[height=0.75\textwidth,width=1.0\textwidth]{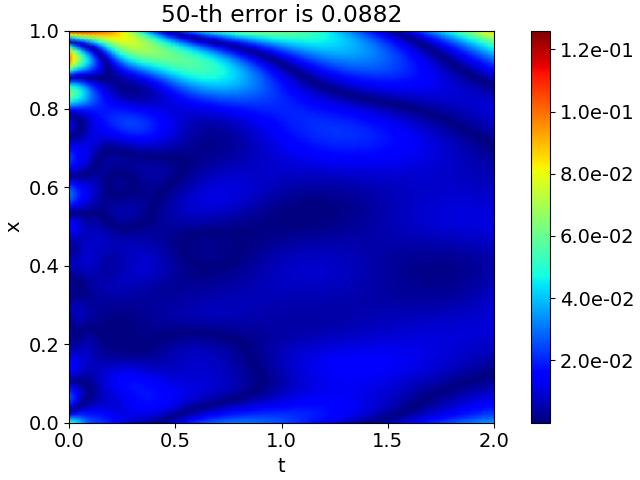}
        \end{subfigure}%
        \begin{subfigure}{.25\textwidth}
            \centering
            \includegraphics[height=0.75\textwidth,width=1.0\textwidth]{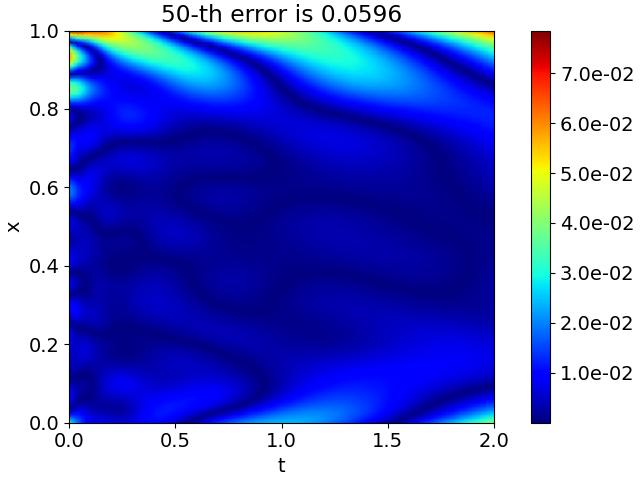}
        \end{subfigure}%
        \newline
        \raggedleft
        \begin{subfigure}{.25\textwidth}
            \centering
            \includegraphics[height=0.75\textwidth,width=1.0\textwidth]{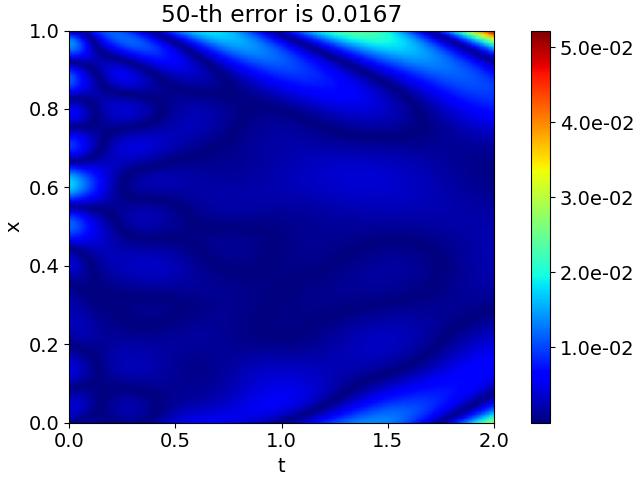}
        \end{subfigure}%
        \begin{subfigure}{.25\textwidth}
            \centering
            \includegraphics[height=0.75\textwidth,width=1.0\textwidth]{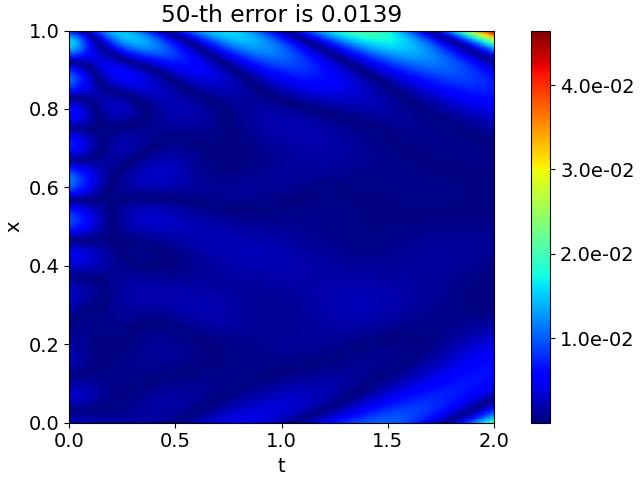}
        \end{subfigure}%
        \begin{subfigure}{.25\textwidth}
            \centering
            \includegraphics[height=0.75\textwidth,width=1.0\textwidth]{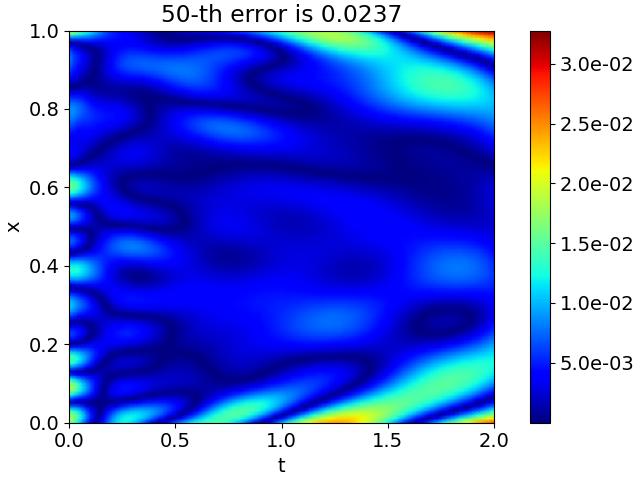}
        \end{subfigure}%
        \begin{subfigure}{.25\textwidth}
            \centering
            \includegraphics[height=0.75\textwidth,width=1.0\textwidth]{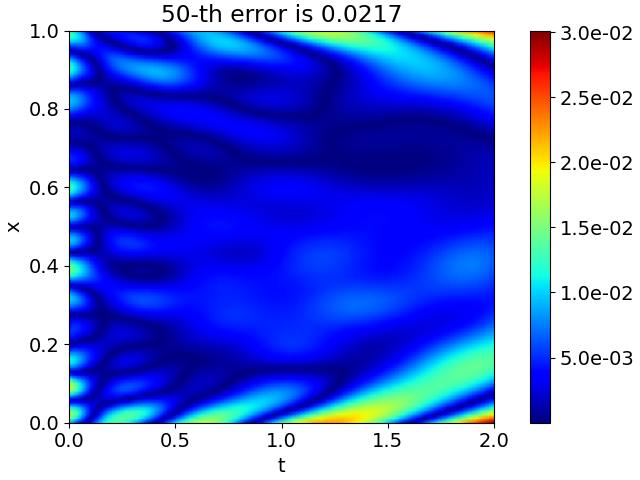}
        \end{subfigure}%
        \newline
        \raggedleft
        \begin{subfigure}{.25\textwidth}
            \centering
            \includegraphics[height=0.75\textwidth,width=1.0\textwidth]{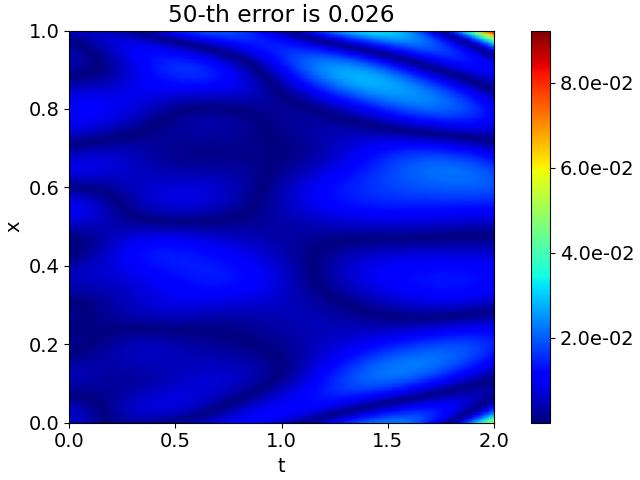}
            \caption{\textbf{MC-II}, $k=4$, $\lambda=4$.}
        \end{subfigure}%
        \begin{subfigure}{.25\textwidth}
            \centering
            \includegraphics[height=0.75\textwidth,width=1.0\textwidth]{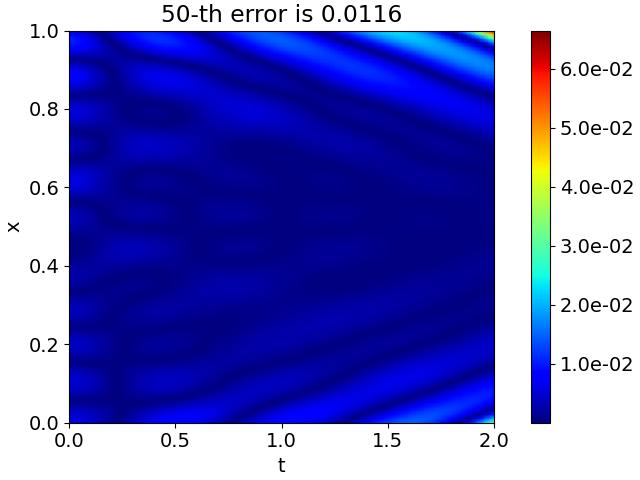}
            \caption{\textbf{GJ-II}, $k=4$, $\lambda=4$.}
        \end{subfigure}%
        \begin{subfigure}{.25\textwidth}
            \centering
            \includegraphics[height=0.75\textwidth,width=1.0\textwidth]{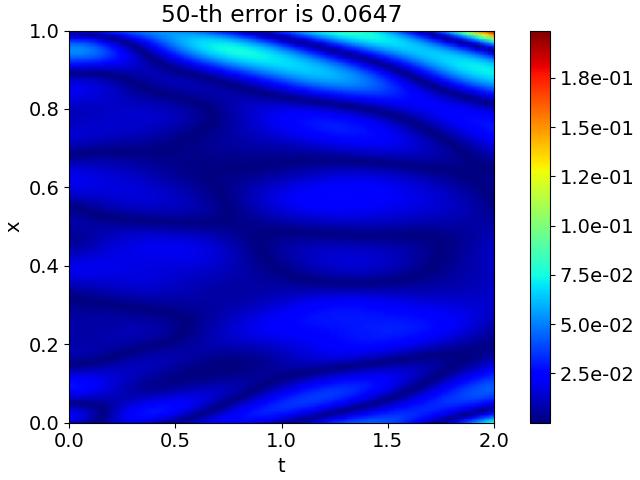}
            \caption{\textbf{MC-II}, $k=6$, $\lambda=6$.}
        \end{subfigure}%
        \begin{subfigure}{.25\textwidth}
            \centering
            \includegraphics[height=0.75\textwidth,width=1.0\textwidth]{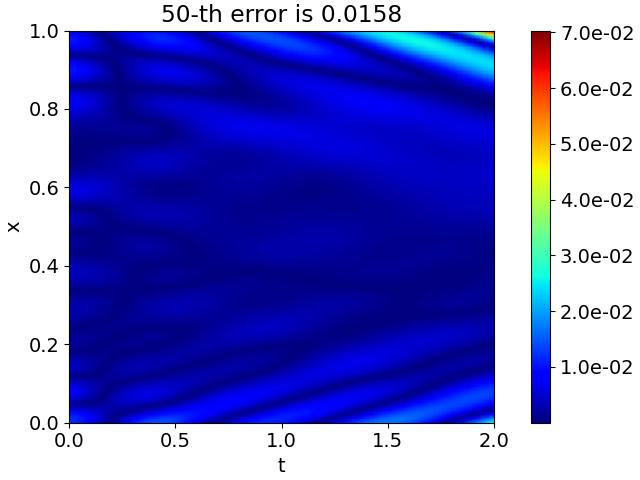}
            \caption{\textbf{GJ-II}, $k=6$, $\lambda=6$.}
        \end{subfigure}%
        \caption{Absolute errors of cases in Table \ref{tab:largeL} with sufficient trainings of different $\alpha$, $k$ and $\lambda$. For each sub-figure, in the first row, $\alpha=1.25$, in the second row, $\alpha=1.5$, in the last row $\al=1.75$.}
        \label{fig:err-largeL}
    \end{figure}

To further investigate the observed performance discrepancies among different configurations, we conducted extended training experiments with sufficient iterations. Specifically, for the cases presented in Table \ref{tab:eg2-1}, comprehensive training was implemented using the Adam optimizer with a learning rate of $10^{-5}$ and $L=50$. The corresponding numerical outcomes are documented in Table \ref{tab:largeL}. With sufficient training, both methods demonstrate enhanced solution accuracy when $\lambda$ is relatively small; however, as $\lambda$ increases, performance degradation becomes evident, highlighting the influence of the diffusion component on solution scaling. The Gauss-Jacobi method consistently outperforms the Monte Carlo approach under these extended training conditions. This suggests that for smaller M, the Monte Carlo integration lacks sufficient resolution to capture solution features effectively, leading to suboptimal training convergence. Consequently, when employing Monte Carlo-based integration within PINNs, a significantly larger number of quadrature samples is advisable to ensure acceptable accuracy. Figure \ref{fig:err-largeL} illustrates the absolute errors, revealing that the error patterns of \textbf{GJ-II} continue to present the high-frequency error structure, implying the better resolution of the Gauss-Jacobi method.
\begin{table}[htbp]
    \centering

    \begin{tabular}{|c|cc|cc|}
        \hline
                      & \multicolumn{2}{c|}{$k=1$, $\lambda=1$}                & \multicolumn{2}{c|}{$k=2$, $\lambda=4$}              \\ \hline
                      & \multicolumn{1}{c|}{\textbf{MC-II}}  & \textbf{GJ-II}  & \multicolumn{1}{c|}{\textbf{MC-II}} & \textbf{GJ-II} \\ \hline
        $\alpha=1.25$ & \multicolumn{1}{c|}{3.30e-3}          & 2.70e-3      & \multicolumn{1}{c|}{3.76e-2}        & 3.12e-2        \\ \hline
        $\alpha=1.5$ & \multicolumn{1}{c|}{3.20e-3}          & 1.60e-3         & \multicolumn{1}{c|}{2.65e-2}        & 1.66e-2        \\ \hline
        $\alpha=1.75$ & \multicolumn{1}{c|}{4.20e-3}          & 9.00e-4         & \multicolumn{1}{c|}{2.85e-2}        & 8.40e-3        \\ \hline
                      & \multicolumn{2}{c|}{$k=4$, $\lambda=4$} & \multicolumn{2}{c|}{$k=6$, $\lambda=6$}              \\ \hline
        $\alpha=1.25$ & \multicolumn{1}{c|}{2.07e-2}         & 2.01e-2         & \multicolumn{1}{c|}{8.82e-2}        & 5.96e-2        \\ \hline
        $\alpha=1.5$ & \multicolumn{1}{c|}{1.67e-2}          & 1.39e-2         & \multicolumn{1}{c|}{2.37e-2}        & 2.17e-2        \\ \hline
        $\alpha=1.75$ & \multicolumn{1}{c|}{2.60e-2}         & 1.16e-2        & \multicolumn{1}{c|}{6.47e-2}        & 1.58e-2        \\ \hline
        \end{tabular}
    \caption{Relative errors under sufficient training with $L=50$ and learning rate $10^{-5}$.}
    \label{tab:largeL}
    \end{table}

Extensive numerical experiments demonstrate that both parameters of the Mittag-Leffler function $E_{\alpha,\beta}(-\lambda t^\alpha)$ significantly influence the performance of PINNs. Specifically, $\lambda$ characterizes the decay properties of the solution, while $\alpha$ governs its wave-like characteristics. Larger values of $\lambda$ lead to more rapid solution decay, which can degrade the fidelity of the learned solution. Conversely, increasing $\alpha$ tends to counteract this decay effect by enhancing the wave component, thereby improving the overall solution quality. This behavior is reflected in our numerical results, where models trained with $\alpha = 1.75$ consistently outperform those with $\alpha = 1.25$, particularly in scenarios involving large $\lambda$. In this investigation, we focus on the effects of $\lambda$ and $\alpha$ while maintaining a fixed final time $T=2$, thereby isolating these parameters from temporal influences.

As both $\lambda$ and $\alpha$ increase, the complexity of the solution dynamics imposes greater challenges on the training process. Notably, as $\alpha$ approaches 2, the performance of the Monte Carlo-based method degrades significantly. This degradation is attributed to sensitivity to the hyperparameter $\varepsilon$, as summarized in Section~\ref{subsec:num-valid} and detailed in Appendix~\ref{app:num-valid}. While increases in $\lambda$ negatively affect the accuracy of all tested methods, the Gauss-Jacobi approach demonstrates superior robustness, particularly for larger values of $\alpha$, where the mitigated decay yields more stable solution representations.

These findings highlight the necessity of selecting an appropriate quadrature strategy, either Gauss-Jacobi or Monte Carlo, based on the specific parameter regimes characterizing the problem. The choice of numerical integration method plays a critical role in the efficiency and accuracy of PINN training, especially in the context of time-fractional partial differential equations.

\begin{figure}[htbp]
    \centering
    \includegraphics[scale=0.5]{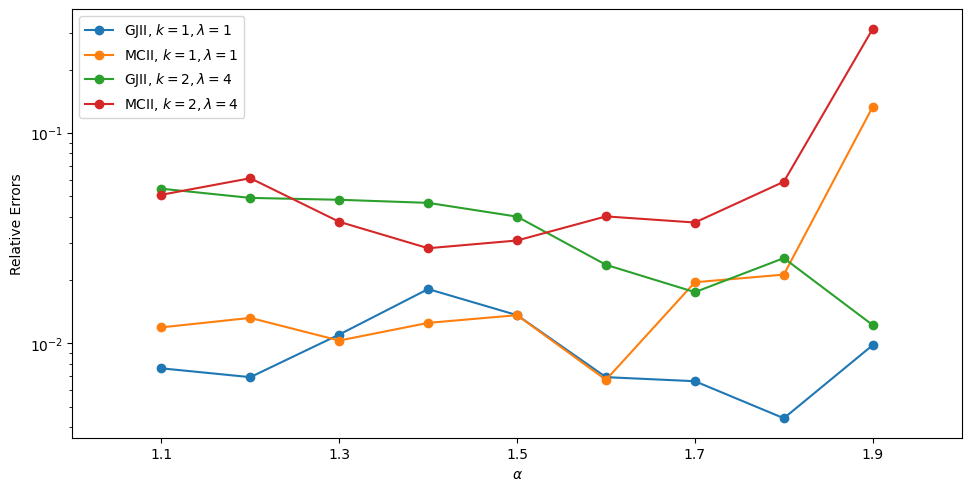}
    \caption{Relative errors of \textbf{MC-II} and \textbf{GJ-II} with different $\alpha$ for $k=1, \lambda=1$ and $k=2$, $\lambda=4$.}
    \label{fig:al-vary}
\end{figure}

A series of experiments were performed for different values of $\alpha$. Figure \ref{fig:al-vary} presents the numerical errors for $\alpha\in[1.1,1.9]$ using both Monte Carlo and Gauss-Jacobi methods under various $k$ and $\lambda$ settings. As $\lambda$ increases, the numerical performance deteriorates across all methods. As $\al$ approaches 2, the performance of the Monte Carlo method deteriorates primarily due to the hyperparameter setting $\varepsilon$.

\begin{figure}[htbp]
    \centering
    \includegraphics[scale=0.5]{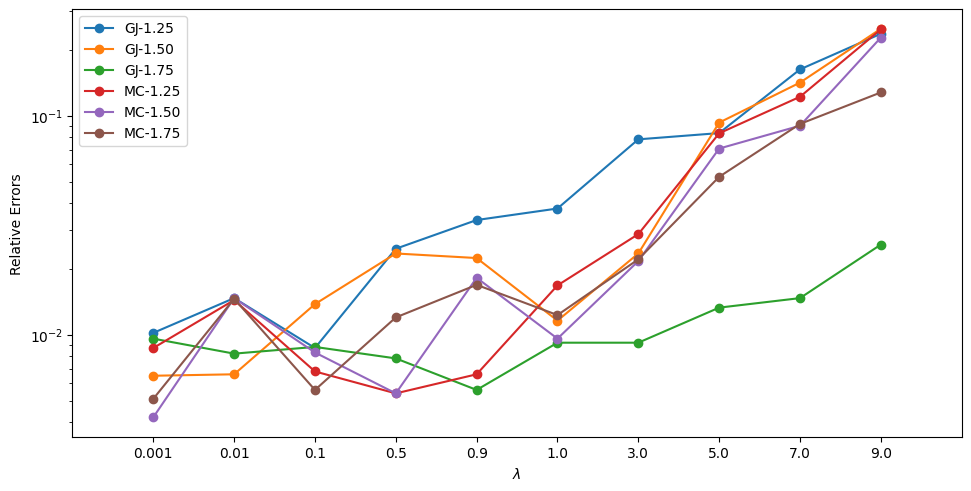}
    \caption{Relative errors of \textbf{MC-II} and \textbf{GJ-II} with different $\lambda$ for $k=2$, $\alpha\in \{1.25,1.5,1.75\}$.}
    \label{fig:lam-vary}
\end{figure}

Furthermore, we conducted additional experiments with a fixed $k=2$ and varying decay parameter $\lambda$ for $\alpha \in {1.25, 1.5, 1.75}$. The results are illustrated in Figure \ref{fig:lam-vary}. These experiments reveal that increasing the decay parameter $\lambda$ substantially impacts the numerical performance of both integration methods. For the special case of \textbf{GJ-II} with $\al=1.75$, with the slow decay of the exact solution, PINNs could still simulate the solution satisfactorily.

\begin{figure}[htbp]
    \centering
    \includegraphics[scale=0.5]{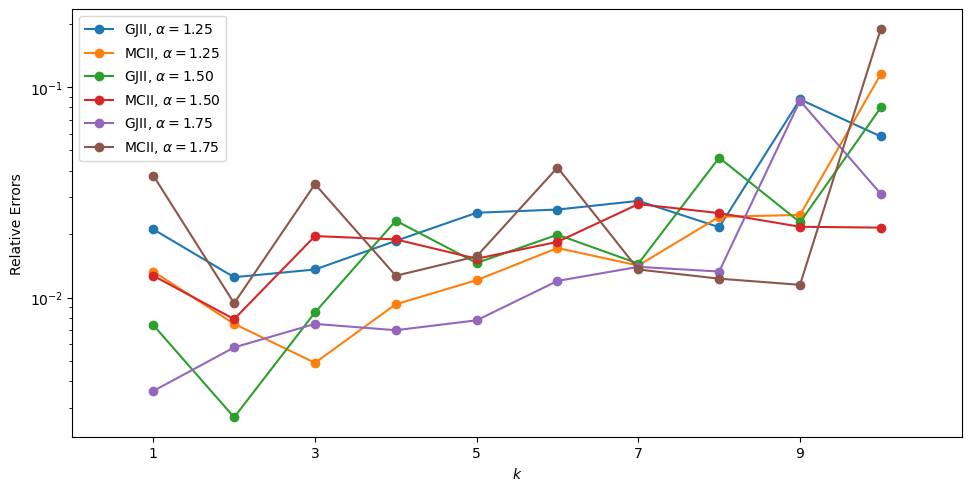}
    \caption{Relative errors of \textbf{MC-II} and \textbf{GJ-II} with different $k$ for $\lambda=1$, $\alpha\in \{1.25,1.5,1.75\}$.}
    \label{fig:k-vary}
\end{figure}

Additionally, we examine the influence of parameter $k$ on solution behavior with fixed $\lambda=1.0$ to isolate it from the effects of large $\lambda$ values on numerical performance. The results, presented in Figure \ref{fig:k-vary}, indicate that as $k$ increases substantially, solution accuracy slightly deteriorates, suggesting the need for enhanced training strategies when dealing with large $k$ values.

In general, the Mittag-Leffler function embedded within the solution of the time-fractional diffusion-wave equation presents significant training challenges. The parameter $\lambda$, which typically originates from the eigenvalue component of the spatial operator and characterizes the decay behavior of the solution, exerts a relatively major influence on the solution behavior for both Monte Carlo and Gauss-Jacobi methods due to the fast decay of the solution. Conversely, the parameter $\alpha$, which derives from the time-fractional component of the temporal differential operator and governs the wave characteristics of the solution, substantially influences the solution behavior for the Monte Carlo method due to hyperparameter $\varepsilon$ when $\al$ approaches 2, but not for the Gauss-Jacobi method. Large $\al$ enhances the wave characteristics of the solution and delays the solution decay, thus facilitating the PINN simulation of the solution.

It is important to highlight that the main challenges in training time-fractional diffusion-wave equations are closely related to the broader issue of training PINNs for wave equations. For future work, relevant insights can be drawn from the published results in \cite{WangYuParis:2022:NTK,Wang2022:causality}.

\subsection{Time-fractional Burgers equation with $\alpha \in (1, 2)$}
\label{subsec:burgers}
In this subsection we consider the time-fractional  Burgers equation with $\alpha \in (1, 2)$, given by
\begin{equation}
    \label{eqn:burgers}
    \begin{aligned}
    \dta u+uu_x &= \frac{0.01}{\pi}u_{xx},\text{ in }\Omega\times(0,T],\\
    u &=0, \text{ on }\partial\Omega,\\
    u(0) = &-\sin(\pi x),\, \partial_t u(0) = \sin(\pi x), \text{ in }\Omega,
    \end{aligned}
\end{equation}
where $\Omega=(-1,1)$, $T=1$.

\begin{figure}[htbp]
    \centering
    \begin{subfigure}{.5\textwidth}
        \centering
        \includegraphics[scale=0.35]{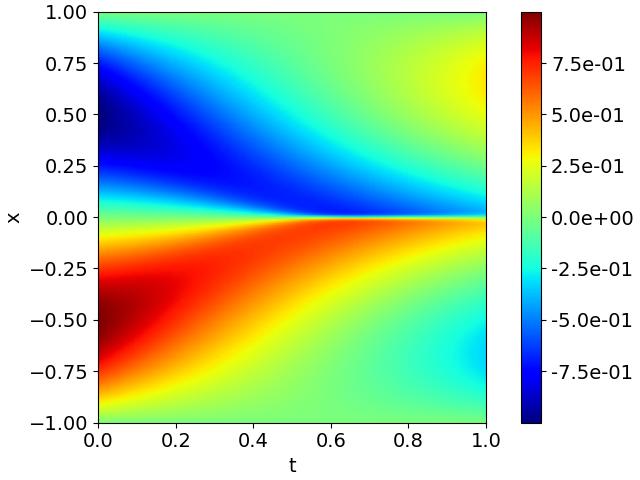}
        \caption{$\al=1.1$.}
        \end{subfigure}%
        \begin{subfigure}{.5\textwidth}
        \centering
        \includegraphics[scale=0.35]{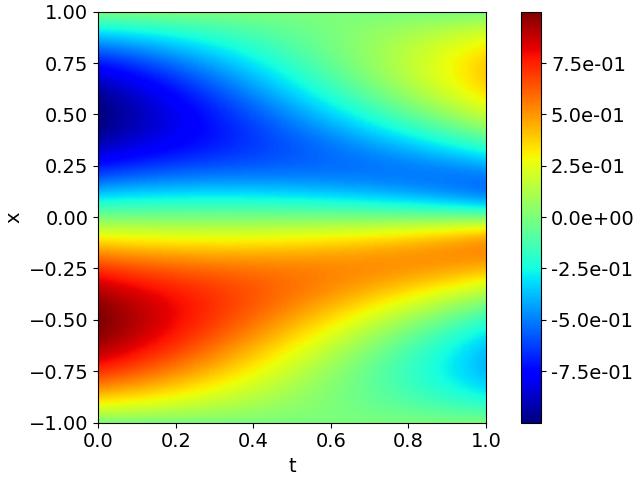}
        \caption{$\al=1.8$.}
        \end{subfigure}%
        \caption{Exact solutions of Burgers equation \eqref{eqn:burgers} with $\al=1.1$ and $\al=1.8$.}
    \label{fig:burgers-exact}
\end{figure}

Building upon the analysis presented in Section \ref{subsec:eg2}, we observe that the wave component establishes a dynamic equilibrium with the diffusion component of the solution. This phenomenon is particularly evident in the fractional Burgers equation when examining different $\alpha$. Figure \ref{fig:burgers-exact} illustrates the exact solution profiles for two representative cases: $\alpha=1.1$ and $\alpha=1.8$. The results clearly demonstrate that shock formation is progressively delayed as $\alpha$ increases, providing further evidence of the complex interplay between diffusive and wave-like behaviors in diffusion-wave systems.

In the following experiments, we employ a minimal dataset comprising 1500 collocation points within the domain and 500 points along the spatial and temporal boundaries. The Adam optimizer is employed with a learning rate of $0.0001$, training for 10,000 epochs with a maximum of $L=20$ iterations. The initial loss weight is set to $\lambda_{init}=100$. We implement the RAD method for adaptive sampling as proposed by \cite{LuLu:2023:RAD}. We set the adaptive resampling ratio to $0.3$, e.g., in each iteration we substitute $30\%$ of the points with RAD sampling points in the domain. Our evaluation of the RAD application focuses on three key metrics: relative errors, maximum absolute errors, and the frequency distribution of absolute error patterns.

Throughout these experiments, we utilize our transformed integration scheme as the default computational approach with $M=80$, based on our findings that these methods can achieve comparable accuracy with significantly reduced computational costs. Based on the findings presented in Section \ref{subsec:eg2}, we demonstrate a clear preference for the Gauss-Jacobi method (\textbf{GJ-II}) over the  Monte Carlo method (\textbf{MC-II}).

\begin{figure}[htbp]
    \centering
    \begin{subfigure}{.25\textwidth}
        \centering
        \includegraphics[height=0.75\textwidth,width=1.0\textwidth]{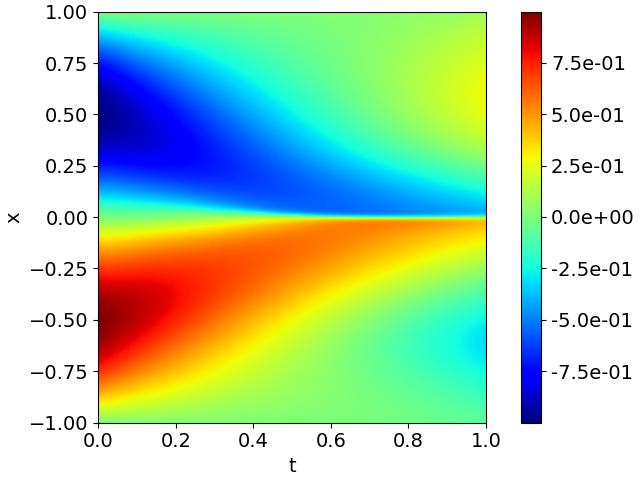}
    \end{subfigure}%
    \begin{subfigure}{.25\textwidth}
        \centering
        \includegraphics[height=0.75\textwidth,width=1.0\textwidth]{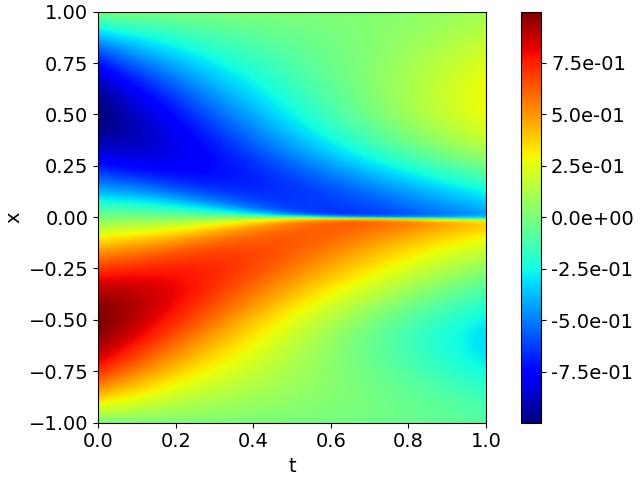}
    \end{subfigure}%
    \begin{subfigure}{.25\textwidth}
        \centering
        \includegraphics[height=0.75\textwidth,width=1.0\textwidth]{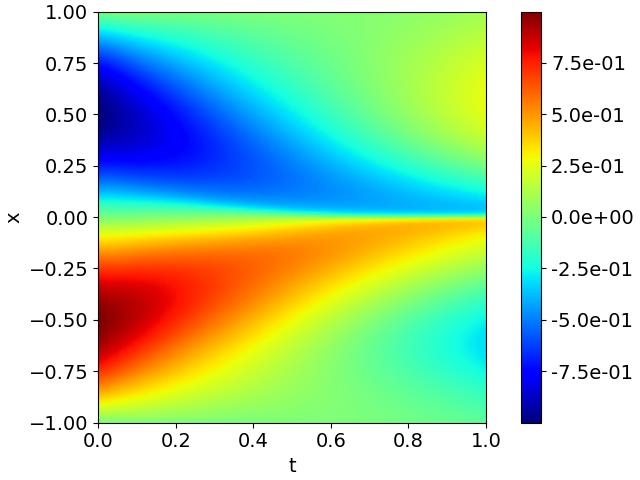}
    \end{subfigure}%
    \begin{subfigure}{.25\textwidth}
        \centering
        \includegraphics[height=0.75\textwidth,width=1.0\textwidth]{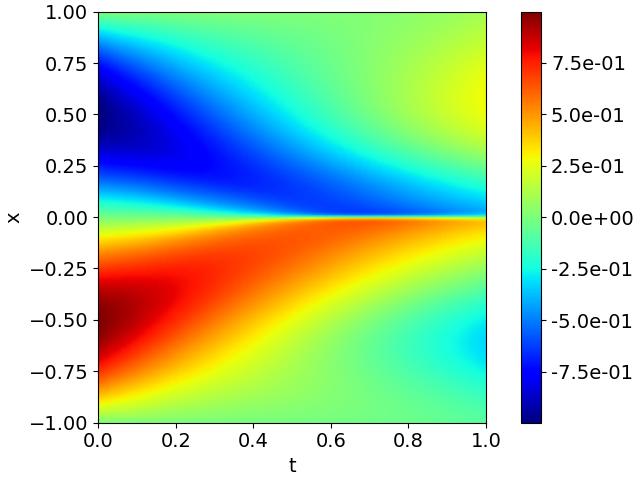}
    \end{subfigure}%
    \newline
    \raggedleft
    \begin{subfigure}{.25\textwidth}
        \centering
        \includegraphics[height=0.75\textwidth,width=1.0\textwidth]{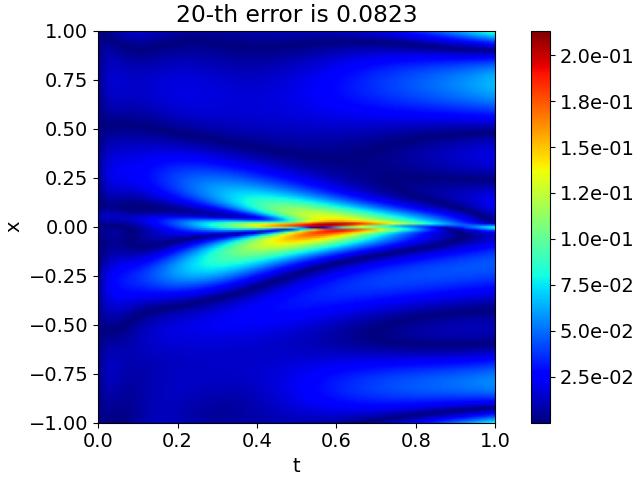}
        \caption{\textbf{MC-II} without RAD.}
    \end{subfigure}%
    \begin{subfigure}{.25\textwidth}
        \centering
        \includegraphics[height=0.75\textwidth,width=1.0\textwidth]{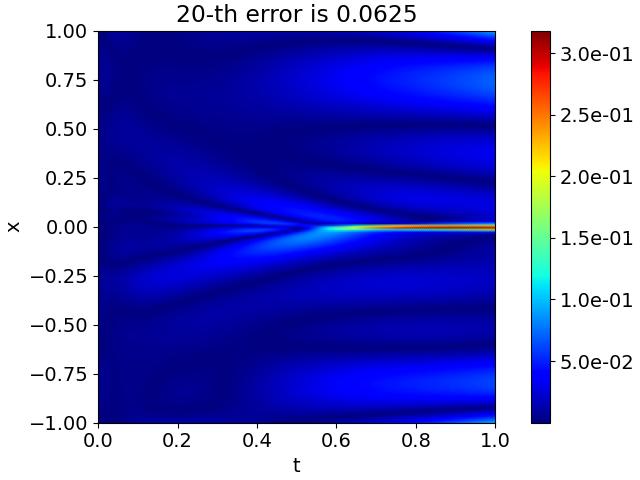}
        \caption{\textbf{MC-II} with RAD.}
    \end{subfigure}%
    \begin{subfigure}{.25\textwidth}
        \centering
        \includegraphics[height=0.75\textwidth,width=1.0\textwidth]{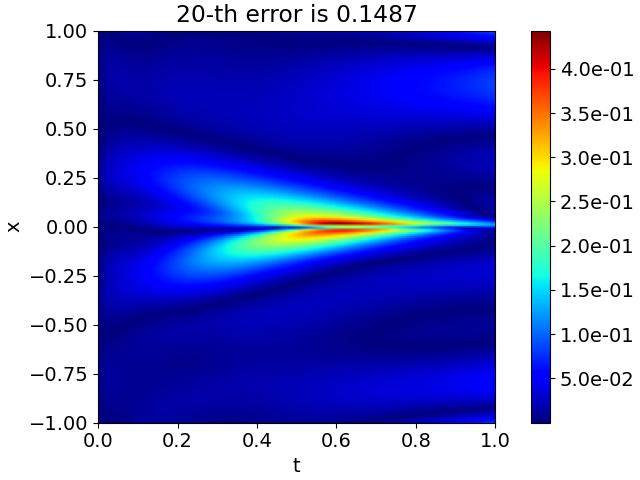}
        \caption{\textbf{GJ-II} without RAD.}
    \end{subfigure}%
    \begin{subfigure}{.25\textwidth}
        \centering
        \includegraphics[height=0.75\textwidth,width=1.0\textwidth]{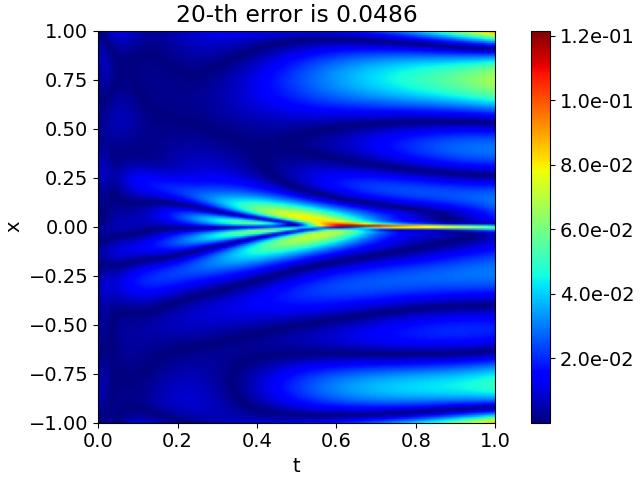}
        \caption{\textbf{GJ-II} with RAD.}
    \end{subfigure}%
    \caption{Profiles of absolute error and neural network solutions for Burgers equation \eqref{eqn:burgers} at $\al=1.1$. First row: numerical solutions. Second row: absolute error.}
    \label{fig:burgers-11}
\end{figure}

Figure \ref{fig:burgers-11} presents the numerical solutions and associated error distributions for the case where $\alpha=1.1$, comparing different numerical integration methods. The results demonstrate that the RAD method achieves significantly higher accuracy in solution approximation compared to the uniform sampling method.  This behavior indicates that RAD helps the optimization process escape local minima during training.

\begin{figure}[htbp]
    \centering
    \begin{subfigure}{.25\textwidth}
        \centering
        \includegraphics[height=0.75\textwidth,width=1.0\textwidth]{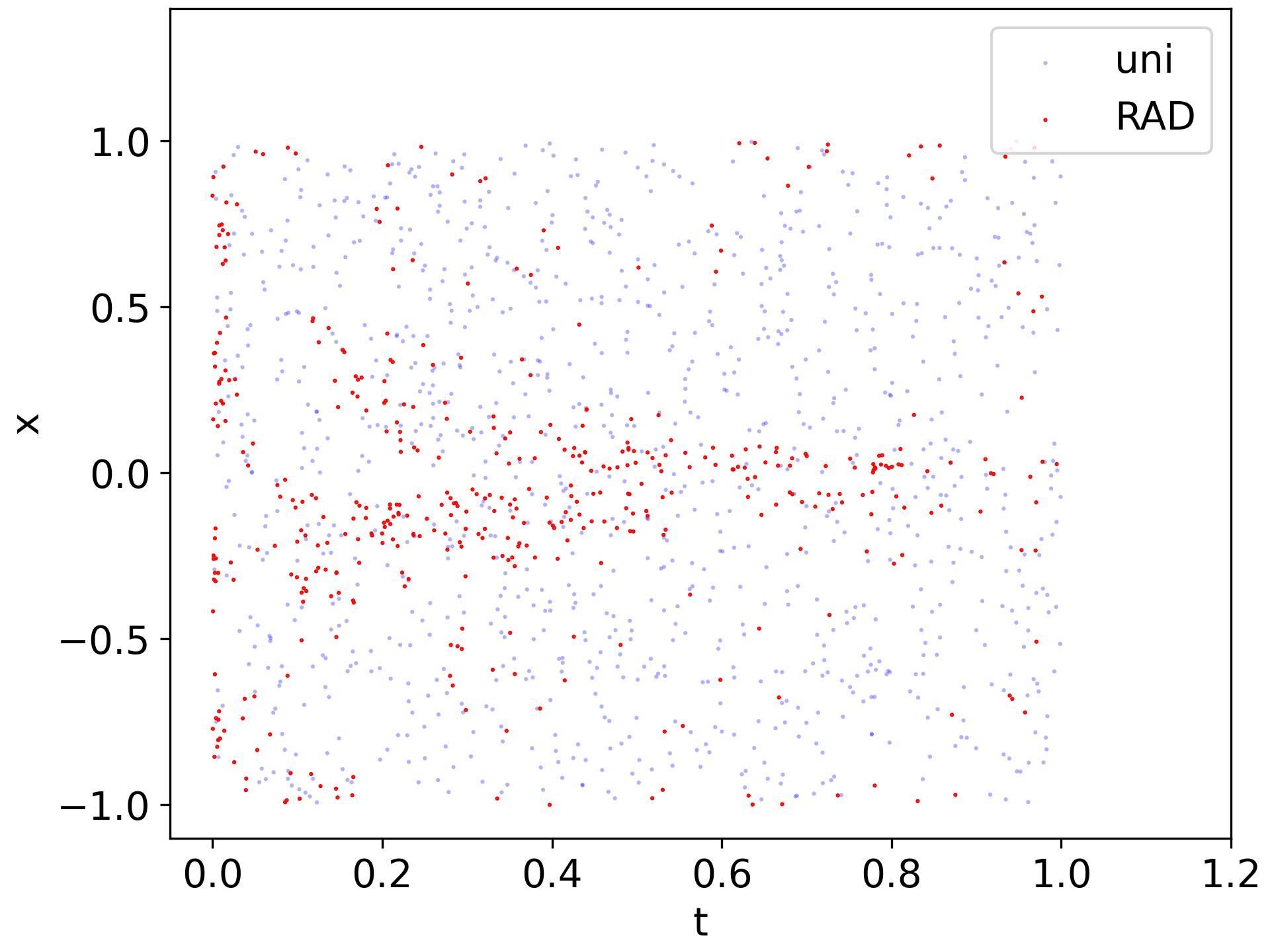}
        \caption{\textbf{MC-II}, after 3rd iteration.}
    \end{subfigure}%
    \begin{subfigure}{.25\textwidth}
        \centering
        \includegraphics[height=0.75\textwidth,width=1.0\textwidth]{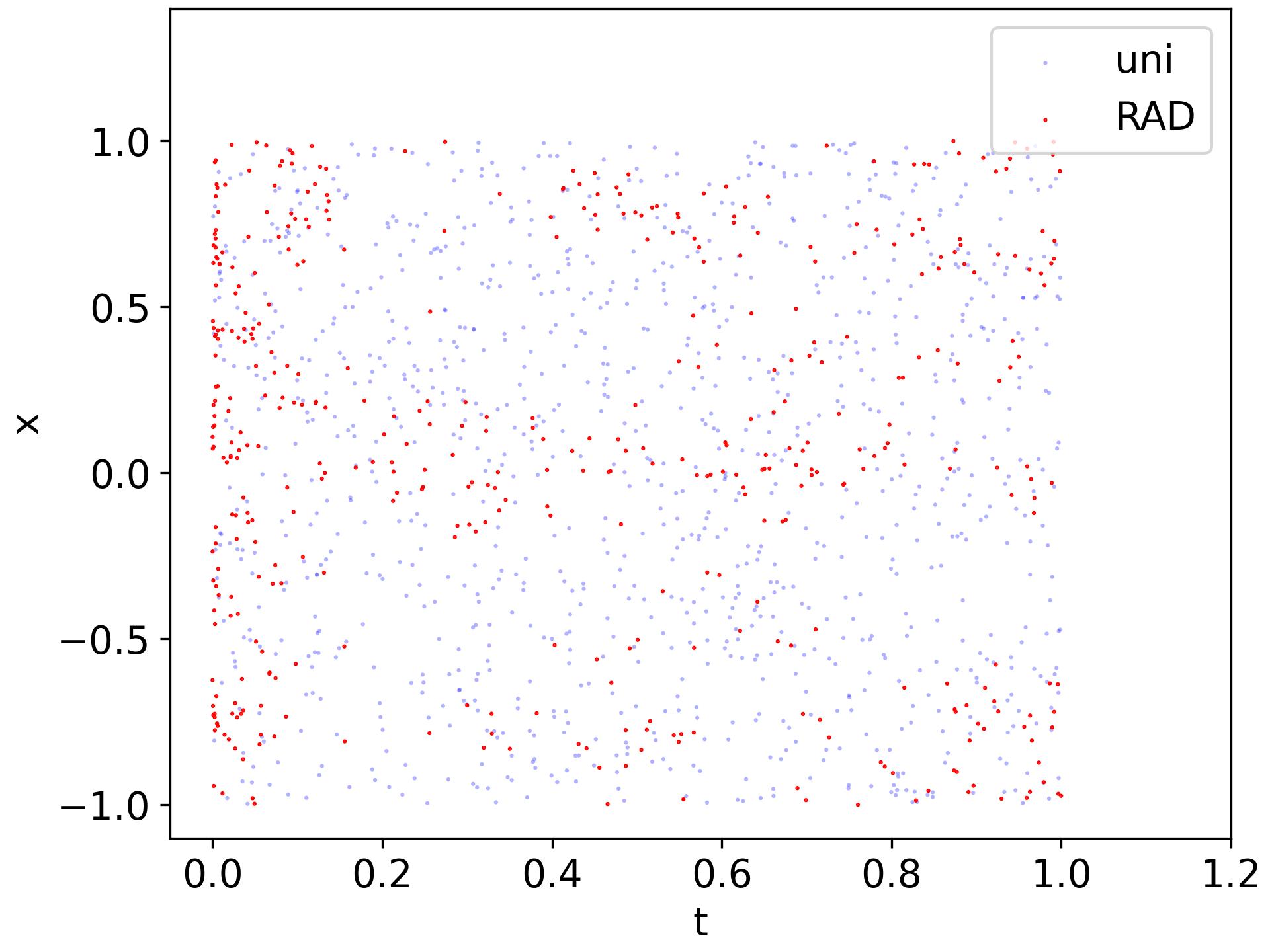}
        \caption{\textbf{MC-II}, after 15th iteration.}
    \end{subfigure}%
    \begin{subfigure}{.25\textwidth}
        \centering
        \includegraphics[height=0.75\textwidth,width=1.0\textwidth]{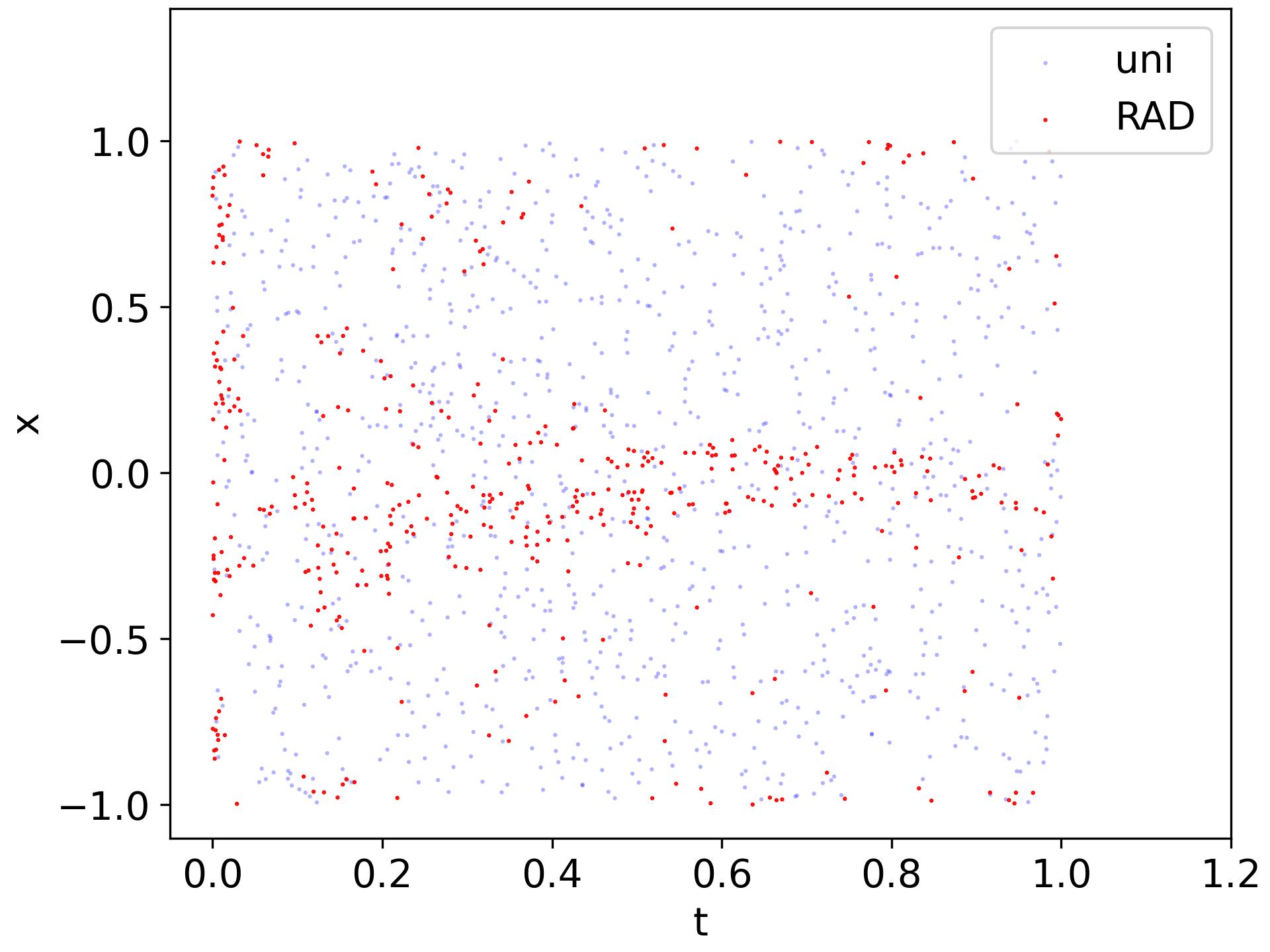}
        \caption{\textbf{GJ-II}, after 3rd iteration.}
    \end{subfigure}%
    \begin{subfigure}{.25\textwidth}
        \centering
        \includegraphics[height=0.75\textwidth,width=1.0\textwidth]{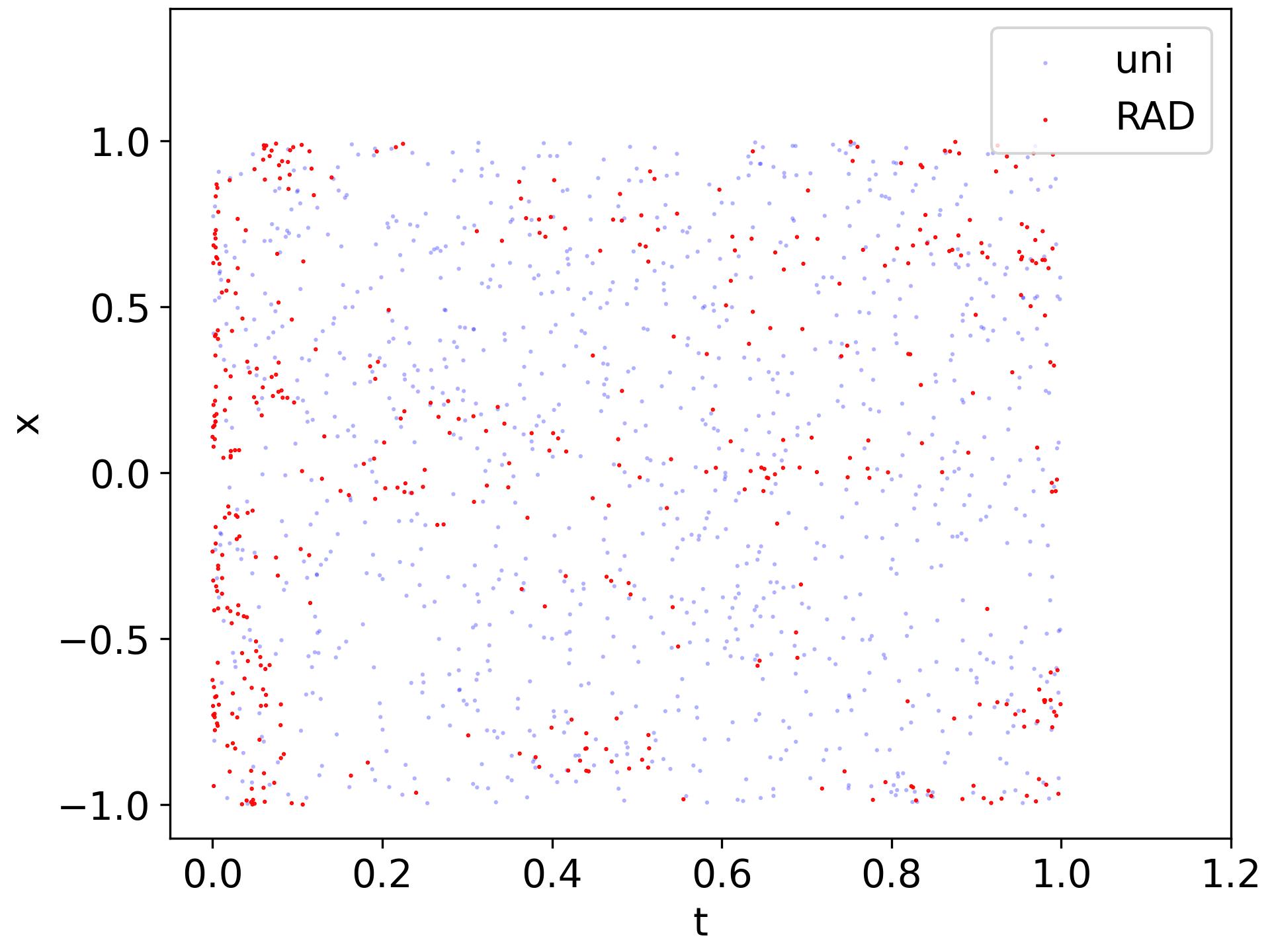}
        \caption{\textbf{GJ-II}, after 15th iteration.}
    \end{subfigure}%
    \caption{Profiles of node distributions for Burgers equation \eqref{eqn:burgers} with $\al=1.1$.}
    \label{fig:burgers-11-node}
\end{figure}

Figure \ref{fig:burgers-11-node} illustrates the distribution of collocation points sampled by the RAD method at different training stages. Our analysis reveals a distinct pattern in the sampling distribution: during the initial training phases, the loss function predominantly concentrates on regions near the temporal boundary ($t=0$) and shock formation zones. As training progresses, however, the spectral characteristics of the loss function exhibit increased frequency components, causing the RAD method to redistribute its sampling focus across a broader domain.

Notably, this behavior differs significantly from RAD's performance in integer-order Burgers equations (see details in \cite{Zhang2025}). In the fractional-order case, even when the adaptive sampling correctly identifies shock regions, the PINNs struggle to achieve satisfactory training outcomes due to the inherent wave component in the solution structure. This observation suggests that effective training strategies for fractional-order Burgers equations should prioritize accurately capturing the wave characteristics of the solution before addressing the shock regions that arise from the nonlinear advection term in the Burgers equation.

\begin{figure}[htbp]
    \centering
    \begin{subfigure}{.25\textwidth}
        \centering
        \includegraphics[height=0.75\textwidth,width=1.0\textwidth]{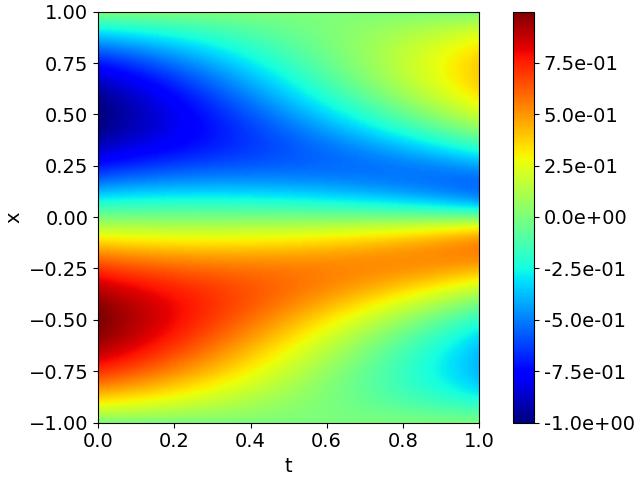}
    \end{subfigure}%
    \begin{subfigure}{.25\textwidth}
        \centering
        \includegraphics[height=0.75\textwidth,width=1.0\textwidth]{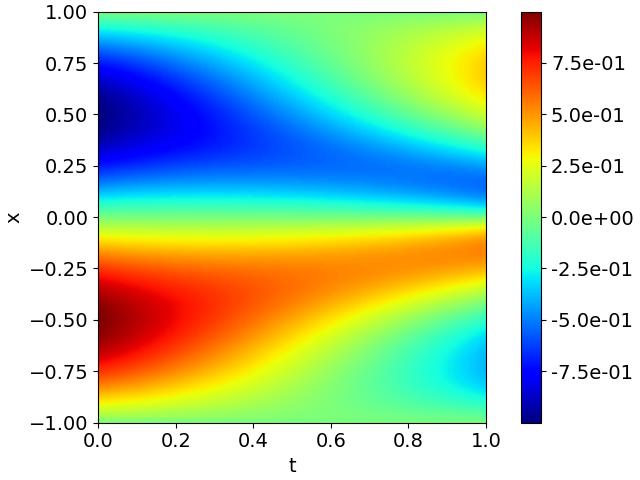}
    \end{subfigure}%
    \begin{subfigure}{.25\textwidth}
        \centering
        \includegraphics[height=0.75\textwidth,width=1.0\textwidth]{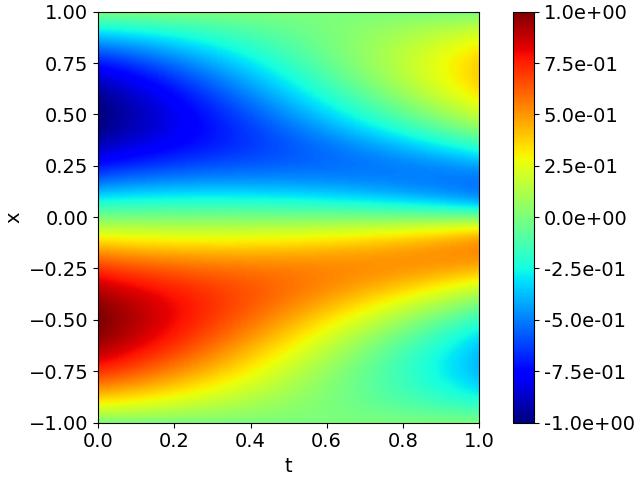}
    \end{subfigure}%
    \begin{subfigure}{.25\textwidth}
        \centering
        \includegraphics[height=0.75\textwidth,width=1.0\textwidth]{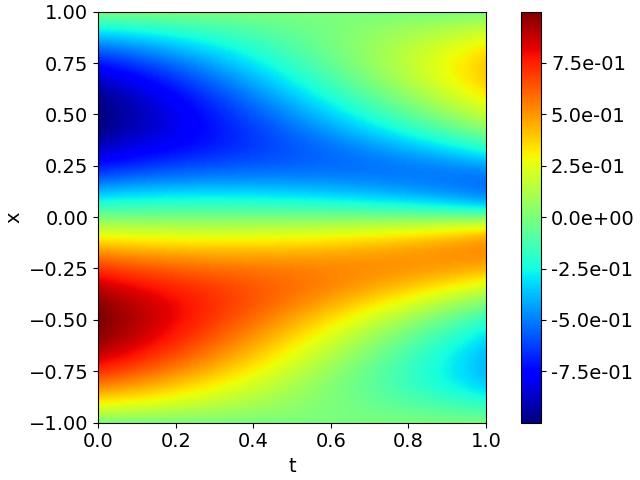}
    \end{subfigure}%
    \newline
    \raggedleft
    \begin{subfigure}{.25\textwidth}
        \centering
        \includegraphics[height=0.75\textwidth,width=1.0\textwidth]{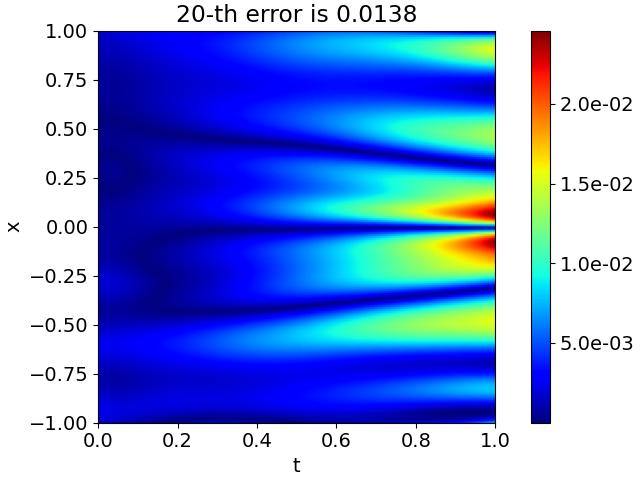}
        \caption{\textbf{MC-II} without RAD.}
    \end{subfigure}%
    \begin{subfigure}{.25\textwidth}
        \centering
        \includegraphics[height=0.75\textwidth,width=1.0\textwidth]{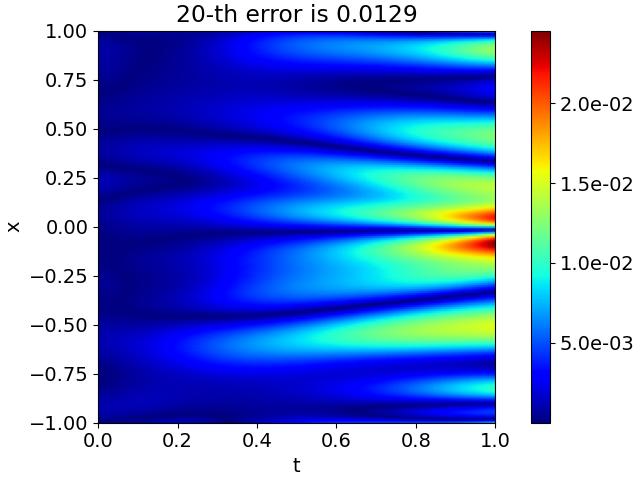}
        \caption{\textbf{MC-II} with RAD.}
    \end{subfigure}%
    \begin{subfigure}{.25\textwidth}
        \centering
        \includegraphics[height=0.75\textwidth,width=1.0\textwidth]{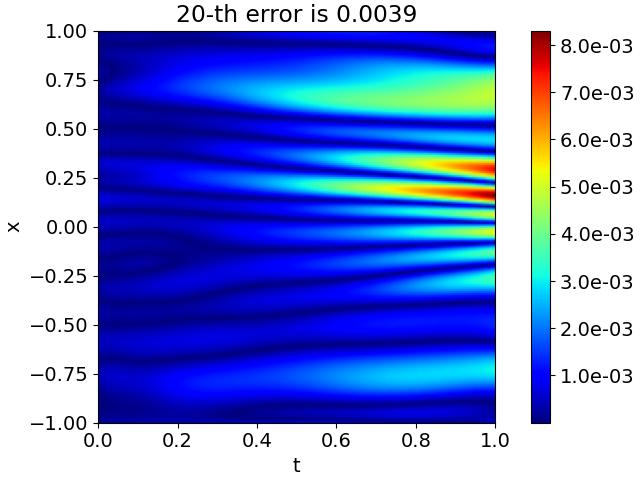}
        \caption{\textbf{GJ-II} without RAD.}
    \end{subfigure}%
    \begin{subfigure}{.25\textwidth}
        \centering
        \includegraphics[height=0.75\textwidth,width=1.0\textwidth]{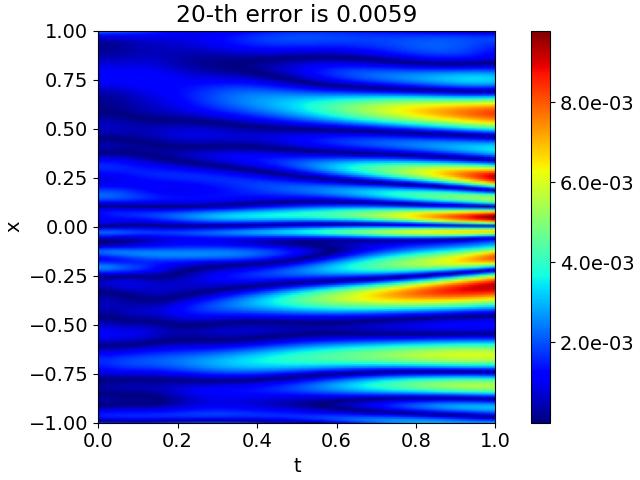}
        \caption{\textbf{GJ-II} with RAD.}
    \end{subfigure}%
    \caption{Profiles of absolute error and neural network solutions for Burgers equation \eqref{eqn:burgers} at $\al=1.8$. First row: numerical solutions. Second row: absolute error.}
    \label{fig:burgers-18}
\end{figure}

For the case where $\alpha=1.8$, Figure \ref{fig:burgers-18} demonstrates that PINNs generate satisfactory numerical approximations, primarily because shock formation does not occur before the terminal time $T=1$. Consequently, the RAD method provides negligible improvement to the solution behaviors. However, we observe that the Gauss-Jacobi method (\textbf{GJ-II}) significantly outperforms the Monte Carlo method (\textbf{MC-II}). This performance gap is attributed to the dominance of wave characteristics in the PDE dynamics when $\alpha > 1.5$, a regime where the RAD method proves less effective. Furthermore, as evidenced by Figure \ref{fig:diff_al}, the superiority of the Gauss-Jacobi method over Monte Carlo integration becomes increasingly pronounced as $\alpha$ approaches $2.0$. For potential improvements, one may refer to adaptive weighting strategies proposed in \cite{WangParis:2021:SISC,Wang2022:causality}.

\begin{figure}[htbp]
    \centering
    \begin{subfigure}{.25\textwidth}
        \centering
        \includegraphics[height=0.75\textwidth,width=1.0\textwidth]{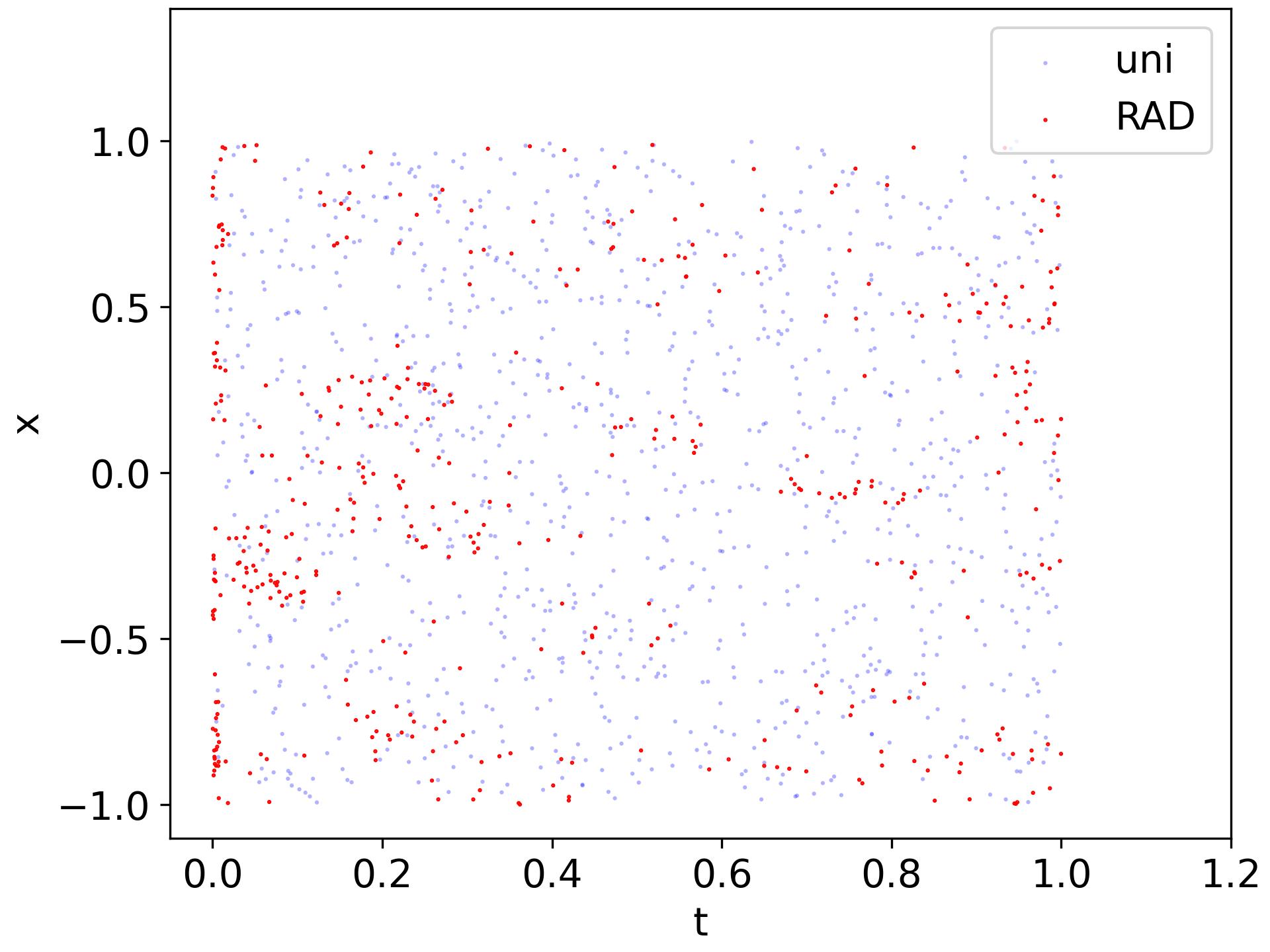}
        \caption{\textbf{MC-II}, after 2nd iteration.}
    \end{subfigure}%
    \begin{subfigure}{.25\textwidth}
        \centering
        \includegraphics[height=0.75\textwidth,width=1.0\textwidth]{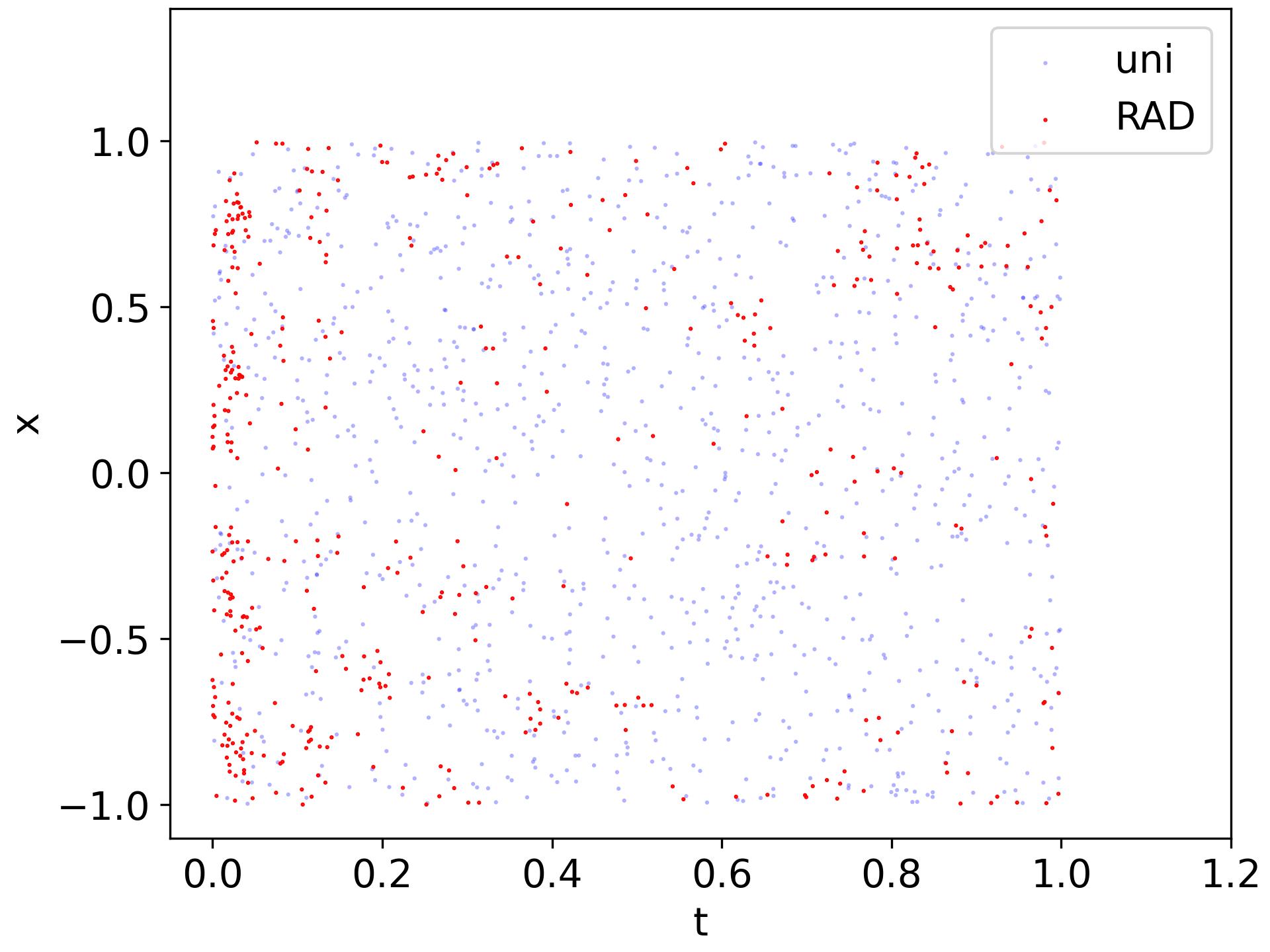}
        \caption{\textbf{MC-II}, after 20th iteration.}
    \end{subfigure}%
    \begin{subfigure}{.25\textwidth}
        \centering
        \includegraphics[height=0.75\textwidth,width=1.0\textwidth]{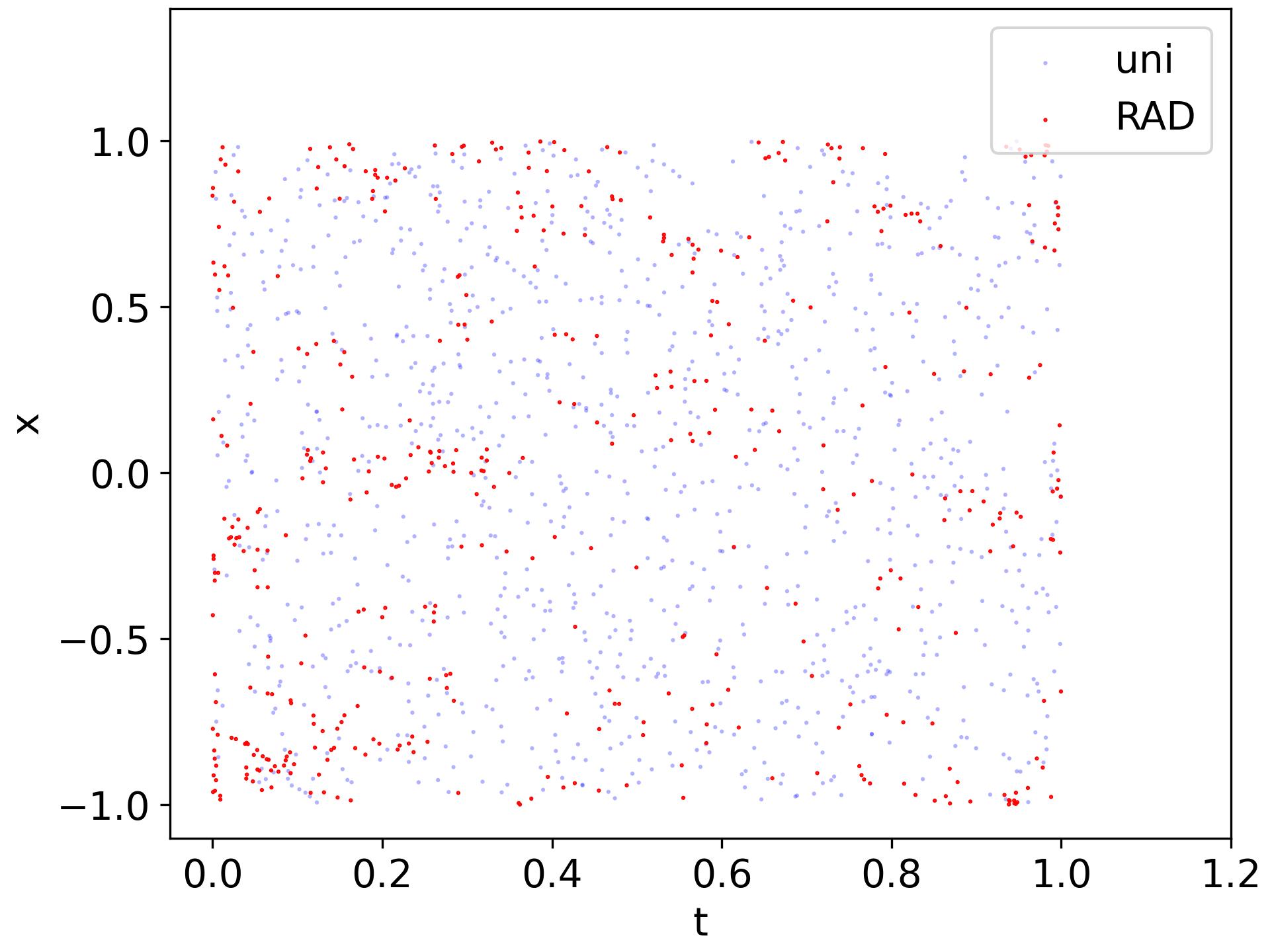}
        \caption{\textbf{GJ-II}, after 2nd iteration.}
    \end{subfigure}%
    \begin{subfigure}{.25\textwidth}
        \centering
        \includegraphics[height=0.75\textwidth,width=1.0\textwidth]{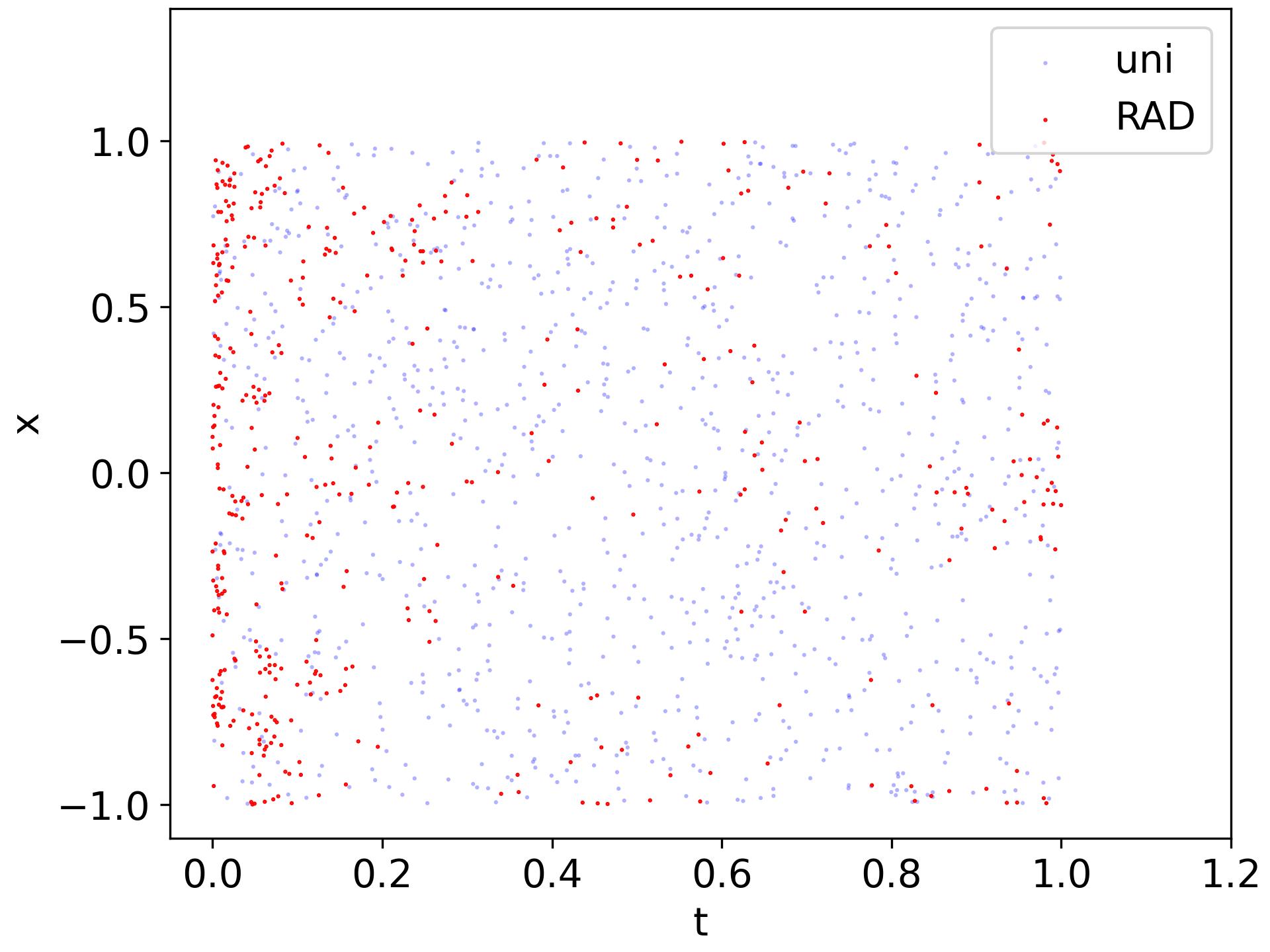}
        \caption{\textbf{GJ-II}, after 20th iteration.}
    \end{subfigure}%
    \caption{Profiles of node distributions for Burgers equation \eqref{eqn:burgers} with $\al=1.8$.}
    \label{fig:burgers-18-node}
\end{figure}

The distinct numerical behaviors observed when applying the RAD to Monte Carlo and Gauss-Jacobi integration methods can be directly attributed to the fundamental influence of adaptive sampling on the optimization landscape. Figure \ref{fig:burgers-18-node} illustrates the distribution of sampling points at both early and final training stages. For both methods, RAD effectively concentrates computational resources in critical regions near the temporal boundaries ($t=0$ and $t=T$), where the PINNs must accurately approximate the initial condition near $t=0$ and anticipate potential shock formation near $t=T$.

The numerical behaviors observed across these experiments align with our expectations. The RAD methodology demonstrates significant efficacy in addressing problems characterized by solution singularities, effectively redistributing computational datasets to critical regions. However, while the Gauss-Jacobi method generally outperforms the Monte Carlo approach, PINNs still face challenges in focusing on shock regions when the wave component in the solution structure is significant. This highlights the need for specialized techniques when applying adaptive sampling frameworks to time-fractional Burgers equations with $\alpha \in (1, 2)$.

These findings lead to an important conclusion regarding numerical approaches for time-fractional derivatives: mesh-free numerical integration techniques offer substantially greater flexibility than their mesh-based counterparts for specific problem classes. The capacity to incorporate adaptive sampling strategies represents a significant advantage of these mesh-free approaches, enabling efficient resolution of complex solution features without the constraints imposed by fixed discretization schemes. However, the unique properties of time-fractional diffusion-wave equations, which exhibit both diffusive and wave-like behaviors, require careful consideration when applying some advanced training strategies like adaptive sampling methods and adaptive weighting methods (see more details in \cite{Wang2022:causality,WangYuParis:2022:NTK}).

\captionsetup{font={color=blue}}
\subsection{A fractional ODE test for mesh-free adaptive sampling}
\label{subsec:ode-rad}
The Burgers experiment indicates that RAD does not always lead to a clear improvement for diffusion-wave dynamics: even when adaptive sampling identifies important regions, the nonlocal memory and wave-like components may still dominate the training difficulty. To isolate the benefit of mesh-free quadrature under adaptive sampling, we next consider a simple fractional ODE,
\begin{equation}
    \label{eqn:ode-rad}
    \begin{aligned}
        \dta u(t) &= f_\lambda(t), \quad t\in(0,T],\\
        u(0)&=0,\quad \partial_t u(0)=0,
    \end{aligned}
\end{equation}
where
\[
    u(t)=1-E_{\alpha,1}(-\lambda t^\alpha),\qquad
    f_\lambda(t)=\lambda E_{\alpha,1}(-\lambda t^\alpha).
\]
The test uses $\alpha=1.5$, $\lambda=100$, and $T=1$. Both methods are trained for 10000 Adam steps with the same DNN backbone of 7 hidden layers and 20 neurons per hidden layer. The fPINNs baseline uses a fixed temporal mesh with $\tau=0.01$, giving 100 time points. The \textbf{GJ-II}+RAD run uses 100 mesh-free collocation points, $M=16$ Gauss--Jacobi nodes, and RAD: every 1000 steps, $10\%$ of the collocation points are resampled from 8192 residual-based candidates.

\begin{figure}[htbp]
    \centering
    \begin{subfigure}{.33\textwidth}
        \centering
        \includegraphics[width=1.0\textwidth]{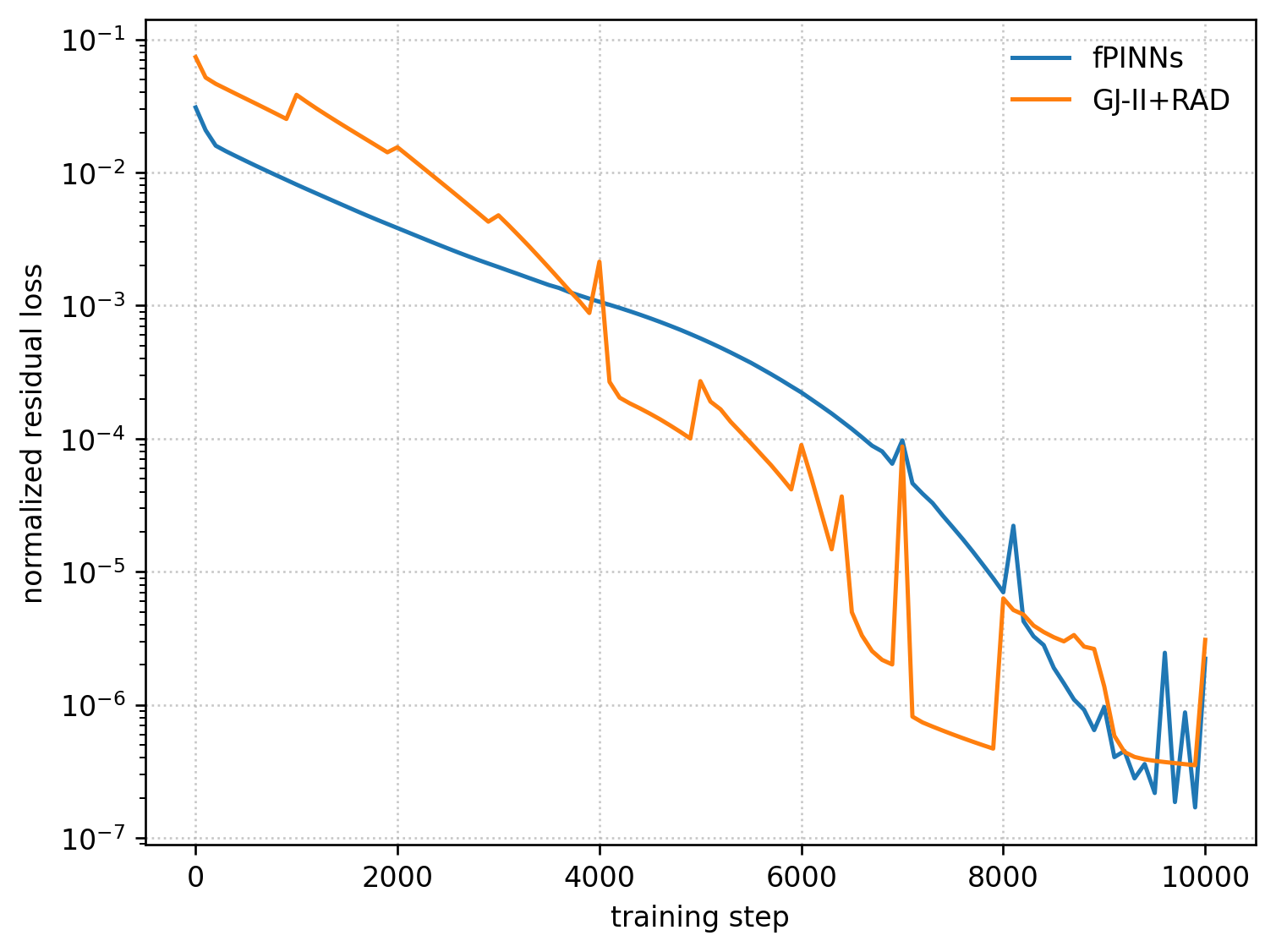}
        \caption{Loss.}
    \end{subfigure}%
    \begin{subfigure}{.33\textwidth}
        \centering
        \includegraphics[width=1.0\textwidth]{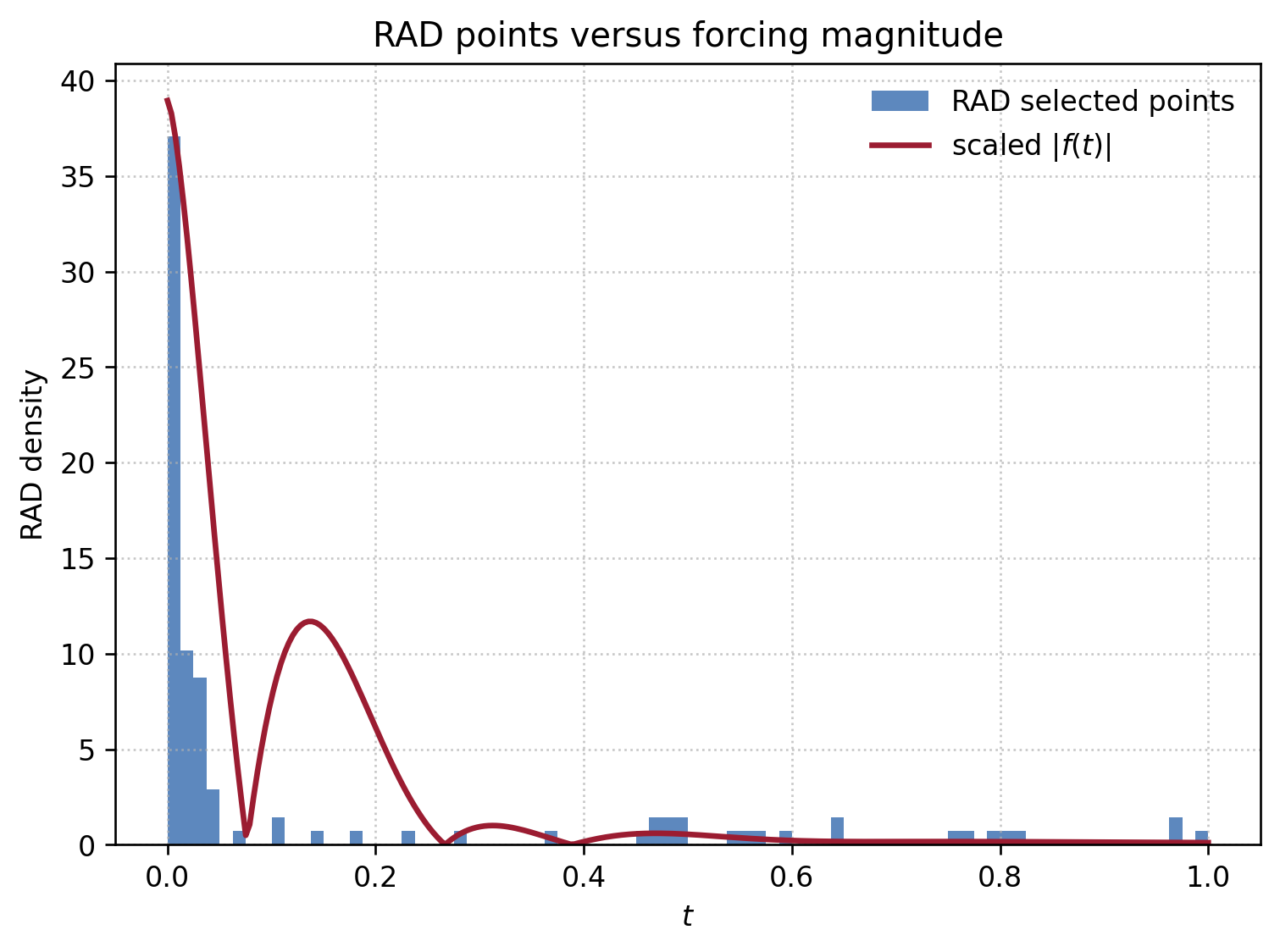}
        \caption{RAD points.}
    \end{subfigure}%
    \begin{subfigure}{.33\textwidth}
        \centering
        \includegraphics[width=1.0\textwidth]{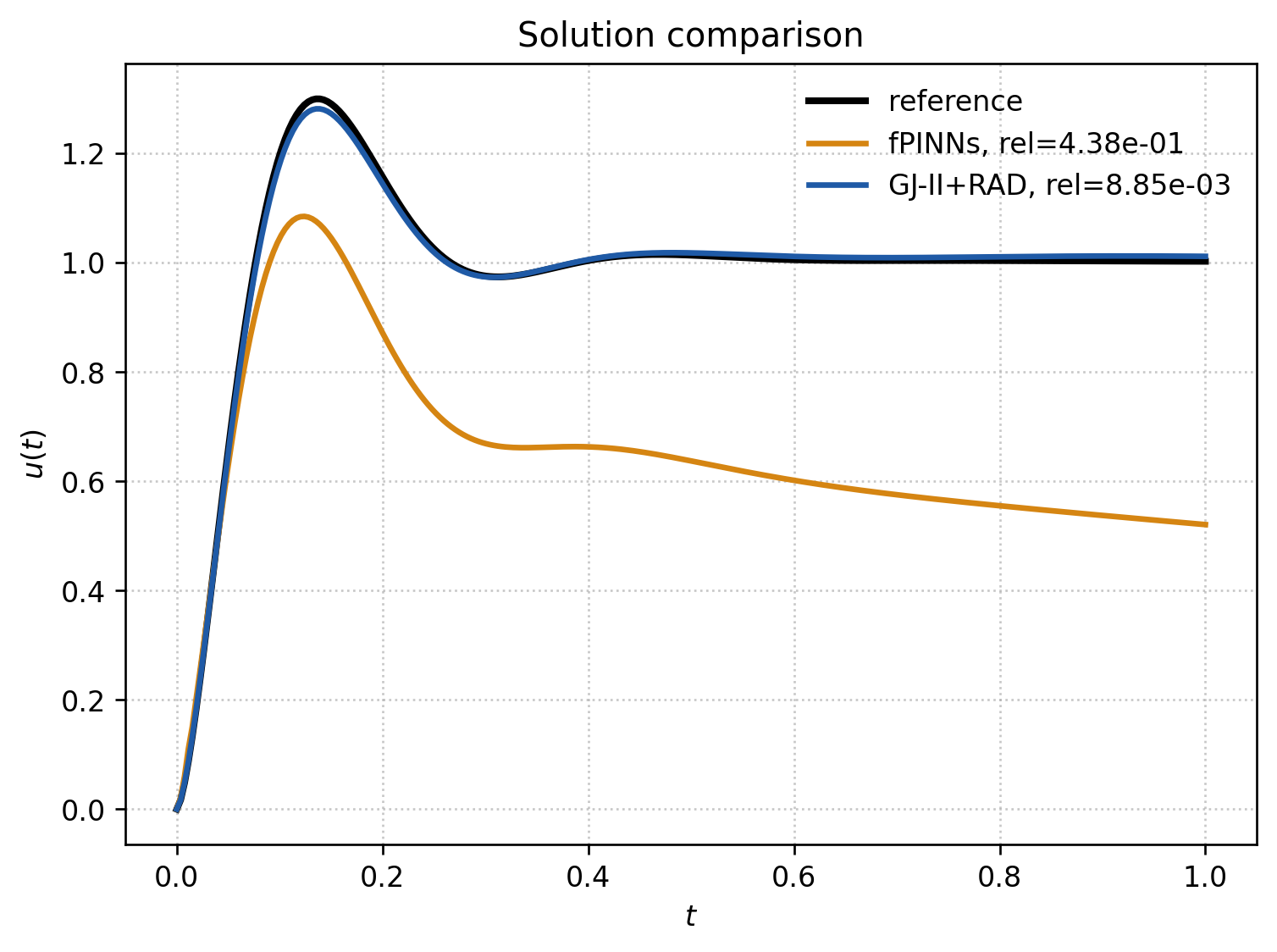}
        \caption{Solution.}
    \end{subfigure}%
    \caption{Fractional ODE experiment \eqref{eqn:ode-rad}. Left: normalized residual loss histories. Middle: RAD-selected points and the scaled forcing magnitude. Right: reference and learned solutions.}
    \label{fig:ode-rad-results}
\end{figure}

Figure~\ref{fig:ode-rad-results} summarizes the training behavior, adaptive sampling pattern, and final solution comparison. Both methods reduce the normalized residual loss, but this residual comparison alone does not reflect the final solution quality. Under the same 100-point setting, fPINNs requires 343.6 seconds, whereas \textbf{GJ-II}+RAD requires 117.0 seconds, so the mesh-free quadrature run is more than twice as fast. The RAD-selected points concentrate in the region where the forcing magnitude is large, which illustrates a practical advantage of the mesh-free \textbf{GJ-II} formulation: the residual can be evaluated directly at adaptively chosen points. In contrast, fPINNs is tied to a fixed temporal stencil, making this type of adaptive point redistribution less natural. The fPINNs baseline gives a relative error of $4.38\times10^{-1}$, while \textbf{GJ-II}+RAD reduces the relative error to $8.85\times10^{-3}$. This ODE test therefore highlights two advantages of the mesh-free \textbf{GJ-II} formulation over fPINNs in the same point-budget setting: lower computational time and direct compatibility with adaptive sampling.

\FloatBarrier

\captionsetup{font={color=blue}}
\subsection{Time-fractional diffusion-wave equation on a two-dimensional L-shaped domain}
\label{subsec:lshape}
We next consider a two-dimensional time-fractional diffusion-wave equation on the L-shaped domain
\[
    \Omega_L = (-1,1)^2\setminus ((0,1)\times(0,1)).
\]
The equation is
\begin{equation}
    \label{eqn:lshape}
    \begin{aligned}
        \dta u(t,x,y) - 0.25\Delta u(t,x,y) &= 0,~~ (t,x,y)\in (0,1]\times\Omega_L,\\
        u(t,x,y)&=0,~~(x,y)\in\partial\Omega_L,~~t\in(0,1],\\
        u(0,x,y)&=g(x,y),~~(x,y)\in\Omega_L,\\
        \partial_t u(0,x,y)&=0.2g(x,y),~~(x,y)\in\Omega_L,
    \end{aligned}
\end{equation}
where
\[
\begin{aligned}
g(x,y)=&
\exp\left(-\frac{(x+0.55)^2+(y+0.45)^2}{0.08}\right)
-0.85\exp\left(-\frac{(x+0.55)^2+(y-0.45)^2}{0.06}\right)\\
&+0.60\exp\left(-\frac{(x-0.45)^2+(y+0.55)^2}{0.06}\right).
\end{aligned}
\]
For this two-dimensional irregular-domain example, the relative errors are computed against numerical reference solutions generated on the same L-shaped domain.

For this experiment, we train to 100000 steps and record timing over windows of 5000 steps. The network has 6 hidden layers and 64 neurons per hidden layer. The domain, boundary, and initial batches contain 20000, 5000, and 5000 points, respectively. We use $M=256$ quadrature points for the Monte Carlo schemes and $M=64$ quadrature points for the Gauss-Jacobi schemes. The tested fractional orders are $\alpha=1.25,1.5,1.75$, and the four schemes \textbf{MC-I}, \textbf{MC-II}, \textbf{GJ-I}, and \textbf{GJ-II} are compared under the same training configuration. We also include two conventional fPINNs baselines. To match the domain size of the Monte Carlo and Gauss-Jacobi experiments, the first fPINNs setting uses $N_T=100$ temporal mesh points and $N_x=200$ randomly sampled spatial points, while the second uses $N_T=200$ and $N_x=100$. Both choices give $N_TN_x=20000$ interior points, and the boundary and initial batches are kept unchanged.

\begin{table}[htbp]
    \centering

    \begin{tabular}{|c|c|c|c|c|c|c|}
        \hline
        $\alpha$/Type & \textbf{MC-I} & \textbf{MC-II} & \textbf{GJ-I} & \textbf{GJ-II}
        & \begin{tabular}[c]{@{}c@{}}\textbf{fPINNs}\\($N_T=100$)\end{tabular}
        & \begin{tabular}[c]{@{}c@{}}\textbf{fPINNs}\\($N_T=200$)\end{tabular} \\ \hline
        $1.25$ & \begin{tabular}[c]{@{}c@{}}1.13e-1\\ 37m10s\end{tabular}
        & \begin{tabular}[c]{@{}c@{}}9.47e-2\\ 16m05s\end{tabular}
        & \begin{tabular}[c]{@{}c@{}}7.59e-2\\ 11m07s\end{tabular}
        & \begin{tabular}[c]{@{}c@{}}7.32e-2\\ 5m34s\end{tabular}
        & \begin{tabular}[c]{@{}c@{}}3.37e-1\\ 5m03s\end{tabular}
        & \begin{tabular}[c]{@{}c@{}}3.37e-1\\ 7m07s\end{tabular} \\ \hline
        $1.50$ & \begin{tabular}[c]{@{}c@{}}6.40e-2\\ 37m10s\end{tabular}
        & \begin{tabular}[c]{@{}c@{}}7.81e-2\\ 16m05s\end{tabular}
        & \begin{tabular}[c]{@{}c@{}}8.30e-2\\ 11m07s\end{tabular}
        & \begin{tabular}[c]{@{}c@{}}8.05e-2\\ 5m34s\end{tabular}
        & \begin{tabular}[c]{@{}c@{}}1.43e-1\\ 5m12s\end{tabular}
        & \begin{tabular}[c]{@{}c@{}}4.27e-1\\ 6m48s\end{tabular} \\ \hline
        $1.75$ & \begin{tabular}[c]{@{}c@{}}8.70e-2\\ 37m12s\end{tabular}
        & \begin{tabular}[c]{@{}c@{}}1.67e-1\\ 16m06s\end{tabular}
        & \begin{tabular}[c]{@{}c@{}}3.70e-2\\ 11m06s\end{tabular}
        & \begin{tabular}[c]{@{}c@{}}2.04e-1\\ 5m34s\end{tabular}
        & \begin{tabular}[c]{@{}c@{}}6.88e-2\\ 4m57s\end{tabular}
        & \begin{tabular}[c]{@{}c@{}}4.58e-1\\ 6m45s\end{tabular} \\ \hline
    \end{tabular}
    \caption{$L^2$ relative error and computational time for equation \eqref{eqn:lshape}. In each entry, the first row is the relative error against the numerical reference solution, and the second row is the recorded computational time for one 5000-step timing window. For fPINNs, $N_T=100$ uses $N_x=200$, and $N_T=200$ uses $N_x=100$.}
    \label{tab:lshape}
\end{table}

\begin{figure}[htbp]
    \centering
    \includegraphics[scale=0.35]{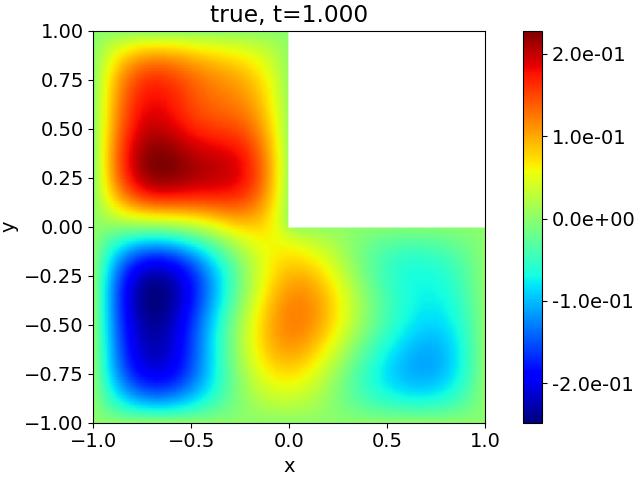}
    \caption{Reference solution for equation \eqref{eqn:lshape} with $\alpha=1.75$ at $t=1.0$.}
    \label{fig:lshape-alpha175-true}
\end{figure}

\begin{figure}[htbp]
    \centering
    \begin{subfigure}{.25\textwidth}
        \centering
        \includegraphics[height=0.75\textwidth,width=1.0\textwidth]{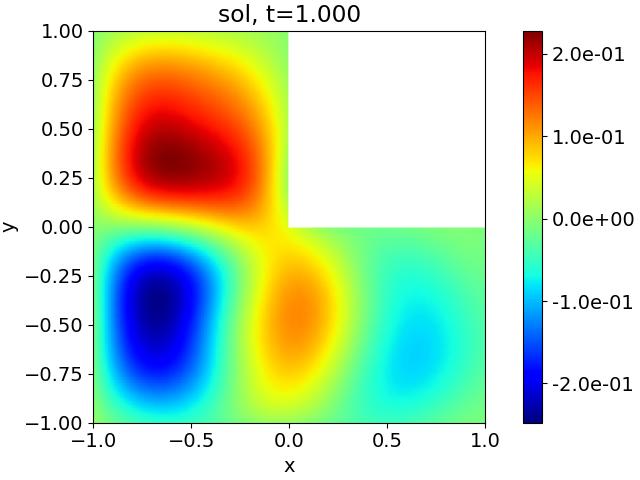}
    \end{subfigure}%
    \begin{subfigure}{.25\textwidth}
        \centering
        \includegraphics[height=0.75\textwidth,width=1.0\textwidth]{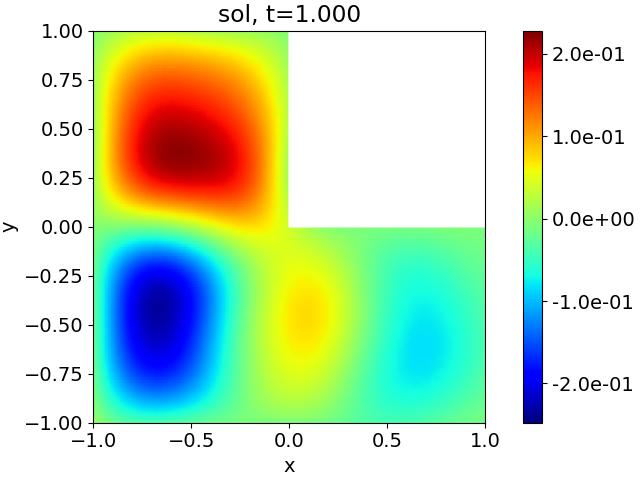}
    \end{subfigure}%
    \begin{subfigure}{.25\textwidth}
        \centering
        \includegraphics[height=0.75\textwidth,width=1.0\textwidth]{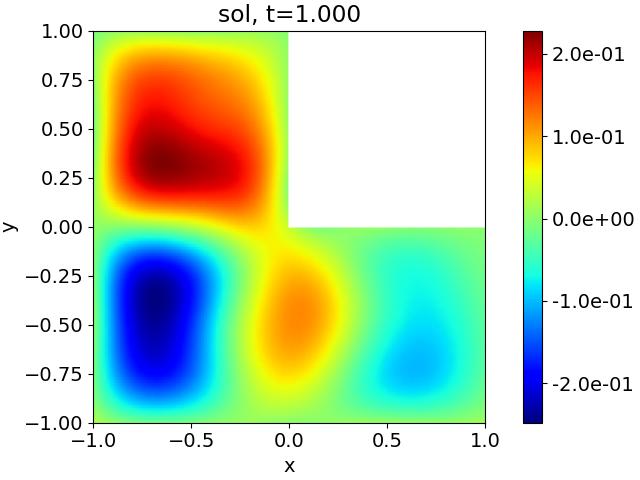}
    \end{subfigure}%
    \begin{subfigure}{.25\textwidth}
        \centering
        \includegraphics[height=0.75\textwidth,width=1.0\textwidth]{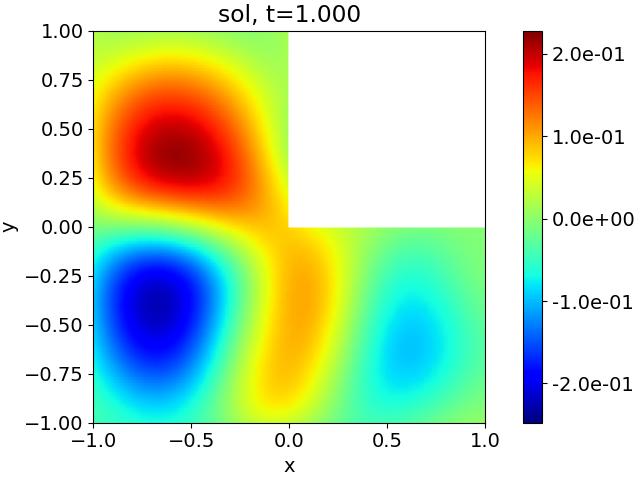}
    \end{subfigure}%
    \newline
    \raggedleft
    \begin{subfigure}{.25\textwidth}
        \centering
        \includegraphics[height=0.75\textwidth,width=1.0\textwidth]{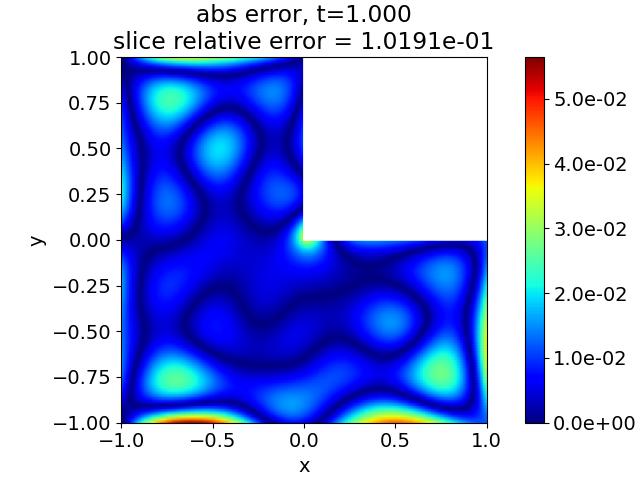}
        \caption{\textbf{MC-I}.}
    \end{subfigure}%
    \begin{subfigure}{.25\textwidth}
        \centering
        \includegraphics[height=0.75\textwidth,width=1.0\textwidth]{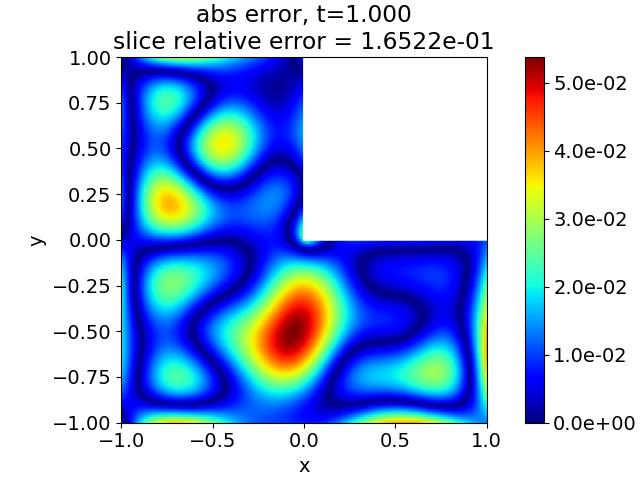}
        \caption{\textbf{MC-II}.}
    \end{subfigure}%
    \begin{subfigure}{.25\textwidth}
        \centering
        \includegraphics[height=0.75\textwidth,width=1.0\textwidth]{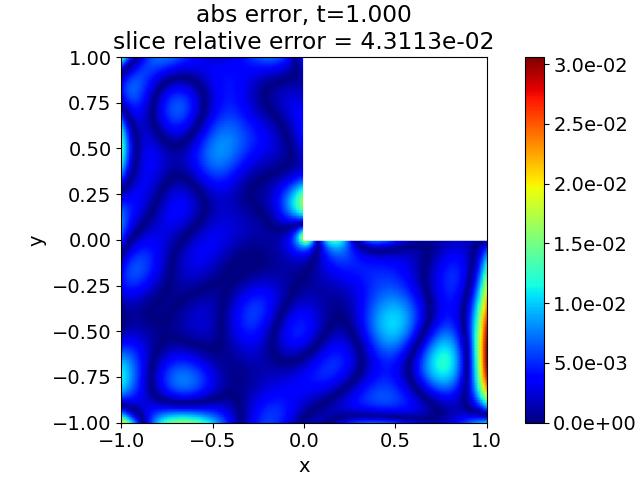}
        \caption{\textbf{GJ-I}.}
    \end{subfigure}%
    \begin{subfigure}{.25\textwidth}
        \centering
        \includegraphics[height=0.75\textwidth,width=1.0\textwidth]{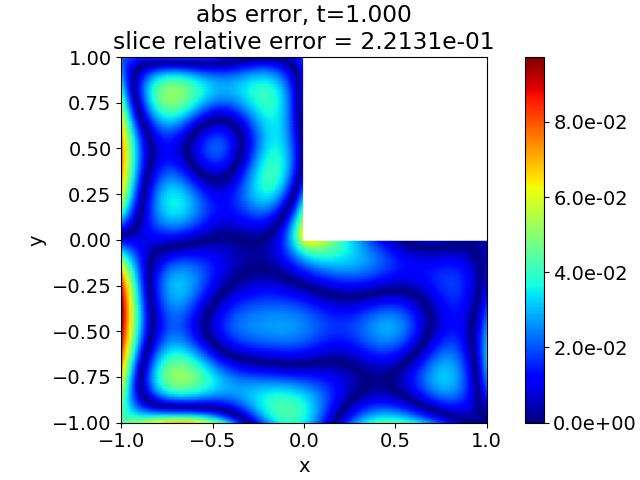}
        \caption{\textbf{GJ-II}.}
    \end{subfigure}%
    \caption{Numerical solutions and time-slice error profiles for equation \eqref{eqn:lshape} with $\alpha=1.75$ at $t=1.0$. First row: learned numerical solutions. Second row: absolute error time-slice plots from the corresponding experiments. In error plots, slice relative error means the relative error at $t=1.0$.}
    \label{fig:lshape-alpha175}
\end{figure}

\begin{figure}[htbp]
    \centering
    \begin{subfigure}{.25\textwidth}
        \centering
        \includegraphics[height=0.75\textwidth,width=1.0\textwidth]{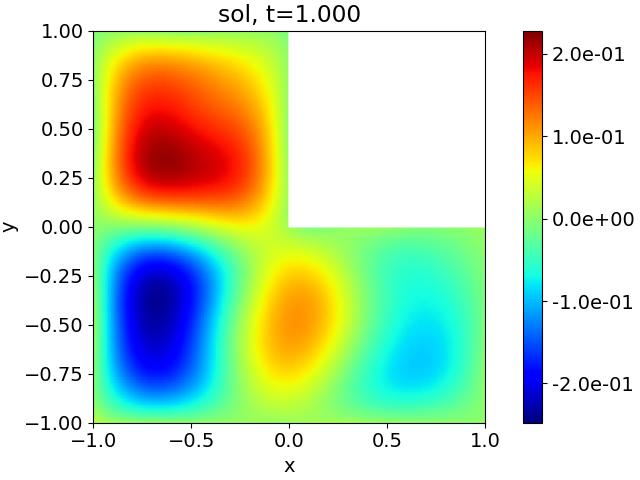}
        \caption{$N_T=100$, sol.}
    \end{subfigure}%
    \begin{subfigure}{.25\textwidth}
        \centering
        \includegraphics[height=0.75\textwidth,width=1.0\textwidth]{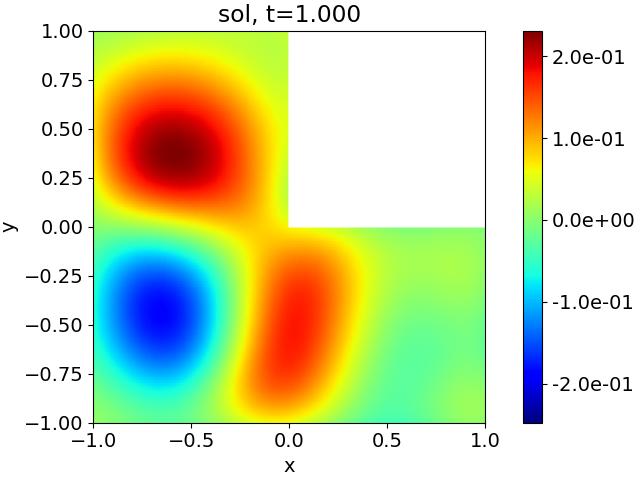}
        \caption{$N_T=200$, sol.}
    \end{subfigure}%
    \begin{subfigure}{.25\textwidth}
        \centering
        \includegraphics[height=0.75\textwidth,width=1.0\textwidth]{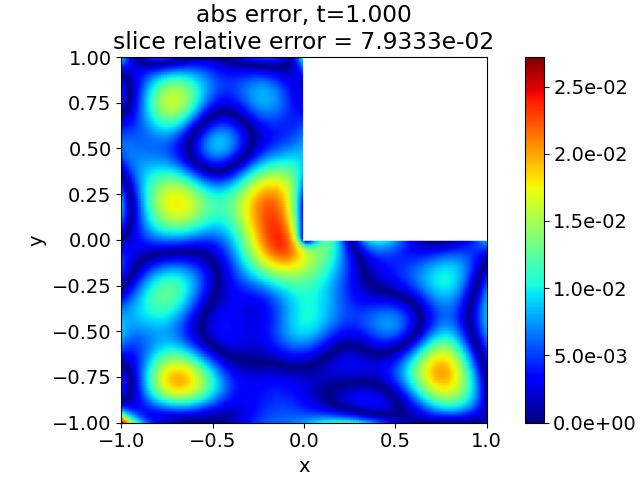}
        \caption{$N_T=100$, error.}
    \end{subfigure}%
    \begin{subfigure}{.25\textwidth}
        \centering
        \includegraphics[height=0.75\textwidth,width=1.0\textwidth]{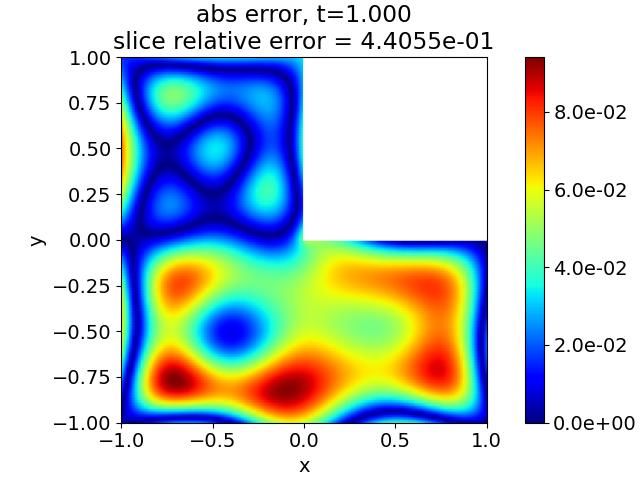}
        \caption{$N_T=200$, error.}
    \end{subfigure}%
    \caption{Numerical solutions and time-slice absolute error profiles of standard fPINNs for equation \eqref{eqn:lshape} with $\alpha=1.75$ at $t=1.0$. Here $N_T=100$ uses $N_x=200$, and $N_T=200$ uses $N_x=100$.}
    \label{fig:lshape-alpha175-fpinns}
\end{figure}

For visualization, Figure~\ref{fig:lshape-alpha175-true} shows the numerical reference solution at $\alpha=1.75$ and $t=1.0$. The corresponding learned solutions and error profiles for the transformed schemes are shown in Figure~\ref{fig:lshape-alpha175}. The fPINNs time-slice error profiles are reported separately in Figure~\ref{fig:lshape-alpha175-fpinns}.

\FloatBarrier

Table \ref{tab:lshape} shows that the Gauss-Jacobi implementation is substantially faster than the Monte Carlo implementation under the tested two-dimensional setting. For the Type-I schemes, the average time over 5000 steps is reduced from about 37 minutes for \textbf{MC-I} to about 11 minutes for \textbf{GJ-I}. For the Type-II schemes, the corresponding time is reduced from about 16 minutes for \textbf{MC-II} to about 5.5 minutes for \textbf{GJ-II}. This confirms the practical advantage of Gauss-Jacobi quadrature in this irregular-domain experiment.

The same table also highlights the computational gain of Type-II over Type-I. For Monte Carlo quadrature, \textbf{MC-II} reduces the 5000-step time from about 37 minutes to about 16 minutes. For Gauss-Jacobi quadrature, \textbf{GJ-II} reduces the time from about 11 minutes to about 5.5 minutes. These savings are consistent with the transformed formulation, which avoids differentiating the neural network at shifted quadrature points. The resulting error levels remain problem-dependent, especially near $\alpha=1.75$, but the timing comparison clearly supports the computational efficiency of the Type-II formulation.

The fPINNs baselines further clarify the cost of using a mesh-based fractional derivative approximation. With $N_T=100$, fPINNs can be inexpensive because the temporal mesh is relatively coarse. However, when the temporal resolution is increased to $N_T=200$ (equivalently, when the temporal mesh size $\tau$ is reduced from $0.01$ to $0.005$), the 5000-step time rises to about 6.8--7.1 minutes, already exceeding the about 5.5 minutes required by \textbf{GJ-II}. Thus, if high-resolution fPINNs calculations are pursued by taking larger $N_T$ or smaller $\tau$, their computational cost can overtake the proposed \textbf{GJ-II} scheme. This supports the effectiveness of \textbf{GJ-II}, which keeps the sampling formulation mesh-free while avoiding the finite-difference temporal stencil used by fPINNs.

This comparison also explains why the timing gap is much clearer in the ODE diagnostic of Section~\ref{subsec:ode-rad} than in the present two-dimensional PDE experiment. In the ODE test, fPINNs uses $N_T=100$ temporal grid points, while \textbf{GJ-II}+RAD is also evaluated on 100 collocation points; the two methods therefore have a directly comparable point budget, and the cost of the fPINNs temporal stencil becomes visible. In the L-shaped PDE test, the fPINNs baseline with $N_T=100$ still uses only 100 temporal stencil levels, whereas \textbf{GJ-II} evaluates the transformed residual at 20000 mesh-free interior collocation points, each requiring the time derivative at the base point. This larger collocation workload can mask the timing advantage seen in the ODE test. Nevertheless, increasing $N_T$ from 100 to 200 already makes fPINNs slower than \textbf{GJ-II}, showing that the apparent speed of the coarse fPINNs baseline is not robust to temporal refinement.

This mesh-based structure also makes standard fPINNs less compatible with adaptive sampling methods such as RAD. RAD repeatedly changes the domain collocation distribution according to the residual, whereas fPINNs relies on a prescribed temporal grid and associated fractional-difference stencil. Applying RAD would require repeatedly reconciling adaptive samples with this fixed stencil structure, weakening the main practical advantage of adaptive point selection. In contrast, the Monte Carlo and Gauss-Jacobi estimators used in tDWfPINNs can be evaluated directly at arbitrary collocation points.

Taken together with the ODE experiment, these results suggest that fPINNs is not a particularly advantageous strategy within a PINN framework: its fixed temporal stencil weakens the mesh-free character of PINNs, and its computational time does not provide a decisive advantage over \textbf{GJ-II}. The mesh-free quadrature formulation is therefore preferable when adaptive sampling or temporal refinement is important.

\FloatBarrier

\section{Conclusion and future works}
\label{sec:conclusion}
This work develops a transformed Caputo-residual construction for fractional PINNs applied to time-fractional diffusion-wave equations with order $\alpha\in(1,2)$. The central idea is to transfer the shifted time-derivative evaluation in the quadrature integrand to an analytic kernel, so that the residual can be assembled from shifted neural-network values. This directly targets the main automatic-differentiation cost in direct mesh-free quadrature implementations of the Caputo derivative.

At the operator and quadrature levels, we establish the equivalence between the transformed representations and the original Caputo derivative, the removable nature of the endpoint singularities, and a representation-level stability bound. We also derive Monte Carlo convergence estimates, Gauss--Jacobi quadrature error estimates, denominator-regularization bias estimates, and computational-complexity estimates. The storage and profiler-FLOP analyses show that the computational advantage of the transformed formulation is relative to the direct Type-I residual construction and becomes more pronounced as the number of collocation points $N$, quadrature points $M$, input dimension $d$, network depth $L$, or hidden width $H$ increases.

The numerical experiments connect this analysis to PINN training. The one-dimensional diffusion-wave examples show that the transformed schemes preserve comparable accuracy while reducing training cost. The Gauss--Jacobi implementation is particularly effective in the tested settings because it uses fewer quadrature nodes for smooth transformed kernels, while the transformed residual removes shifted time-derivative evaluations at those nodes. The fractional ODE diagnostic and the two-dimensional L-shaped-domain experiment further show that the mesh-free \textbf{GJ-II} implementation can provide a favorable accuracy--cost balance compared with standard fPINNs when adaptive sampling or temporal refinement is important. The added L-shaped-domain experiment also verifies the implementation beyond one-dimensional benchmarks and rectangular geometries.


Several directions remain important for future work. A complete convergence and error analysis for the full PINN training process, including neural-network approximation, nonconvex optimization, sampling, and loss-balancing effects, would complement the operator-level analysis developed here. The present transformation is tailored to the Caputo derivative structure in time-fractional diffusion-wave equations with $\alpha\in(1,2)$; extending similar transformed residual constructions to broader integro-differential equations is a natural direction for further research. Application-oriented extensions involving three-dimensional domains, heterogeneous coefficients, experimentally calibrated material parameters, realistic irregular geometries, and inverse diffusion-wave problems would further develop the path from the present benchmark validation toward complex physical simulations.

\section*{CRediT authorship contribution statement}
\textbf{Zhengqi Zhang}: Writing-review \& editing, Writing-original draft, Visualization, Validation, Software, Methodology, Investigation, Formal analysis, Data curation. \textbf{Jing Li}: Writing-review \& editing, Supervision, Project administration, Methodology, Funding acquisition, Conceptualization.

\section*{Declaration of competing interest}
  The authors declare that they have no known competing financial interests or personal relationships that could have appeared to influence the work reported in this paper.

\section*{Data availability}
The code in this work is available from the GitHub repository \url{https://github.com/Zenki229/tDWfPINN}.

\appendix

\section{Derivative-level validation and implementation diagnostics}
\label{app:num-valid}

This appendix provides derivative-level numerical evidence for the quadrature and implementation consequences of the analysis in Section~\ref{subsec:methodology}. The identities and stability estimates in Theorem~\ref{thm:represent_2} and Corollary~\ref{cor:stability} are established analytically, and the tests below complement them by isolating four effects that enter practical implementations of the estimators in Corollary~\ref{coro:4}: Monte Carlo sampling error, Gauss--Jacobi quadrature error, denominator regularization error, and finite-precision cancellation in endpoint difference quotients. These derivative-evaluation tests separate quadrature effects from neural-network approximation and optimization effects.

\subsection{Reference functions and stable quotient evaluation}

Unless otherwise stated, the smooth benchmark is
\begin{equation}
    f(t)=e^{-t},\qquad \text{at}~t_0=1.5,
    \label{eq:valid-exp-benchmark}
\end{equation}
for which
\begin{equation*}
    \partial_t^\alpha f(t)=t^{2-\alpha}E_{1,3-\alpha}(-t).
\end{equation*}
Monte Carlo samples are drawn from $\mathrm{Beta}(2-\alpha,1)$, and Gauss--Jacobi nodes and weights are computed for the weight $\tau^{1-\alpha}$ on $[0,1]$. If $\{\tau_i\}_{i=1}^M$ denotes the resulting Gauss--Jacobi nodes, equivalently obtained from nodes $\{x_i\}_{i=1}^M$ on $[-1,1]$ by $\tau_i=(1+x_i)/2$, we write
\begin{equation*}
    \tau_{\min}:=\min_{1\le i\le M}\tau_i .
\end{equation*}
Thus $\tau_{\min}$ is the Gauss--Jacobi node closest to the endpoint $\tau=0$, and it depends on both $\alpha$ and $M$. Unless explicitly stated otherwise, errors are relative derivative errors,
\begin{equation*}
    \mathrm{err}=\abs{\dta f(t_0)-\btau f(t_0)}/\abs{\dta f(t_0)}.
\end{equation*}

We distinguish carefully between \emph{raw} and \emph{stable} quotient evaluations. For the exponential benchmark, with $h=t\tau$, the raw kernels are
\begin{equation*}
    K_f(t,\tau)=\frac{f'(t)-f'(t-h)}{h},\qquad
    H_f(t,\tau)=\frac{f(t)-f(t-h)-h f'(t)}{h^2}.
\end{equation*}
Their stable algebraic forms are evaluated as
\begin{equation}
    K_f^{\rm stab}(t,\tau)=e^{-t}\frac{\operatorname{expm1}(h)}{h},
    \qquad
    H_f^{\rm stab}(t,\tau)=e^{-t}\frac{h-\operatorname{expm1}(h)}{h^2},
    \label{eq:valid-stable-exp}
\end{equation}
Here $\operatorname{expm1}(h)$ from Numpy denotes the standard numerically stable
evaluation of $e^h-1$. It is used because, for small $h$, directly forming
$e^h-1$ subtracts two nearly equal numbers and can lose most of the significant
digits before division by $h$ or $h^2$. In the Type-II expression, the remaining
quantity $h-\operatorname{expm1}(h)=O(h^2)$ is still an endpoint-scale
difference, so the stable formula is combined with the Taylor replacement below
when $h$ is below machine-resolution thresholds. We use the Taylor endpoint replacement
\begin{equation*}
    K_f(t,0)=e^{-t},
    \qquad
    H_f(t,0)=-\frac12e^{-t}.
\end{equation*}
This definition is important: Lemma~\ref{lem:type-kernels} removes the endpoint singularity at the continuous level, but a raw floating-point quotient need not evaluate that removable singularity accurately.

\subsection{Sensitivity with respect to the fractional order}
\label{sec:sense_to_frac_order}

\begin{figure}[htbp]
    \centering
    \includegraphics[width=0.92\textwidth]{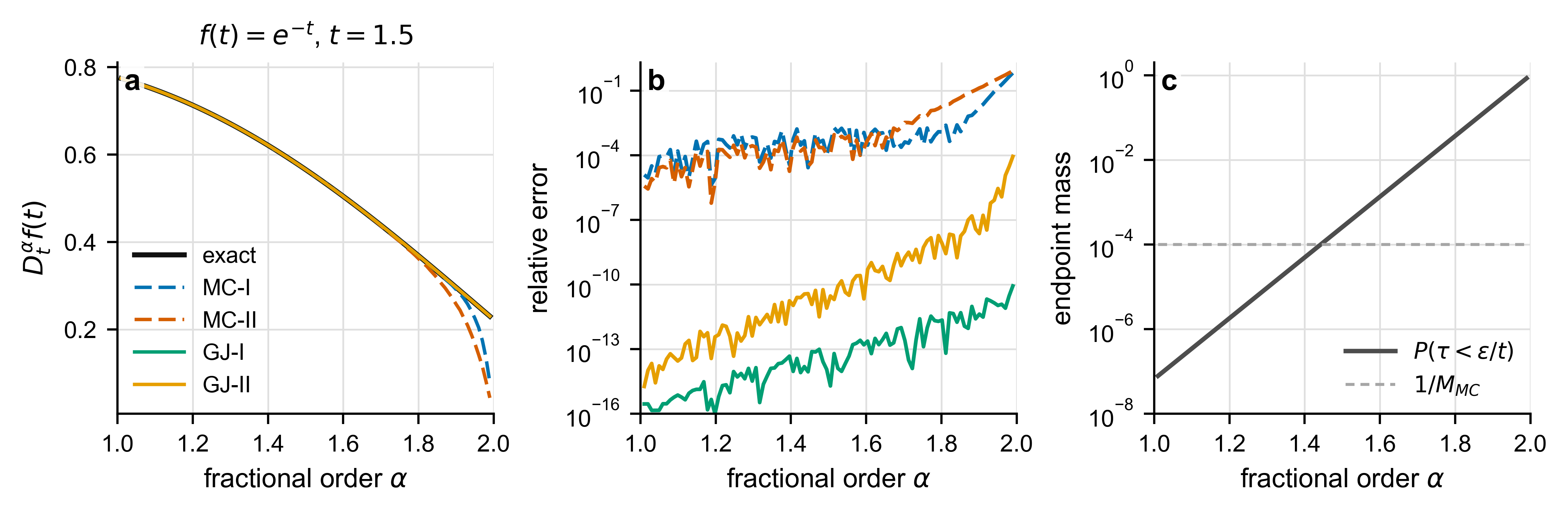}
    \caption{Derivative-level sensitivity with respect to $\alpha$ for the benchmark \eqref{eq:valid-exp-benchmark}. The plot uses $M_{\rm MC}=10^4$, $M_{\rm GJ}=100$, $\varepsilon_I=10^{-16}$, and $\varepsilon_{II}=10^{-7}$. The middle panel uses the same legend as the left panel. The right panel reports the endpoint mass \eqref{eq:valid-endpoint-mass}. The deterioration of the Monte Carlo curves near $\alpha=2$ is consistent with increased endpoint exposure of the practical estimator, not with a change of the theoretical Monte Carlo rate.}
    \label{fig:valid-alpha}
    \label{fig:diff_al}
\end{figure}

For Monte Carlo sampling, the endpoint mass is explicit:
\begin{equation}
    \mathbb P(\tau<\varepsilon/t)=(\varepsilon/t)^{2-\alpha},
    \qquad \tau\sim \mathrm{Beta}(2-\alpha,1).
    \label{eq:valid-endpoint-mass}
\end{equation}
Figure~\ref{fig:valid-alpha} was generated using the benchmark \eqref{eq:valid-exp-benchmark}, with  $100$ equally spaced
orders $\alpha\in[1.01,1.99]$, $M_{\rm MC}=10^4$ Monte Carlo samples, and
$M_{\rm GJ}=100$ Gauss--Jacobi nodes. The denominator
regularization is applied only to the Monte Carlo quotients in this test, with
absolute memory cutoffs $\varepsilon_I=10^{-16}$ for MC-I and
$\varepsilon_{II}=10^{-7}$ for MC-II. The right panel plots
$(\varepsilon_{II}/t)^{2-\alpha}$ together with the reference level
$1/M_{\rm MC}$.

\begin{figure}[htbp]
    \centering
    \includegraphics[scale=0.6]{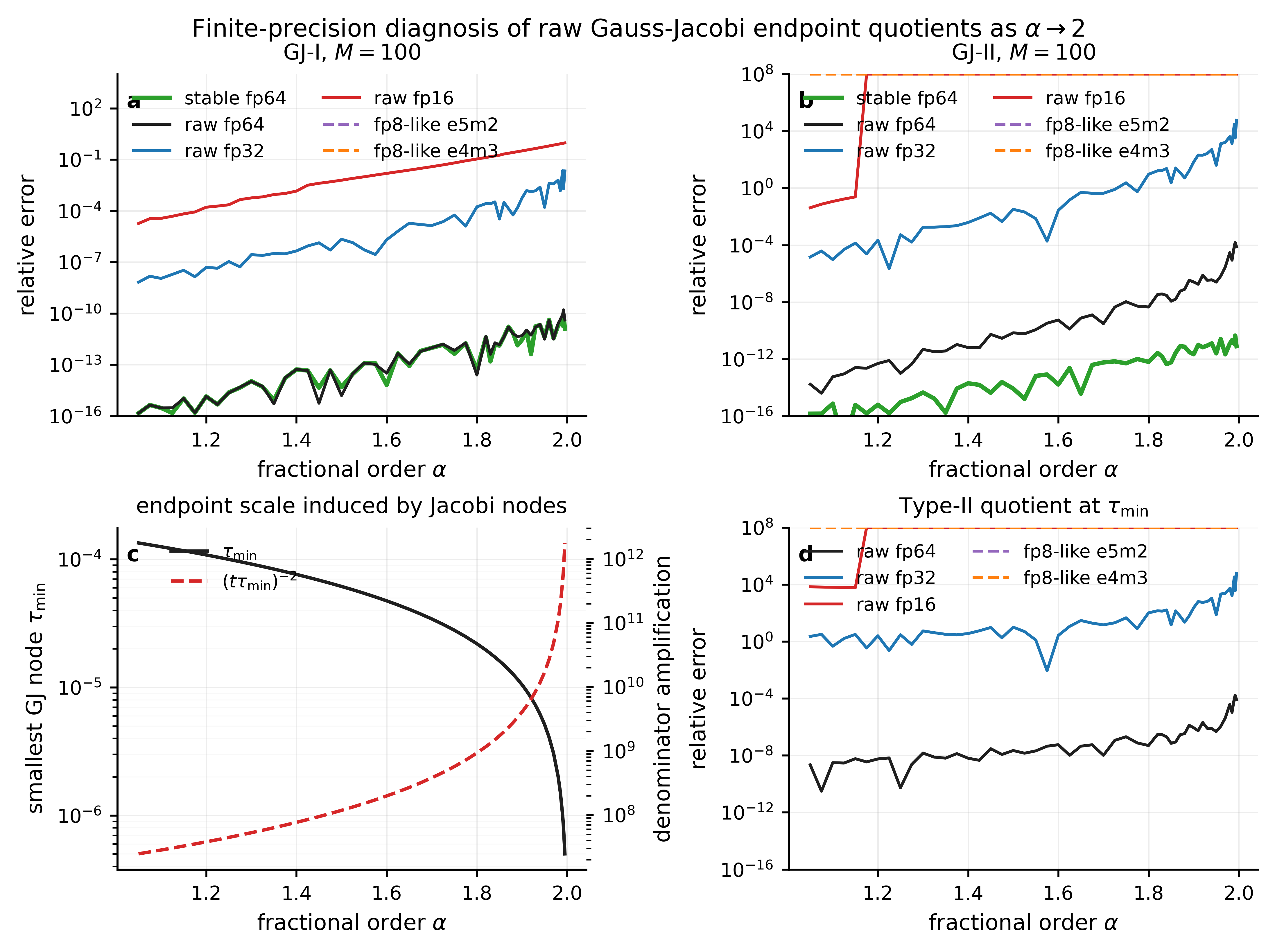}
    \caption{Finite-precision diagnosis of raw Gauss--Jacobi endpoint quotients as $\alpha\to2$ for $f(t)=e^{-t}$, $t=1.5$, and $M=100$. The stabilized fp64 implementation uses \texttt{expm1} and Taylor endpoint replacement. The raw GJ-II error increases because $\tau_{\min}\to0$ and $(t\tau_{\min})^{-2}$ amplifies endpoint cancellation.}
    \label{fig:valid-alpha2-precision}
\end{figure}

Thus, as $\alpha\to2$, a larger fraction of Monte Carlo samples falls in the
interval where raw quotient evaluation and denominator regularization are most
exposed to endpoint perturbations. The Monte Carlo curves deteriorate near $\alpha=2$, whereas the
Gauss--Jacobi curves remain much more accurate for the smooth benchmark. This
Monte Carlo behavior is controlled in practice by the cutoff $\delta$, whose
bias--conditioning effect is examined in
Section~\ref{subsubsec:valid-cutoff-tradeoff}. The small growth of the raw
GJ-II curve near $\alpha=2$ is instead the finite-precision endpoint-quotient
effect diagnosed in Figure~\ref{fig:valid-alpha2-precision}, not a loss of the
Gauss--Jacobi quadrature rate. The Monte Carlo deterioration should not be
interpreted as a failure of the $M^{-1/2}$ Monte Carlo rate in
Proposition~\ref{prop:mc-convergence}; rather, it shows that the practical
regularized estimator becomes increasingly endpoint-sensitive near the
second-order limit.

The same fractional-order sweep also serves as a compact diagnostic for the
Gauss--Jacobi curves. Since the Monte Carlo deterioration near $\alpha=2$ is
controlled in practice by the cutoff $\delta$, as examined later in
Section~\ref{subsubsec:valid-cutoff-tradeoff}, we separately checked whether the
small raw GJ-II growth is caused by quadrature or by floating-point quotient
evaluation. Using the benchmark \eqref{eq:valid-exp-benchmark}, $M=100$
Gauss--Jacobi nodes, and raw quotient evaluations in fp64, fp32, fp16, and
fp8-like quantized arithmetic, Figure~\ref{fig:valid-alpha2-precision} shows
that the growth is a finite-precision endpoint-quotient effect. As
$\alpha\to2$, the Jacobi exponent $1-\alpha$ approaches $-1$, so
$\tau_{\min}\to0$ and the Type-II factor $(t\tau_{\min})^{-2}$ amplifies
cancellation in the $O(h^2)$ numerator. Type-I is less exposed because its raw
quotient has only first-order endpoint amplification.

\subsection{Sensitivity with respect to the number of quadrature nodes}
\label{subsubsec:valid-quadrature-nodes}
This experiment varies the number of quadrature nodes $M$ while keeping the
benchmark \eqref{eq:valid-exp-benchmark} fixed. Its purpose is to separate two
different effects. For Monte Carlo quadrature, increasing $M$ reduces the
sampling error. We therefore draw repeated independent samples for each fixed
$M$ and report the median relative errors in Table~\ref{tab:valid-M-effect}.
The MC-II median error decreases from $1.35\times10^{-2}$ at $M=10$ to
$2.24\times10^{-4}$ at $M=10240$, consistent with the $M^{-1/2}$ behavior in
Proposition~\ref{prop:mc-convergence}.

\begin{figure}[htbp]
    \centering
    \includegraphics[width=0.92\textwidth]{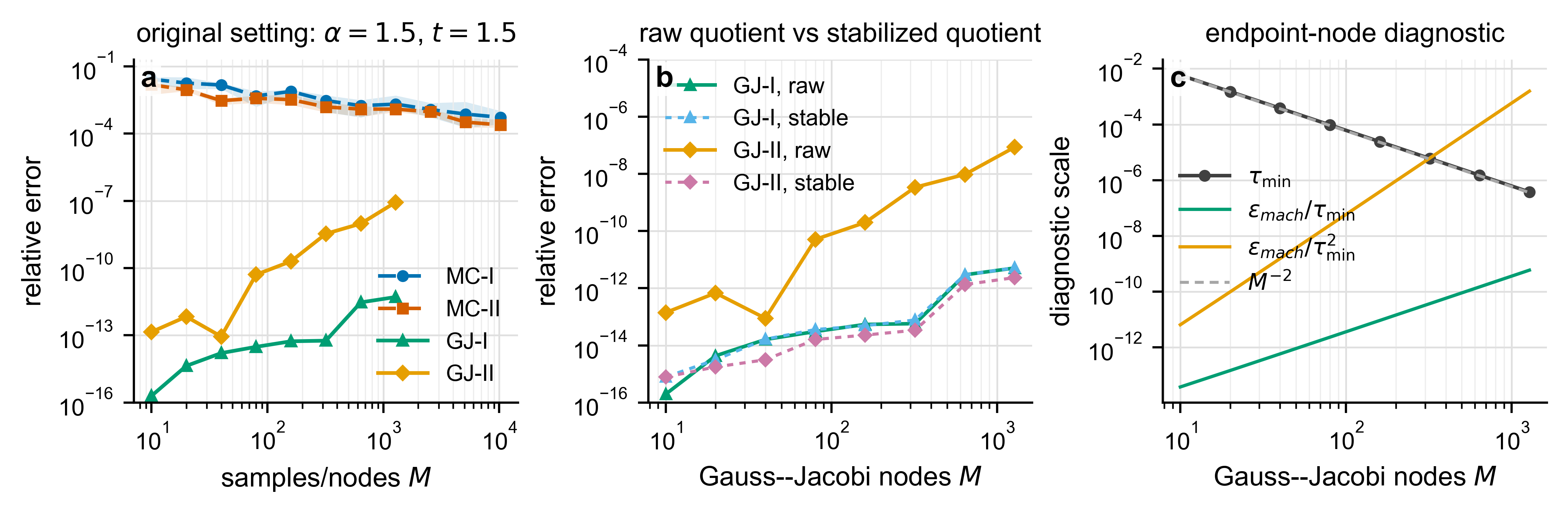}
    \caption{Effect of the number of quadrature points $M$. Left: repeated Monte Carlo runs show decreasing median errors with interquartile variability. Middle: Gauss--Jacobi raw and stable quotients separate when endpoint cancellation dominates. Right: endpoint scale of the smallest Gauss--Jacobi node, $\tau_{\min}\sim M^{-2}$, which explains why increasing $M$ can expose raw endpoint quotients to stronger roundoff amplification.}
    \label{fig:valid-M}
    \label{fig:diff_M_err}
\end{figure}

For Gauss--Jacobi quadrature, increasing $M$ has an additional finite-precision
effect. As shown in Section~\ref{sec:sense_to_frac_order}, when $\alpha\to2$,
endpoint nodes move toward $\tau=0$ and raw GJ-II quotients become sensitive to
cancellation; see Figure~\ref{fig:valid-alpha2-precision}. Figure~\ref{fig:valid-M}
shows the same endpoint mechanism as $M$ increases: for fixed $\alpha$, the
smallest Jacobi node satisfies approximately $\tau_{\min}=O(M^{-2})$. Thus
large $M$ makes $h=t\tau_{\min}$ very small, and the raw Type-II quotient can be
dominated by the amplification factor $(t\tau_{\min})^{-2}$.
Table~\ref{tab:valid-M-effect} and Figure~\ref{fig:valid-M} show this
raw--stable separation. The stable quotient continues to be accurate, so the
effect is a floating-point implementation issue rather than a deterioration of
the Gauss--Jacobi quadrature rule.

\begin{table}[htbp]
\centering
\caption[Representative quadrature-size effects]{Representative quadrature-size effects for $f(t)=e^{-t}$, $t=1.5$, and $\alpha=1.5$. The Monte Carlo columns report the median relative error over repeated independent samples at the same $M$. The Gauss--Jacobi columns compare raw and stable Type-II quotient evaluations.}
\label{tab:valid-M-effect}
\begin{tabular}{c|cc|ccc}
\hline
$M$ & \shortstack{MC-I median\\rel. error} & \shortstack{MC-II median\\rel. error} & $\tau_{\min}$ & \shortstack{GJ-II raw\\rel. error} & \shortstack{GJ-II stable\\rel. error} \\
\hline
$10$    & $2.67\times10^{-2}$ & $1.35\times10^{-2}$ & $5.86\times10^{-3}$ & $1.40\times10^{-13}$ & $7.84\times10^{-16}$ \\
$80$    & $7.19\times10^{-3}$ & $3.49\times10^{-3}$ & $9.58\times10^{-5}$ & $5.09\times10^{-11}$ & $1.67\times10^{-14}$ \\
$320$   & $4.07\times10^{-3}$ & $1.43\times10^{-3}$ & $6.01\times10^{-6}$ & $3.39\times10^{-9}$  & $3.37\times10^{-14}$ \\
$1280$  & $1.85\times10^{-3}$ & $1.01\times10^{-3}$ & $3.76\times10^{-7}$ & $8.66\times10^{-8}$  & $2.36\times10^{-12}$ \\
$10240$ & $8.85\times10^{-4}$ & $2.24\times10^{-4}$ & -- & -- & -- \\
\hline
\end{tabular}
\end{table}

\begin{figure}[htbp]
    \centering
    \includegraphics[width=0.92\textwidth]{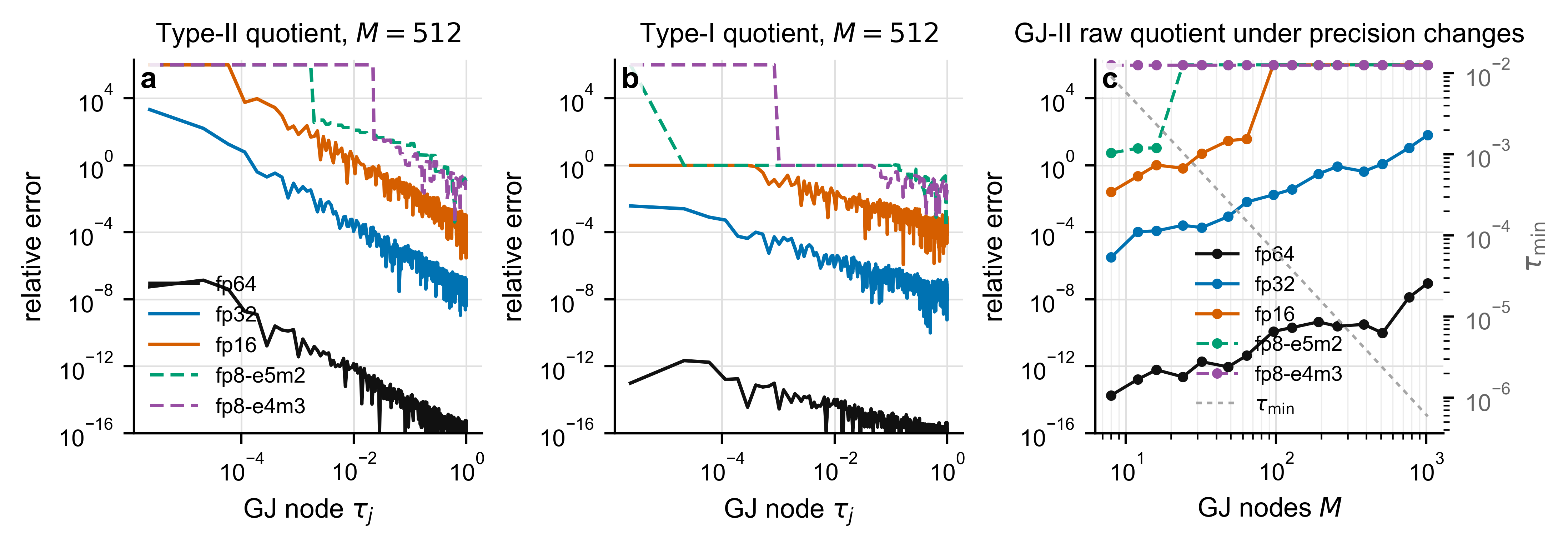}
    \caption{Finite-precision diagnosis of raw endpoint quotients. Left and middle: local relative errors of the raw Type-II and Type-I quotients at all Gauss--Jacobi nodes $\tau_j$ of a fixed $M=512$ rule. The middle panel uses the same legend as the left panel. Right: resulting GJ-II derivative error as $M$ changes. The reference is the stabilized fp64 quotient. Type-II is more vulnerable to precision loss because the removable singularity is evaluated through a second-order denominator. The fp8 curves are quantization diagnostics, not full mixed-precision training experiments.}
    \label{fig:valid-precision}
    \label{fig:precision-raw-quotients}
\end{figure}

Figure~\ref{fig:valid-precision} gives a node-level finite-precision check for
the raw Gauss--Jacobi quotients. For a fixed $M=512$ rule, the raw Type-II
quotient loses accuracy first near the endpoint nodes, and the effect is much
stronger in lower precision because the denominator scales as $(t\tau)^2$. The
right panel of Figure~\ref{fig:valid-precision} shows how this local endpoint
loss propagates to the GJ-II derivative error as $M$ changes. The fp8-like curves
are controlled scalar-quantization diagnostics only, not a full mixed-precision
PINN training experiment. Figure~\ref{fig:valid-precision} therefore reinforces
the conclusion from Figure~\ref{fig:valid-M}: when raw GJ-II quotients are used, high precision floating type and  moderate
Gauss--Jacobi $M$ values are important.

\subsection{Cutoff bias--conditioning trade-off}
\label{subsubsec:valid-cutoff-tradeoff}

We next examine how the cutoff $\delta$ balances regularization bias and
endpoint conditioning. For the smooth benchmark \eqref{eq:valid-exp-benchmark},
let $\mathcal D_\delta^{\rm stab}$ denote the regularized operator evaluated
with the expm1/Taylor formulas in \eqref{eq:valid-stable-exp}, and let
$\mathcal D_\delta^{\rm raw}$ denote the same regularized operator evaluated by
direct subtraction. We define the stabilized error and the raw--stable gap
\begin{equation*}
    E_{\rm stab}(\delta)
    =
    \frac{\abs{\mathcal D_\delta^{\rm stab} f(t_0)-\partial_t^\alpha f(t_0)}}
    {\abs{\partial_t^\alpha f(t_0)}},
    \qquad
    E_{\rm gap}(\delta)
    =
    \frac{\abs{\mathcal D_\delta^{\rm raw} f(t_0)
    -\mathcal D_\delta^{\rm stab} f(t_0)}}
    {\abs{\partial_t^\alpha f(t_0)}}
\end{equation*}
as the two diagnostic quantities.

\begin{figure}[htbp]
    \centering
    \includegraphics[width=0.96\textwidth]{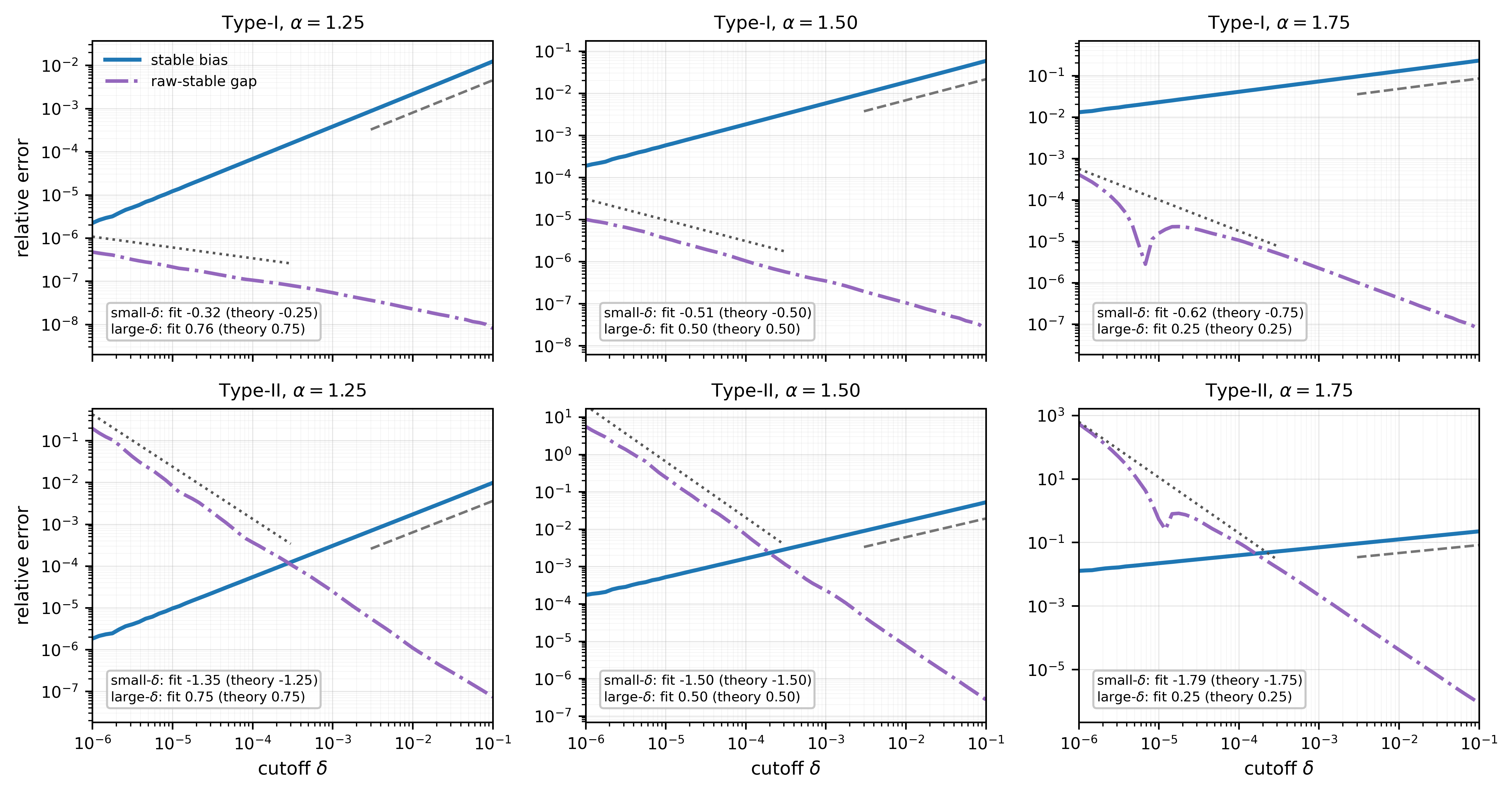}
    \caption{Natural finite-precision cutoff trade-off for the smooth benchmark
    \eqref{eq:valid-exp-benchmark}. The blue curves use stabilized quotient
    evaluation and the purple
    curves show the raw--stable gap. For large $\delta$, the stable error follows
    the regularization-bias exponent $2-\alpha$. For small $\delta$, the
    raw--stable gap follows the endpoint-conditioning exponent $1-\alpha$ for
    Type-I and $-\alpha$ for Type-II, up to finite-precision transition effects.}
    \label{fig:valid-fp-window-tradeoff}
\end{figure}

Remark~\ref{rem:delta-tradeoff} predicts two opposite effects. For large
$\delta$, the stable error is dominated by denominator-regularization bias and
scales like $O(\delta^{2-\alpha})$. For small $\delta$, the raw--stable gap exposes
endpoint conditioning: Type-I has the order $O(\delta^{1-\alpha})$, while Type-II has the stronger value-cancellation order $O(\delta^{-\alpha})$.
Figure~\ref{fig:valid-fp-window-tradeoff} and
Table~\ref{tab:valid-fp-window-slopes} report the corresponding fp32 slope. Thus $\delta$ cannot be chosen too large, because this increases
regularization bias, or too small, because this amplifies finite-precision
endpoint cancellation.

\begin{table}[htbp]
\centering
\caption{Slope fits for the finite-precision cutoff diagnostic in
Figure~\ref{fig:valid-fp-window-tradeoff}. Small-cutoff slopes are fitted from
the raw fp32--stable gap on $10^{-6}\le\delta\le10^{-4}$; large-cutoff slopes are
fitted from the stable error on $3\times10^{-3}\le\delta\le10^{-1}$.}
\label{tab:valid-fp-window-slopes}
\begin{tabular}{c|c|c|cc|cc}
\hline
$\alpha$ & Type & precision & small fit & expected & large fit & expected \\
\hline
1.25 & I  & fp32 & -0.32 & -0.25 & 0.76 & 0.75 \\
1.25 & II & fp32 & -1.35 & -1.25 & 0.75 & 0.75 \\
1.50 & I  & fp32 & -0.51 & -0.50 & 0.50 & 0.50 \\
1.50 & II & fp32 & -1.50 & -1.50 & 0.50 & 0.50 \\
1.75 & I  & fp32 & -0.62 & -0.75 & 0.25 & 0.25 \\
1.75 & II & fp32 & -1.79 & -1.75 & 0.25 & 0.25 \\
\hline
\end{tabular}
\end{table}

The fp32 fits in Table~\ref{tab:valid-fp-window-slopes} are consistent with
Remark~\ref{rem:delta-tradeoff}. The purpose of this diagnostic is therefore practical: it
identifies the window where $\delta$ is small enough to avoid excessive bias but
not so small that raw endpoint cancellation dominates.

\subsection{Monte Carlo and Gauss--Jacobi rate checks}
\label{subsubsec:valid-rate-checks}

\begin{figure}[htbp]
\centering
    \includegraphics[width=0.92\textwidth]{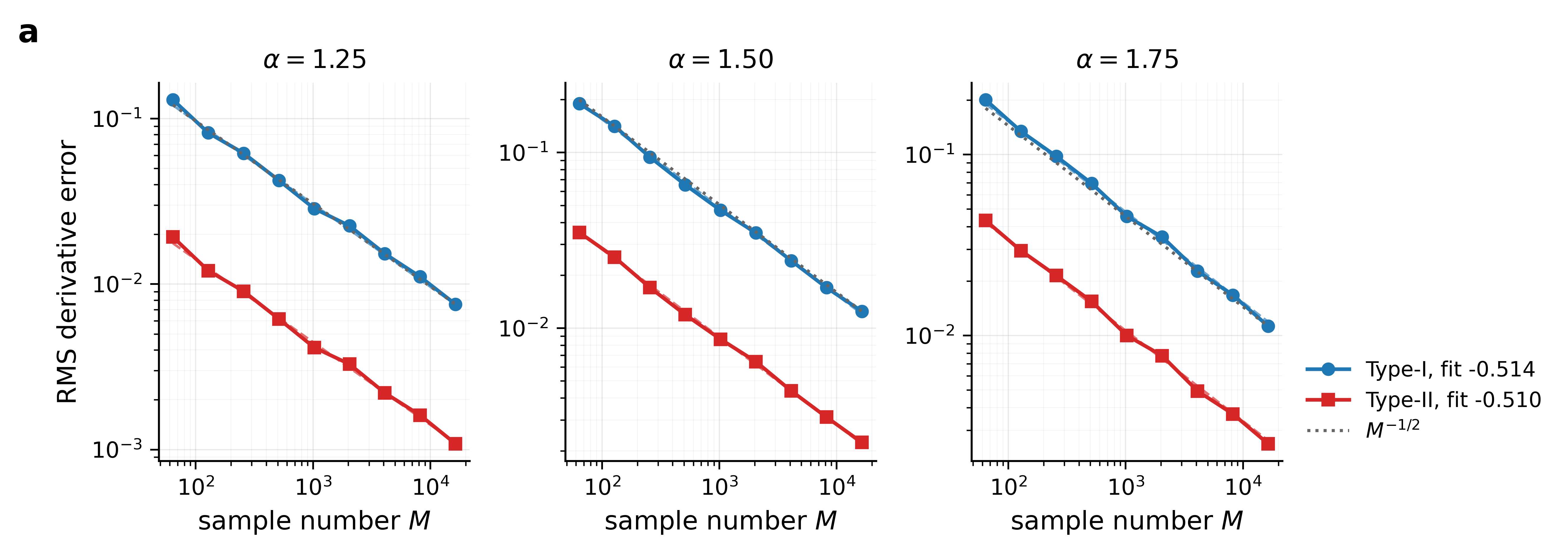}
    \caption{Monte Carlo RMS convergence. Across $\alpha=1.25,1.50,1.75$, the empirical slopes remain close to $-1/2$ for both Type-I and Type-II, consistent with Proposition~\ref{prop:mc-convergence}.}
    \label{fig:valid-mc-rate}
    \label{fig:mc-rate}
\end{figure}

\paragraph*{Monte Carlo}
The Monte Carlo experiment is designed to isolate sampling error and affect from the cutoff. For each
$\alpha\in\{1.25,1.50,1.75\}$ and for both Type-I and Type-II, we compute the
empirical root-mean-square derivative error over 5-times independent samples
using the exact bounded kernels, so the denominator cutoff does not enter this
rate check. Table~\ref{tab:valid-mc-rate-summary} and
Figure~\ref{fig:valid-mc-rate} report fitted log--log slopes close to $-1/2$.
These results support the Monte Carlo convergence estimate in
Proposition~\ref{prop:mc-convergence}.

\begin{table}[htbp]
\centering
\caption{Monte Carlo RMS rate check. Slopes are fitted in log--log scale from empirical RMS derivative errors over repeated independent samples.}
\label{tab:valid-mc-rate-summary}
\begin{tabular}{c|c|c|c}
\hline
$\alpha$ & type & expected slope & fitted slope \\
\hline
1.25 & Type-I  & $-0.50$ & $-0.5005$ \\
1.25 & Type-II & $-0.50$ & $-0.5050$ \\
1.50 & Type-I  & $-0.50$ & $-0.4950$ \\
1.50 & Type-II & $-0.50$ & $-0.4959$ \\
1.75 & Type-I  & $-0.50$ & $-0.5139$ \\
1.75 & Type-II & $-0.50$ & $-0.5103$ \\
\hline
\end{tabular}
\end{table}

\paragraph*{Gauss--Jacobi}
The Gauss--Jacobi experiment checks the two kernel-regularity regimes in
Theorem~\ref{thm:GJ-convergence}: algebraic convergence when
$K_f(t,\cdot)$ and $H_f(t,\cdot)$ have finite smoothness, and spectral
convergence when they are analytic.

\begin{figure}[htbp]
\centering
    \includegraphics[scale=0.6]{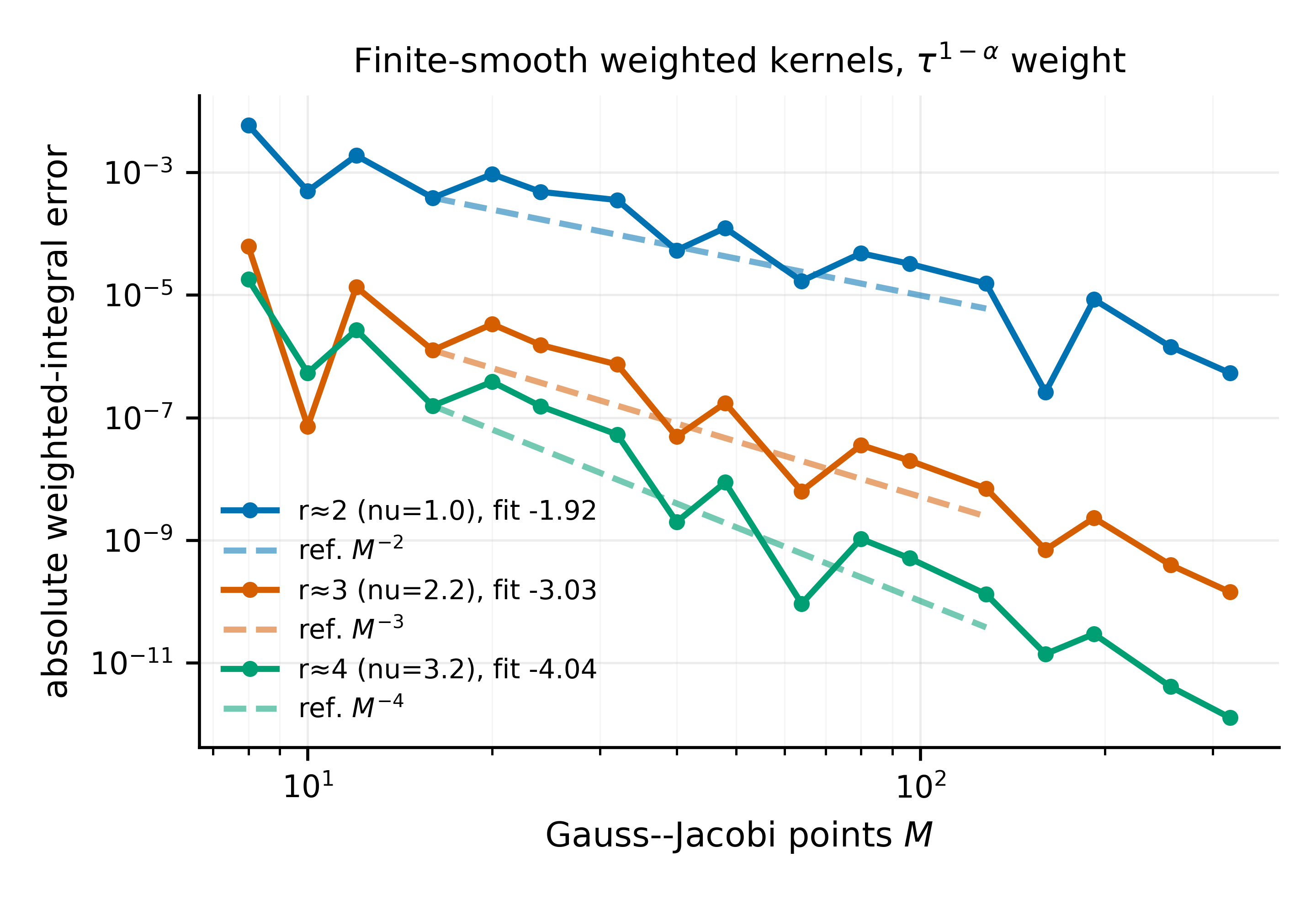}
    \includegraphics[scale=0.6]{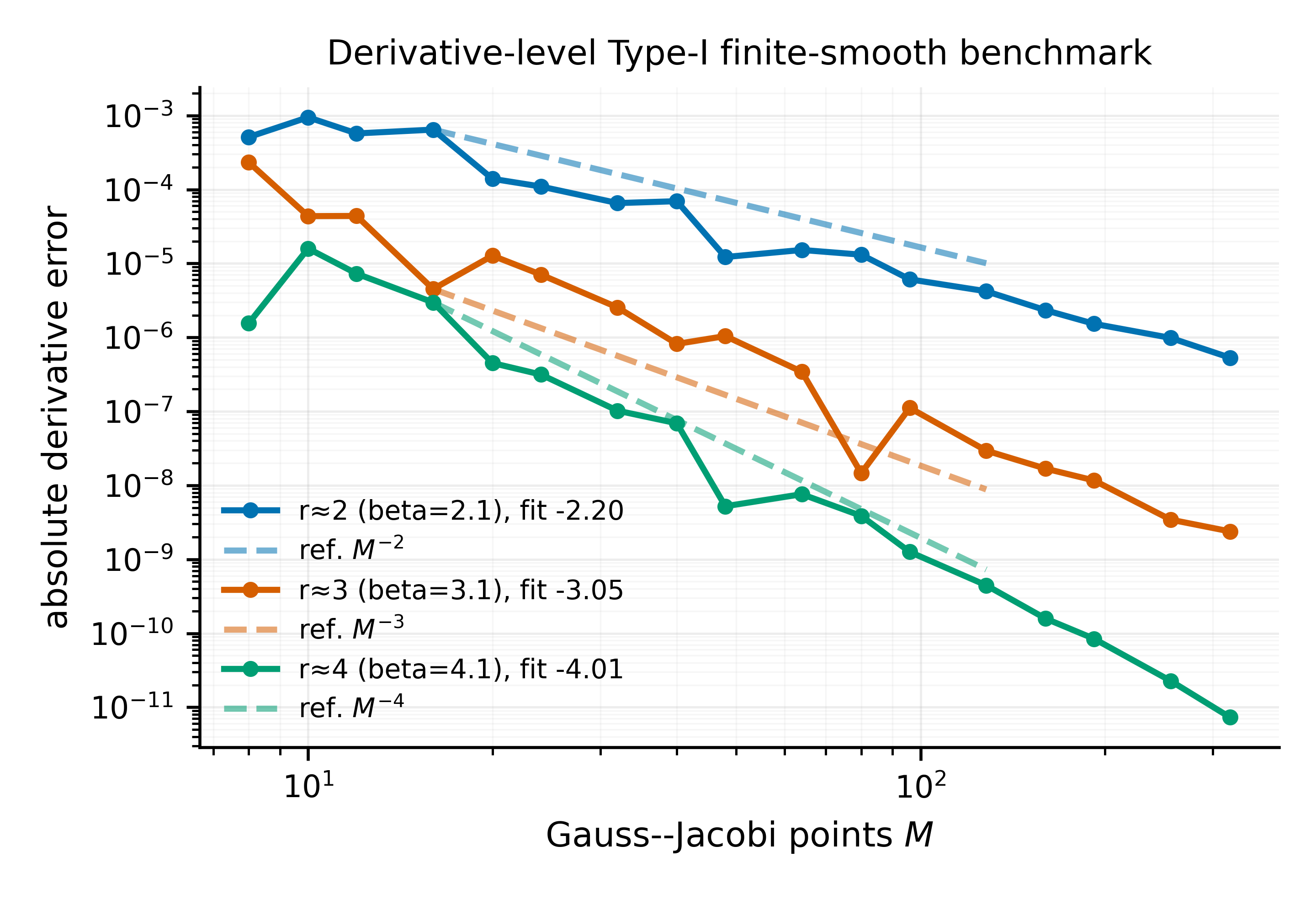}
    \caption{Gauss--Jacobi finite-smoothness convergence at $\alpha=1.5$. Left: the weighted functional \eqref{eq:valid-gj-weighted-functional} with model kernels $\phi_\nu(\tau)=|\tau-\tau_c|^\nu$. Right: the derivative-level Type-I shifted-power benchmark \eqref{eq:valid-shifted-power-benchmark}. The dashed guide lines show the algebraic reference rates $M^{-r}$.}
    \label{fig:valid-gj-mr-rate}
\end{figure}

For the finite-smoothness regime, we first apply the Gauss--Jacobi rule
directly to
\begin{equation}
    I_\alpha[\phi]=\int_0^1 \tau^{1-\alpha}\phi(\tau)\,\d\tau,
    \label{eq:valid-gj-weighted-functional}
\end{equation}
using model kernels
\[
    \phi_\nu(\tau)=|\tau-\tau_c|^\nu,\qquad 0<\tau_c<1 .
\]
This isolates the weighted quadrature functional with controlled finite
regularity. To connect this kernel-level test with an actual transformed
fractional derivative, we also use the shifted-power benchmark
\begin{equation}
    p_{\beta,t_c}(t)=(t-t_c)_+^\beta,\qquad t\ge0 .
    \label{eq:valid-shifted-power-benchmark}
\end{equation}
For $t>t_c$ and $\beta>\alpha-1$, its Caputo derivative is
\begin{equation*}
    \partial_t^\alpha p_{\beta,t_c}(t)
    =\frac{\Gamma(\beta+1)}{\Gamma(\beta+1-\alpha)}(t-t_c)_+^{\beta-\alpha},
    \qquad t>t_c.
\end{equation*}
In the derivative-level test, we set $t=1.5$ and $t_c=0.55$, so the loss of smoothness is mapped to an interior
point of the memory interval. For $\alpha=1.5$, the fitted slopes over
$M=16,\ldots,128$ in Table~\ref{tab:gj-mr-rate} and
Figure~\ref{fig:valid-gj-mr-rate} agree with the algebraic reference slopes.
These results support the finite-smoothness part of
Theorem~\ref{thm:GJ-convergence}.

\begin{table}[htbp]
\centering
\caption{Finite-smoothness Gauss--Jacobi rate check at $\alpha=1.5$. Slopes are fitted in log--log scale over $M=16,\ldots,128$. In the weighted-kernel rows, the tested non-weight integrand is $\phi_\nu(\tau)=|\tau-\tau_c|^\nu$, and $\nu$ controls the finite regularity at the interior point $\tau_c$. The derivative-level rows use the shifted-power benchmark \eqref{eq:valid-shifted-power-benchmark}.}
\label{tab:gj-mr-rate}
\begin{tabular}{c|c|c|c|c}
\hline
Test & parameter & type & reference slope & fitted slope \\
\hline
Weighted kernel & $\nu=1.0$ & $r\simeq2$ & $-2$ & $-1.9180$ \\
Weighted kernel & $\nu=2.2$ & $r\simeq3$ & $-3$ & $-3.0321$ \\
Weighted kernel & $\nu=3.2$ & $r\simeq4$ & $-4$ & $-4.0435$ \\
Type-I derivative & $\beta=2.1$ & $r\simeq2$ & $-2$ & $-2.1953$ \\
Type-I derivative & $\beta=3.1$ & $r\simeq3$ & $-3$ & $-3.0514$ \\
Type-I derivative & $\beta=4.1$ & $r\simeq4$ & $-4$ & $-4.0067$ \\
\hline
\end{tabular}
\end{table}

\begin{figure}[htbp]
    \centering
    \includegraphics[scale=0.6]{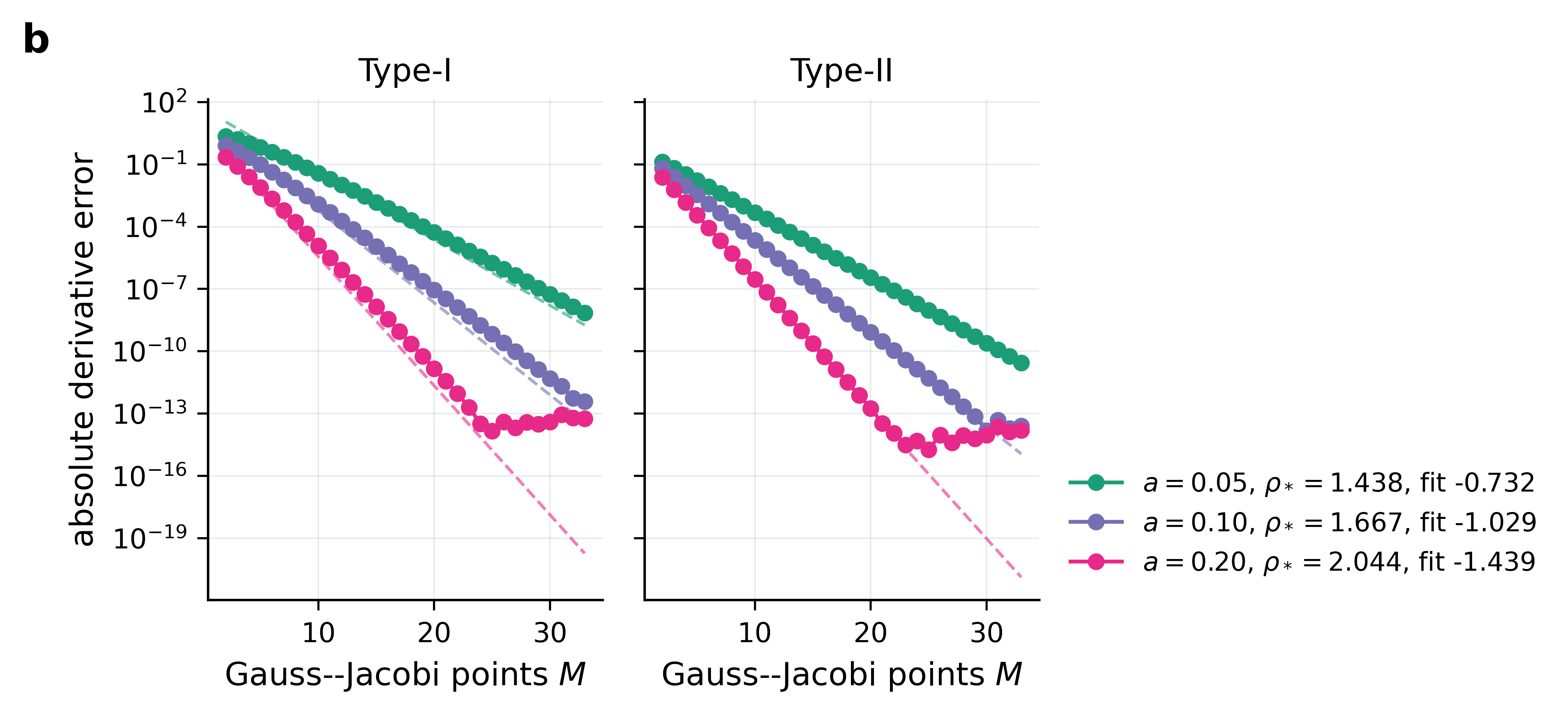}
    \caption{Gauss--Jacobi spectral convergence for analytic rational kernels. The expected pre-saturation exponential slope is determined by the nearest singularity through $-2\log\rho_*$. The fit excludes the machine-precision plateau.}
    \label{fig:valid-gj-rate}
    \label{fig:gj-rate}
\end{figure}

For the analytic regime, the exponential benchmark reaches machine precision too
quickly, so we use the rational benchmark
\begin{equation}
    f_a(t)=\frac{1}{a+t},
    \label{eq:valid-rational-benchmark}
\end{equation}
at $t_0=1.5$. The nearest singularity in the $\tau$-plane is
$\tau_*=1+a/t_0>1$. Under the Bernstein map $x=2\tau-1$, the associated ellipse
parameter is
\begin{equation*}
    \rho_*=x_*+\sqrt{x_*^2-1},\qquad x_*=2\tau_*-1.
\end{equation*}
The expected exponential slope is therefore $s=-2\log\rho_*$. The fitted slopes
in Table~\ref{tab:valid-gj-rho-rate-summary} and
Figure~\ref{fig:valid-gj-rate} match this prediction before machine-precision
saturation, supporting the analytic part of Theorem~\ref{thm:GJ-convergence}. The exeriments are conducted under $\al=1.5$.

\begin{table}[htbp]
\centering
\caption{Analytic Gauss--Jacobi spectral rate check. Slopes are fitted in the pre-saturation regime using the rational benchmark \eqref{eq:valid-rational-benchmark}.}
\label{tab:valid-gj-rho-rate-summary}
\begin{tabular}{c|c|c|c}
\hline
case & type & expected $s$ & fitted $s$ \\
\hline
$a=0.05$, $\rho_*=1.438$ & Type-I  & $-0.7263$ & $-0.7285$ \\
$a=0.05$, $\rho_*=1.438$ & Type-II & $-0.7263$ & $-0.7318$ \\
$a=0.10$, $\rho_*=1.667$ & Type-I  & $-1.0217$ & $-1.0203$ \\
$a=0.10$, $\rho_*=1.667$ & Type-II & $-1.0217$ & $-1.0289$ \\
$a=0.20$, $\rho_*=2.044$ & Type-I  & $-1.4299$ & $-1.4322$ \\
$a=0.20$, $\rho_*=2.044$ & Type-II & $-1.4299$ & $-1.4387$ \\
\hline
\end{tabular}
\end{table}

\subsection{Loss of smoothness}
\label{subsubsec:valid-loss-smoothness}

We try to verify the role of kernel smoothness in
Theorem~\ref{thm:GJ-convergence}. The experiment reuses the shifted-power
benchmark \eqref{eq:valid-shifted-power-benchmark} with $\beta=1.35$,
$t=1.5$, and $\alpha=1.5$. The case $t_c=0$ gives an endpoint function
$t^{1.35}$, which is $C^1$ but not $C^2$ at the left endpoint. The case
$t_c=0.7$ moves the same nonsmooth feature into the memory interval, because
\begin{equation*}
    \tau_c=1-\frac{t_c}{t},
\end{equation*}
so $t_c=0.7$ corresponds to an interior point $\tau_c\approx0.533$. In this
case the induced kernels $K_f(t,\cdot)$ and $H_f(t,\cdot)$ are nonsmooth inside
$[0,1]$, which is more restrictive for Gauss--Jacobi quadrature than an endpoint
loss of smoothness.

\begin{figure}[htbp]
    \centering
    \includegraphics[width=0.92\textwidth]{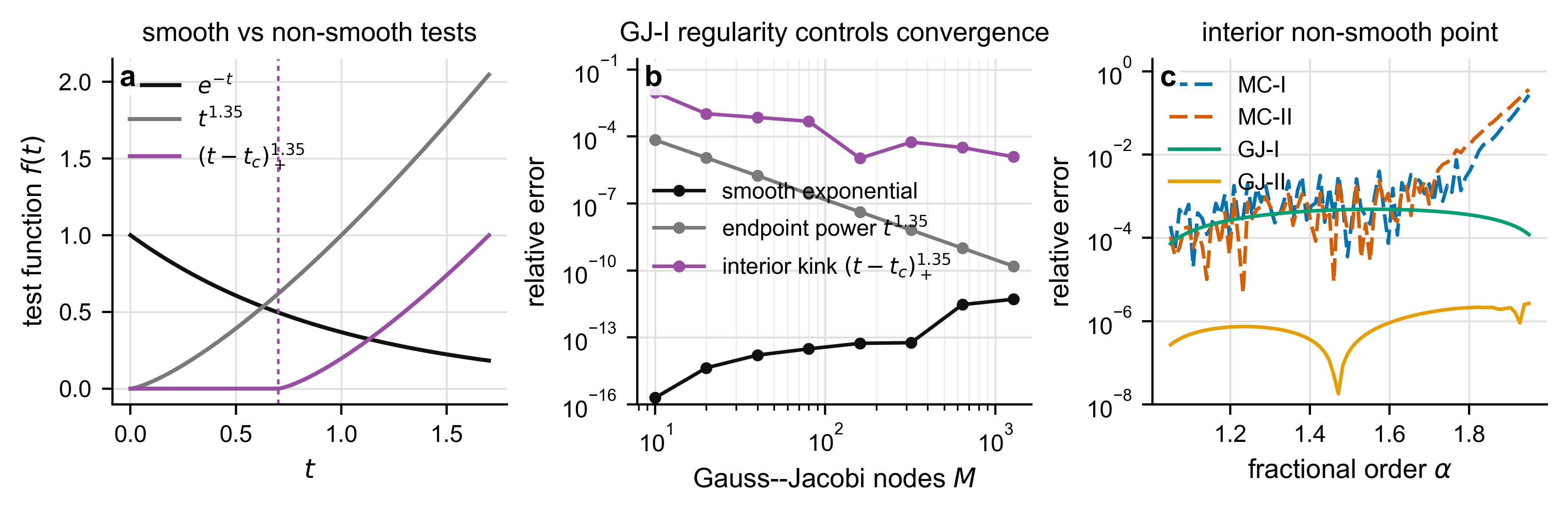}
    \caption{Nonsmooth shifted-power benchmark \eqref{eq:valid-shifted-power-benchmark} with $\beta=1.35$. The endpoint case $t_c=0$ is $C^1$ but not $C^2$ at the left endpoint, while the interior case $t_c=0.7$ places the nonsmooth feature inside the memory interval. The relative-error ordering confirms that loss of smoothness in the induced kernels $K_f(t,\cdot)$ and $H_f(t,\cdot)$ reduces the Gauss--Jacobi convergence rate predicted by Theorem~\ref{thm:GJ-convergence}.}
    \label{fig:valid-nonsmooth}
    \label{fig:nonsmooth-validation}
\end{figure}

Figure~\ref{fig:valid-nonsmooth} confirms this ordering. The smooth exponential
benchmark has the smallest relative error, the endpoint power $t^{1.35}$ gives a
larger but still steadily decreasing error, and the interior-kink case
$(t-t_c)_+^{1.35}$ gives the largest relative error. This behavior is consistent
with Theorem~\ref{thm:GJ-convergence}: once the smoothness of $f$ lowers the
regularity of the induced kernels $K_f(t,\cdot)$ and $H_f(t,\cdot)$, the
Gauss--Jacobi convergence rate is reduced accordingly.

\subsection{Empirical complexity verification}
\label{subsubsec:empirical-complexity-verification}

We verify the complexity estimates in
Section~\ref{subsubsec:computational-complexity} by testing the actual
Type-I/Type-II difference predicted from the computation in
\ref{app:gpu-storage} and \ref{app:gpu-flops}. For graph
storage, \ref{app:gpu-storage} gives the fp64 tensor-data cost of one
ordinary value graph and one time-derivative graph as
\[
    S_u(1)=8(d+c_\sigma LH+1),\qquad
    S_{\partial_tu}(1)=8[2d+(c_\sigma+2)LH+1].
\]
For each collocation point, Type-I stores $M$ shifted time-derivative graphs
and two endpoint time-derivative graphs, whereas Type-II stores $M$ shifted
ordinary value graphs and the same two endpoint time-derivative graphs. Hence
\begin{align*}
    S_I^{\rm graph}
    &=N(M+2)S_{\partial_tu}(1),\\
    S_{II}^{\rm graph}
    &=NM S_u(1)+2N S_{\partial_tu}(1).
\end{align*}
Subtracting the two formulas removes the common endpoint contribution and gives
\[
    \Delta S
    :=S_I^{\rm graph}-S_{II}^{\rm graph}
    =NM\{S_{\partial_tu}(1)-S_u(1)\}
    =8NM(d+2LH).
\]
Thus the predicted graph-storage saving has order
$O(NM(d+LH))$.

The FLOP prediction is assembled in the same way from
\ref{app:gpu-flops}. With
$A_{\rm mac}=Hd+(L-1)H^2+H$, the one-graph training-backward counts are
\[
    F_u^{\rm bwd}(1)=4A_{\rm mac}+3LH,\qquad
    F_{\partial_tu}^{\rm bwd}(1)=8A_{\rm mac}+8LH.
\]
Therefore
\begin{align*}
    F_I^{\rm bwd}
    &=N(M+2)F_{\partial_tu}^{\rm bwd}(1),\\
    F_{II}^{\rm bwd}
    &=NM F_u^{\rm bwd}(1)+2N F_{\partial_tu}^{\rm bwd}(1),
\end{align*}
and the predicted training-backward FLOP saving is
\begin{align*}
    \Delta F
    &:=F_I^{\rm bwd}-F_{II}^{\rm bwd}\\
    &=NM\{F_{\partial_tu}^{\rm bwd}(1)-F_u^{\rm bwd}(1)\}\\
    &=NM(4A_{\rm mac}+5LH).
\end{align*}
Equivalently, the leading order is $O(NM(dH+LH^2))$; this predicts first-order
dependence on $N$, $M$, $d$, and $L$, and second-order dependence on $H$ when
the hidden-to-hidden term dominates.

The empirical benchmark constructs the generalized Gauss--Jacobi Type-I and
Type-II residual graphs directly, using the settings in
\ref{app:empirical-complexity-settings}. For each sweep, only one of
$N$, $M$, $d$, $L$, and $H$ is varied while the others are fixed. The plotted
differences are
\[
    \Delta S
    =S_I^{\rm graph}-S_{II}^{\rm graph},
    \qquad
    \Delta F
    =F_I^{\rm bwd}-F_{II}^{\rm bwd}.
\]
The storage experiment records the peak CUDA allocation during residual-graph
construction, and the FLOP experiment profiles only the training backward pass
after the residual graph has been built.

\begin{figure}[htbp]
    \centering
    \includegraphics[width=\textwidth]{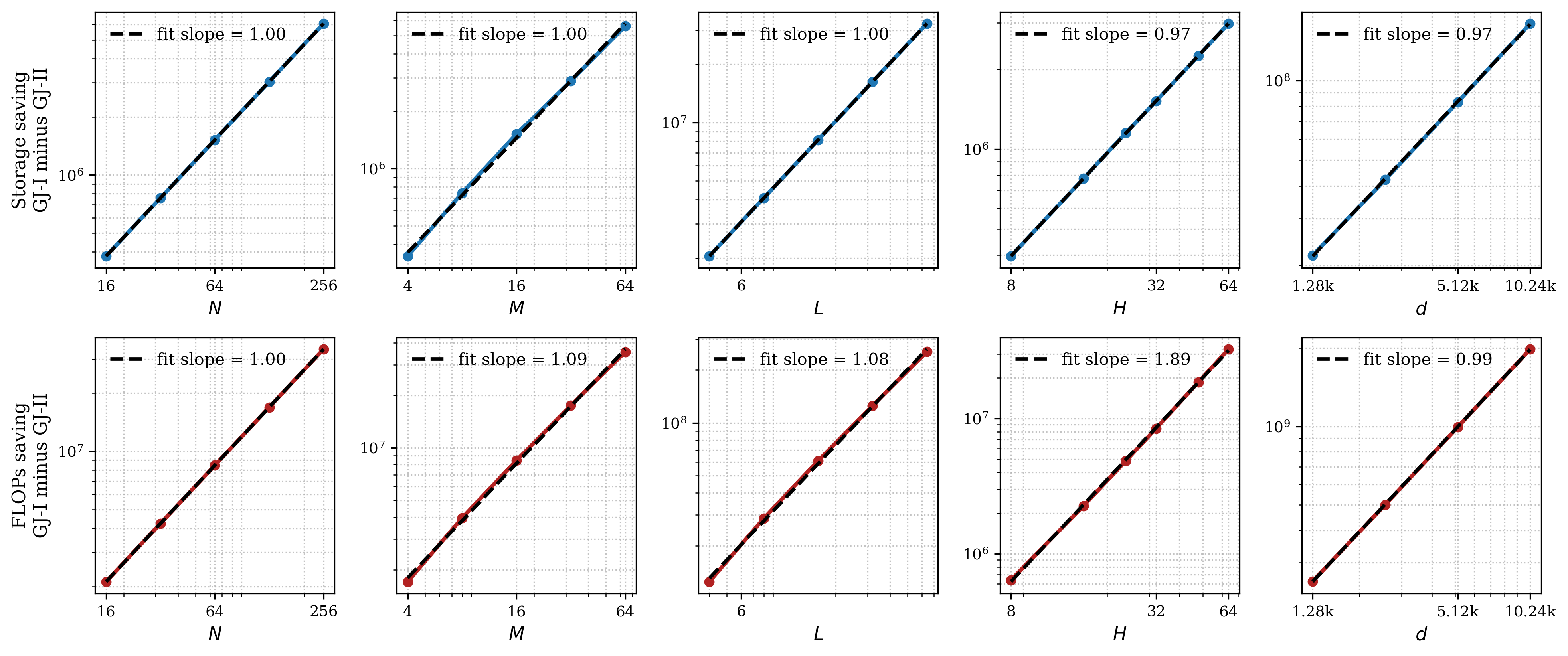}
    \caption{Empirical storage and training-backward FLOP savings for
    Gauss--Jacobi Type-I versus Type-II residuals. The top row reports measured
    peak CUDA graph-storage saving, and the bottom row reports profiler FLOP
    saving. Columns correspond to sweeps over
    $N$, $M$, $L$, $H$, and $d$. Each panel plots the median measured saving from three repeated measurements
    $($Type-I minus Type-II$)$ and its log--log fit.}
    \label{fig:empirical-complexity}
\end{figure}

Figure~\ref{fig:empirical-complexity} and Table~\ref{tab:empirical-complexity-slopes} confirms the predicted orders of these
two differences. The storage row follows
$\Delta S=O(NMd)+O(NMLH)$, with fitted log--log slopes close to the expected
linear dependence in the $N$, $M$, $d$, $L$, and $H$ sweeps. The FLOP row
matches the prediction for $\Delta F$: the measured saving is nearly linear in
$N$, $M$, $d$, and $L$, and grows almost quadratically in $H$. These results verify the order of the Type-I/Type-II
difference, not merely the absolute runtime of a particular implementation.

\begin{table}[htbp]
\centering
\footnotesize
\setlength{\tabcolsep}{4pt}
\renewcommand{\arraystretch}{1.12}
\caption{Expected and fitted log--log slopes for the complexity savings in
Figure~\ref{fig:empirical-complexity}.}
\label{tab:empirical-complexity-slopes}
\begin{tabular}{c|cc|cc}
\hline
& \multicolumn{2}{c|}{$\Delta S$}
& \multicolumn{2}{c}{$\Delta F$} \\
Sweep parameter & Expected & Fitted & Expected & Fitted \\
\hline
$N$ & $1$ & $1.00$ & $1$ & $1.00$ \\
$M$ & $1$ & $1.00$ & $1$ & $1.09$ \\
$L$ & $1$ & $1.00$ & $1$ & $1.08$ \\
$H$ & $1$ & $0.97$ & $2$ & $1.89$ \\
$d$ & $1$ & $0.97$ & $1$ & $0.99$ \\
\hline
\end{tabular}
\end{table}

\subsection{Implications for the PINN experiments}
\label{subsubsec:valid-implications-pinn}

This appendix connects the representation theory with the
implementation behavior of the four estimators used in
Section~\ref{sec:experiments}. The derivative-level tests isolate four factors:
the fractional order $\alpha$, the quadrature size $M$, finite-precision
endpoint cancellation, and the smoothness inherited by the transformed kernels
$K_f(t,\cdot)$ and $H_f(t,\cdot)$. These conclusions concern the numerical
evaluation of the fractional derivative; the full PINN error also includes
network approximation and optimization effects.

For Monte Carlo quadrature, Section~\ref{subsubsec:valid-cutoff-tradeoff} shows
that the cutoff becomes important as $\alpha\to2$, because the
$\mathrm{Beta}(2-\alpha,1)$ samples concentrate near $\tau=0$ and the raw
difference quotients become more sensitive to endpoint cancellation. Once the
cutoff is fixed, the sampling part follows the usual RMS rate
$O(M^{-1/2})$, as verified in Section~\ref{subsubsec:valid-rate-checks}. Thus
Monte Carlo is robust but usually needs larger $M$ than Gauss--Jacobi, together
with a carefully chosen cutoff near $\tau=0$.

For Gauss--Jacobi quadrature, the decisive condition is the regularity of the
induced kernels. Section~\ref{subsubsec:valid-rate-checks} verifies that smooth
kernels give substantially faster convergence than Monte Carlo: finitely smooth
kernels lead to algebraic rates, while analytic kernels exhibit the spectral
$\rho^{-2M}$ behavior in Theorem~\ref{thm:GJ-convergence}. Conversely,
Section~\ref{subsubsec:valid-loss-smoothness} shows that loss of smoothness in
$K_f(t,\cdot)$ or $H_f(t,\cdot)$ directly degrades the observed rate.

Finally, the Type-I/Type-II distinction is mainly computational but has an
implementation caveat. Section~\ref{subsubsec:empirical-complexity-verification}
shows that Type-II is much cheaper because it removes the shifted
time-derivative graphs from the $M$ quadrature nodes. However, the Type-II raw
quotient is also more sensitive to the floating-point type, $\alpha$, and $M$:
near $\alpha=2$ or for very large $M$, small endpoint nodes can make roundoff
visible. In the PINN experiments below, Type-II is therefore used for its
computational advantage, but the quadrature size, cutoff, and fp64/stable
endpoint evaluation must be chosen cautiously.

\section{Derivation of the graph-storage estimates}
\label{app:gpu-storage}

This appendix gives the storage calculation used in
Section~\ref{subsubsec:computational-complexity}. The count deliberately
isolates the fp64 tensor data retained by reverse-mode automatic
differentiation; allocator overhead, bookkeeping data, and other
implementation-dependent constants are not included.

Let $z=(t,\x)\in\mathbb R^d$ and consider a bias-free fully connected network with
$L$ hidden layers of width $H$ and scalar output. We write
\begin{align*}
    a^0 &= z,\\
    h^\ell &= W^\ell a^{\ell-1},
    \qquad
    a^\ell=\sigma(h^\ell),
    \qquad \ell=1,\ldots,L,\\
    u(z;\theta) &= W^{L+1}a^L.
\end{align*}
Thus $a^0\in\mathbb R^d$, $a^\ell,h^\ell\in\mathbb R^H$ for
$\ell=1,\ldots,L$, and $u(z;\theta)\in\mathbb R$.

\paragraph{Ordinary value graph.}
If the scalar value $u(z;\theta)$ enters the residual, then the training
backward pass must keep the usual forward graph. At the tensor-data level, this
requires the input, the hidden activations, and the scalar output. Some
implementations retain only activations, whereas others retain both activations
and pre-activations, or equivalent data for evaluating $\sigma'$. We represent
this by a constant $c_\sigma\in[1,2]$. Hence the fp64 storage for one ordinary
value graph is
\begin{equation}\label{eqn:app-value-storage}
    S_u(1)=8\bigl(d+c_\sigma LH+1\bigr)
    \quad\text{bytes}.
\end{equation}

\paragraph{Time-derivative graph}
When the residual contains $\partial_t u(z;\theta)$, training requires the loss
backward pass to compute
\[
    \nabla_\theta \partial_t u(z;\theta).
\]
Therefore the graph that produces the time derivative must itself be
differentiable. The reverse sweep for the input gradient can be written as
\begin{align*}
    r^L &= (W^{L+1})^T,\\
    s^\ell &= r^\ell\odot\sigma'(h^\ell),
    \qquad \ell=L,\ldots,1,\\
    r^{\ell-1} &= (W^\ell)^T s^\ell,
    \qquad \ell=L,\ldots,1.
\end{align*}
The vector $r^0\in\mathbb R^d$ is $\nabla_z u(z;\theta)$, and
$\partial_t u(z;\theta)$ is its time component. In addition to the ordinary
forward graph, the derivative graph must retain the sensitivity vectors
\[
    r^1,\ldots,r^L,\qquad s^1,\ldots,s^L,\qquad r^0.
\]
These account for $2LH+d$ additional fp64 scalars. Therefore
\begin{equation}\label{eqn:app-dt-storage}
    S_{\partial_tu}(1)
    =8\bigl[d+c_\sigma LH+1+d+2LH\bigr]
    =8\bigl[2d+(c_\sigma+2)LH+1\bigr]
    \quad\text{bytes}.
\end{equation}
Equations~\eqref{eqn:app-value-storage} and \eqref{eqn:app-dt-storage} are the
only graph-storage ingredients needed in the main text. Section~\ref{subsubsec:computational-complexity}
combines them with the number of shifted and endpoint quadrature evaluations to
obtain the Type-I and Type-II storage counts.

\section{Derivation of the training-backward FLOP estimates}
\label{app:gpu-flops}

This appendix derives the backward FLOP count used in
Section~\ref{subsubsec:computational-complexity}. We count one multiply-add as
two floating-point operations. The network is bias-free, so no bias additions
or bias-gradient reductions occur; memory traffic, allocator costs, and
bookkeeping overhead are omitted.

\paragraph{Dense-layer backward.}
For a dense layer in batch form,
\[
    Y=XW,\qquad
    X\in\mathbb R^{B\times n_{\rm in}},
    \quad
    W\in\mathbb R^{n_{\rm in}\times n_{\rm out}},
\]
the forward matrix product contains $B n_{\rm in}n_{\rm out}$ MACs, or
$2B n_{\rm in}n_{\rm out}$ FLOPs. In the ordinary backward pass, the two
dominant matrix products are
\begin{align*}
    \overline W &= X^T\overline Y,\\
    \overline X &= \overline Y W^T .
\end{align*}
Each costs $2B n_{\rm in}n_{\rm out}$ FLOPs. Therefore the dense part of the
ordinary backward pass costs
\begin{equation}\label{eqn:app-dense-ordinary-bwd}
    4B n_{\rm in}n_{\rm out}
\end{equation}
FLOPs.

\paragraph{Ordinary value graph}
For the network architecture $d\to H\to\cdots\to H\to1$, the MAC count of one
plain forward pass is
\begin{equation}\label{eqn:app-A-mac}
    A_{\rm mac}=Hd+(L-1)H^2+H.
\end{equation}
Applying \eqref{eqn:app-dense-ordinary-bwd} to all dense layers gives the dense
backward contribution $4A_{\rm mac}$. With the $\tanh$ activation, the
activation backward may be evaluated from the retained activation $a$ by
forming
\[
    q=1-a^2,\qquad \overline h=\overline a\,q,
\]
which costs three scalar operations per hidden unit. Therefore, for one
ordinary value graph,
\begin{equation}\label{eqn:app-value-bwd-flops}
    F_u^{\rm bwd}(1)=4A_{\rm mac}+3LH.
\end{equation}

\paragraph{Time-derivative graph.}
A residual term containing $\partial_t u(z;\theta)$ is more expensive because
the training backward pass must compute
\[
    \nabla_\theta \partial_t u(z;\theta),
\]
that is, it must differentiate through the reverse-mode sweep that produced the
input derivative. With
\[
    s^\ell=r^\ell\odot\sigma'(h^\ell),
    \qquad
    r^{\ell-1}=(W^\ell)^Ts^\ell,
\]
the backward pass through $r^{\ell-1}=(W^\ell)^Ts^\ell$ produces the same
leading matrix-product cost as an ordinary dense-layer backward. In addition,
because $s^\ell$ depends on the original pre-activation $h^\ell$, the adjoint
generated through this dependence must also be propagated through the original
forward dense layer. Hence the dense-layer contribution is doubled, giving
$8A_{\rm mac}$. For $\tanh$, differentiating the sensitivity recursion
contributes eight scalar operations per hidden unit in the explicit count used
in the main text. Thus
\begin{equation}\label{eqn:app-dt-bwd-flops}
    F_{\partial_tu}^{\rm bwd}(1)=8A_{\rm mac}+8LH.
\end{equation}
Equations~\eqref{eqn:app-value-bwd-flops} and \eqref{eqn:app-dt-bwd-flops} are
the one-graph backward counts used in Section~\ref{subsubsec:computational-complexity}.
The assembly over the $M$ shifted quadrature nodes and the two endpoint terms is
carried out there, not repeated here.

\section{Empirical complexity benchmark settings}
\label{app:empirical-complexity-settings}

Table~\ref{tab:empirical-complexity-settings} lists the benchmark settings used
for Figure~\ref{fig:empirical-complexity}. The implementation constructs the same
generalized Gauss--Jacobi Type-I and Type-II residuals for all sweeps; there is
no change of algorithmic backend between parameters.

\begin{table}[htbp]
\centering
\footnotesize
\setlength{\tabcolsep}{4pt}
\renewcommand{\arraystretch}{1.18}
\caption{Parameter settings for the empirical storage and backward-FLOP
complexity verification.}
\label{tab:empirical-complexity-settings}
\begin{tabular}{p{0.24\textwidth}|p{0.68\textwidth}}
\hline
Item & Setting \\
\hline
Baseline parameters & $N=64$, $M=16$, $d=2$, $L=3$, $H=32$,
$\alpha=1.5$. In each sweep, one of $N,M,d,L,H$ is varied and the remaining
parameters stay at these baseline values. \\
\hline
Sweep values & $N\in\{16,32,64,128,256\}$,
$M\in\{4,8,16,32,64\}$,
$d\in\{1280,2560,5120,10240\}$,
$L\in\{4,8,16,32,64\}$,
$H\in\{8,16,24,32,48,64\}$. \\
\hline
Precision and device & CUDA execution with fp64 tensor data. For storage, the
fp64 choice is the precision used in the graph-storage formulas in
Section~\ref{subsubsec:computational-complexity}; using fp32 would halve the
byte prefactor but leave the parameter scaling unchanged. \\
\hline
Storage recording & CUDA memory is recorded with PyTorch's
\texttt{torch.cuda.memory\_allocated()} before constructing the fractional
residual graph, after the model, collocation points, and quadrature data have
already been allocated. We then call
\texttt{torch.cuda.reset\_peak\_memory\_stats()} and record the peak graph
construction footprint using \texttt{torch.cuda.max\_memory\_allocated()}. The
plotted storage saving is the Type-I minus Type-II peak increment,
$\Delta S=S_I^{\rm graph}-S_{II}^{\rm graph}$. \\
\hline
Timing recording & Graph-construction and backward wall-clock times are recorded
with \texttt{time.perf\_counter()}. For CUDA runs, we call
\texttt{torch.cuda.synchronize()} before stopping the timer, so asynchronous GPU
kernels are included in the measured elapsed time. \\
\hline
Warm-up and aggregation & An unrecorded CUDA/autograd warm-up pair is run
before the recorded sweep. Storage uses three repeated measurements and plots
median savings; FLOPs use three repeated measurements per parameter value. \\
\hline
\end{tabular}
\end{table}

\bibliographystyle{elsarticle-num}
\bibliography{ref}

\end{document}